\documentclass[final,onefignum,onetabnum]{siamart220329}

\usepackage{lipsum}
\usepackage{amsfonts, amssymb, amsmath}
\usepackage{graphicx}
\usepackage{epstopdf}
\usepackage{algorithmic}
\usepackage{bm}
\usepackage[caption=false]{subfig}
\usepackage{xcolor}
\usepackage{tikz}
\usetikzlibrary{arrows,3d,calc,intersections,decorations.markings}

\definecolor{darkred}{rgb}{0.75, 0.0, 0.0}
\definecolor{brightred}{rgb}{0.95, 0.0, 0.0}

\newcommand{\drawarcdelta}[4]{
	\draw[blue, stealth-] ($#1+(#2:#4)$) arc[start angle=#2, delta angle=#3, radius=#4];
	\draw[blue, -stealth] ($#1+(#2:#4)$) arc[start angle=#2, delta angle=#3, radius=#4];
}

\newcommand{\drawlabeledarcdelta}[6]{
	\drawarcdelta{#1}{#2}{#3}{#4}
	\node at ($#1+(#2+#3/2:#6)$) {#5};
}

\ifpdf
\DeclareGraphicsExtensions{.eps,.pdf,.png,.jpg}
\else
\DeclareGraphicsExtensions{.eps}
\fi

\newcommand{\argdot}{\makebox[1.5ex]{\textbf{\(\cdot\)}}}%
\DeclareMathOperator\supp{supp}

\newsiamremark{remark}{Remark}
\newsiamremark{hypothesis}{Hypothesis}
\crefname{hypothesis}{Hypothesis}{Hypotheses}
\newsiamthm{claim}{Claim}

\headers{Imaging of atmospheric dispersion processes with DIAL}{R. Lung, N. Polydorides}

\title{Imaging of atmospheric dispersion processes with Differential Absorption Lidar
	\thanks{\textbf{Funding:} RL was supported by the Maxwell Advanced Technology Scholarship}}

\author{Robert Lung\thanks{School of Engineering, University of Edinburgh, Edinburgh, UK 
		(\email{robert.lung@ed.ac.uk}, \email{n.polydorides@ed.ac.uk}).}
	\and Nick Polydorides\footnotemark[2]}

\usepackage{amsopn}

\makeatletter
\newcommand*{\addFileDependency}[1]{
	\typeout{(#1)}
	\@addtofilelist{#1}
	\IfFileExists{#1}{}{\typeout{No file #1.}}
}
\makeatother

\crefalias{formulation}{equation}
\crefname{formulation}{formulation}{formulations}
\Crefname{formulation}{Formulation}{Formulations}
\creflabelformat{condition}{#2\textup{(#1)}#3} 

\begin{document}

\maketitle
\begin{abstract}
  We consider the inverse problem of fitting atmospheric dispersion parameters based on time-resolved back-scattered differential absorption Lidar (DIAL) measurements. The obvious advantage of light-based remote sensing modalities is their extended spatial range which makes them less sensitive to strictly local perturbations/modelling errors or the distance to the plume source. In contrast to other state-of-the-art DIAL methods, we do not make a single scattering assumption but rather propose a new type modality which includes the collection of multiply scattered photons from wider/multiple fields-of-view and argue that this data, paired with a time dependent radiative transfer model, is beneficial for the reconstruction of certain image features. The resulting inverse problem is solved by means of a semi-parametric approach in which the image is reduced to a small number of dispersion related parameters and high-dimensional but computationally convenient nuisance component. This not only allows us to effectively avoid a high-dimensional inverse problem but simultaneously provides a natural regularisation mechanism along with parameters which are directly related to the dispersion model. These can be associated with meaningful physical units while spatial concentration profiles can be obtained by means of forward evaluation of the dispersion process.
\end{abstract}


\section{Introduction}
Optical measurements have been used successfully in atmospheric science for different purposes \cite{Weitkamp:2005}. For example, differential absorption lidar (DIAL) systems that operate on two adjustable wavelengths, have been used for determining the concentration of gases in the atmosphere or the release rate of gaseous plumes \cite{Robinson:2011, Innocenti:2017}. Classical Lidar methods typically rely on a single scattering assumption for the recorded photons and possibly some multiple scattering approximations/corrections based on simplifying assumptions \cite{Hogan:2008, HoganBattaglia:2008}. Wider fields-of-view (FOVs) have been shown to contain information beyond what is contained in single scattering \cite{Hutt:1994, CalahanEtAl:2005} but are clearly only useful if the optical forward model is chosen appropriately so that this data can be considered as an information signal as opposed to a nuisance. More recently, the authors of \cite{LevisSchechnerDavis:2015, HolodovskyEtAl:2016, LevisSchechnerDavis:2017} demonstrated the superiority of accurate forward models for optical measurements by considering the full radiative transfer equation in applications involving cloud and aerosol imaging. The more realistic assumptions on the light propagation come at a very substantial increase in computing time but can be remedied by using advanced algorithmic methods and low dimensional models for the quantities of interest \cite{BalLangmoreMarzouk:2013, LangmoreDavisBal:2013}. Using the full Radiative Transfer Equation (RTE) as an optical model in the context of DIAL data, without enforcing strong assumptions on the scattering functions, is not straight-forward and, to the best of our knowledge, remains an open problem.

For many practical purposes it is prohibitively expensive to acquire accurate measurements of a 3-dimensional region which means that regularisation based on prior assumptions is necessary in order to obtain meaningful images. To address this issue we propose a regularisation of the optical inverse problem based on a dispersion process which enforces smoothness while at the same time reducing the concentration profile to a low-dimensional set of parameters which benefits the computationally challenging reconstruction process when complex optical forward models are used. The choice of parameterisation is not only suitable for localisation of gas plumes at the time of measurement but also makes future tracking of gas plumes much simpler as the recovered parameter values can be directly used in the associated dispersion process. Optical remote sensing measurements will by design be spatially averaged to some degree and thus less sensitive to  less smooth image features, e.g. local variability caused by turbulence in the wind field, that may alter the shape of the plume \cite{Stull1988}. Our method chooses a parameterisation based on dispersion related quantities that explicitly disregards such non-smooth features and is therefore particularly suitable when wide FOVs are used and the averaging effect becomes more pronounced.

The paper is organised as follows. In \cref{sec:description} we give a formal description of our problem and derive a coupled formulation of dispersion and optical models. In \cref{sec:uniqueness} we show that the extended optical forward model preserves the uniqueness of absorption field under essentially the same assumptions of DIAL and provide insights into when the additional data benefits the image reconstruction. In \cref{sec:computing} we give a description of our algorithm which is capable of efficiently utilising prior knowledge in the form of a suitable dispersion process. Numerical results to validate our findings are provided in section \cref{sec:experiments}.

\section{Problem description \& motivation}\label{sec:description}
We start with a general form of the problem at hand. Assume that the plume which we are trying to image is given by a smooth function, i.e. \(u \in \mathcal{V} \subseteq C^{\infty}(S,\mathbb{R})\) for some \(S \subseteq \mathbb{R}^d\) where the dimension \(d\) depends on, for example, whether there is a temporal component or the plume is stationary. Since the used measurements are virtually always going to be corrupted by a significant amount of noise, we assume that we have access to data in the form of realisations of a random element \(\bm{Y}\) whose distribution \(\mathbb{P}(\bm{Y} \in \argdot  \mid u)\) given the image \(u\) is known and given a realisation \(\bm{Y}(\omega) = \bm{y}\) of \(\bm{Y}\) we may construct a suitable loss functional \(u \mapsto \mathsf{Loss}(u \mid \bm{y})\), e.g. a least squares penalty \(\|\mathbb{E}(\bm{Y} \mid u) - \bm{y}\|^2\) or, as in our case, a negative profile log-likelihood, whose precise description we give in \cref{sec:computing}, of the image \(u\) for the observed data \(\bm{y}\), and find an estimator \(\hat{u}\) for the image by solving the optimisation problem
\begin{align}\label{eq:mle}
    \hat{u} = \arg\min_{u \in  \mathcal{V}} \mathsf{Loss}(u \mid \bm{y}).
\end{align}
There are two competing objectives when dealing with estimators as the one given in \Cref{eq:mle}. On one hand we would like the set \(\mathcal{V}\) to be large, i.e. cover a broad range of images \(u \in C^{\infty}(S, \mathbb{R})\) in order to accommodate as much flexibility in the image reconstruction as possible. On the other hand, any finite parameterisation/discretisation of a large set \(\mathcal{V}\) will require a large amount of parameters resulting in severely ill-posed and, for complex forward operators \(\mathbb{E}(\bm{Y} \mid u)\), computationally intractable problems. The standard approach to such problems typically results in a regularised version of \Cref{eq:mle} taking the form of a high dimensional problem
\begin{align}\label{eq:soft_regularised_reconstruction}
    \hat{u} = \arg\min_{u \in  \mathcal{V}} \mathsf{Loss}(u \mid \bm{y}) + \gamma \mathsf{Reg}(u)
\end{align}
for a \(\mathcal{V}\) parameterised via local basis functions, a regularisation penalty \(\mathsf{Reg}(u)\), e.g. the norm of a discrete differential operator such as \(\|\Delta u\|^2\) with \(\gamma > 0\). Alternative approaches, such as modeling the response as a smooth function (e.g. with low-degree piece-wise polynomials \cite{HarrisEtAl:2016}) have a very similar effect. The smoothness enforcing regularisation term remedies some of the ill-posedness but the problem in \Cref{eq:soft_regularised_reconstruction} is often even more computationally demanding than that of \Cref{eq:mle}. Furthermore, when the data are extremely noisy it may become necessary to increase \(\gamma\) to an extent where too much structure gets lost for it to be of actual value for the intended application.
\paragraph{Dispersion based regularisation terms} Conceptually, our approach can be derived from \Cref{eq:soft_regularised_reconstruction} by first replacing a generic regularisation term such as \(\|\Delta u\|^2\) with \(\|\mathsf{D}_{\theta}(u)\|^2\) where the operator \(\mathsf{D}_{\theta}\) is a differential operator that relates dispersion quantities \(\theta \in \Theta \subseteq \mathbb{R}^d\) to gas concentration over time. This results in an informed regularisation penalty which preserves the structure of a plume even for large regularisation parameters \(\gamma\) and is similar to the parameterisation used in \cite{BalLangmoreMarzouk:2013}. If we take the limit \(\gamma \to \infty\) and assume that for each \(\theta \in \Theta\) there is a unique \(u[\theta]\) so that \(\mathsf{D}_{\theta}(u) = 0\) then \Cref{eq:soft_regularised_reconstruction} is equivalent to solving
\begin{align}\label[formulation]{eq:dsp_regularised_reconstruction}
    \hat{\theta} = \arg\min_{\theta \in  \Theta} \mathsf{Loss}(u[\theta] \mid \bm{y}) \quad \mathrm{and} \quad \hat{u} = u[\hat{\theta}].
\end{align}
The map \(\theta \mapsto u[\theta]\) is, unlike parameterisations based on local (voxel-based) functions, non-linear but it is often possible to select a parameter space \(\Theta\) of rather low dimension and still obtain meaningful images which makes \Cref{eq:dsp_regularised_reconstruction} a low-dimensional problem and, given that we have taken \(\gamma \to \infty\), a good robustness towards noise can be expected for well behaving dispersion operators, i.e. when the parameter \(\theta\) restricts the flexibility of potential solutions. If the latter is not the case, then essentially the same issues that arise in \Cref{eq:dsp_regularised_reconstruction} will be present in any method for plume tracking and localisation based on that dispersion model regardless of how concentration was obtained. The main issue with an approach such as \Cref{eq:dsp_regularised_reconstruction} is placing constraints that are too restrictive and result in unrealistic images with large biases.

\subsection{Dispersion models}\label{sec:dispersion_model}
In this section we introduce a family of widely used gas dispersion models, which under suitable boundary conditions and domains have a closed form solution. Following the developments in \cite{Stockie:2011} we consider the advection-diffusion operator given by 
\begin{align}\label{eq:general_plume_equation}
    \mathsf{D}_{\theta}(u) = \frac{\partial}{\partial t} u + \nabla \cdot (\eta u) - Q + \frac{1}{2}\nabla \cdot (\kappa \nabla u)
\end{align}
with \(\theta = (\eta,Q,\kappa)\). In \Cref{eq:general_plume_equation} we use \(\eta\) to model the drift, \(Q\) is a source term, \(\kappa\) is a diagonal matrix with diffusion coefficients and \(u\) is the gas concentration. We can assume a point source located at \(q = (q_1,q_2,q_3)^\top\) which means that \(Q\) takes the form
\begin{align}\label{eq:point_source}
    Q(x,t) = \rho_Q \cdot \delta(x_1-q_1)\delta(x_2-q_2)\delta(x_3-q_3)\delta(t)
\end{align}
where \(\rho_Q > 0\) models the amount of released gas and \(\delta\) are single variable delta distributions. It is also assumed that \(\kappa=\kappa(t)\) is a function of time. This is similar to assuming dependence on downwind distance \(\|x-q\|\) since particles move downwind with time according to the wind speed. We shall also assume position-independent wind although we may want to allow time-dependence, i.e. \(\eta=\eta(t)\) but we do not assume that diffusion in the wind direction is much smaller than advection and can thus be ignored, which is assumed in order to derive the steady-state Gaussian plume models \cite{Stockie:2011, SeinfeldPandis:2016}, as this would prevent us from using the dispersion model in situations with negligible wind. Despite their simplicity and potential inaccuracy in complex environments, modified variants of Gaussian plumes are still used for regulatory purposes and, due to the excessive computational cost of many alternative approaches, are the preferred choice for time critical applications and long-term average loads \cite{LeelossyEtAl:2014}. Nevertheless, we chose a somewhat more general approach to model steady-state plumes via superpositions of instantaneous releases, which is also commonly used approach for modelling of dispersion in practice \cite{LeelossyEtAl:2014, LepicardHelingMaderich:2004, ZhouEtAl:2003, GhannamElFadel:2013, LevyEtAl:2002, PrueksakornEtAl:2012}. In that situation, assuming an infinite domain without boundary we can solve \Cref{eq:general_plume_equation}
\begin{align}\label{eq:general_plume_solution}
    u[\theta](x,t) = \prod_{j=1}^3 \frac{1}{\sqrt{2\pi \int_0^t \kappa_{jj}(s)ds}}\exp\left(-\frac{\left(x_j-\left[q_j + \int_0^t w_j(s) ds\right]\right)^2}{2\int_0^t \kappa_{jj}(s)ds}\right)
\end{align}
Assuming that the diffusion is isotropic, i.e. \(\kappa\) is a scalar function \(\kappa_0\) multiplied by the identity matrix, and setting \(h(t) = \int_0^t \kappa_{0}(s)ds\) as well as writing \(m_q(t) = q + \int_0^t w(s)ds\) we can simplify \Cref{eq:general_plume_solution} to
\begin{align}\label{eq:simpler_plume_solution}
    u[\theta](x,t) = \left(\frac{1}{\sqrt{2\pi h(t)}}\right)^3 \exp\left( -\frac{\|x-m_q(t)\|^2}{2h(t)}\right).
\end{align}
In situations where the boundary is not negligible, e.g. when the plume is close to a reflecting or absorbing flat ground at \(x_3=0\), similar expressions to \Cref{eq:general_plume_solution} and \Cref{eq:simpler_plume_solution} can be obtained.
\paragraph{Steady-state models and super-positions} 
The Gaussian puff model from the previous section considers an instantaneous release and as such it is not suitable for modelling longer or ongoing releases. Instead of using the solution presented in the previous paragraph we may consider a continuous release as the integrated superposition of instantaneous releases at different times. In fact, by allowing a sufficient amount of independent spatially distributed sources we could essentially approximate any smooth function with puffs such as in \Cref{eq:simpler_plume_solution}, regardless of boundary conditions or even dispersion dynamics. This means that the complexity of the dispersion controls the flexibility of the images and thus the amount of regularisation in \Cref{eq:dsp_regularised_reconstruction}.

In this paper we assume that \(h(t)\) as well as \(m_q(t)\) are piece-wise linear and that we can use a small number of independent sources to model the structure of the dispersion. In other words, we consider \(u = \sum_{j=1}^K u_i\) where each \(u_i\) has the form of \Cref{eq:simpler_plume_solution}. For the sake of consistency with our later developments we shall also include a constant, i.e. \((x,t)\)-independent, term \(u_0 = u_0[\theta]\) to model a homogeneous concentration field. If we assume that each source releases gas continuously in an interval \([0, T]\) for some \(T \in (0,\infty]\) and use multiple puffs to approximate the continuous release then, for \(t \leq T\), we end up with a dispersion model of the form
\begin{align}\label{eq:kernel_parameterisation}
    u_0[\theta]+\int_{0}^{t} u[\theta](x,s) ds \approx \sum_{j=1}^N w_j \phi\left( \frac{\|x-m_j\|_2}{h_j}\right)
\end{align}
for some suitable, possibly \(t\)-dependent, values for \(w_j, h_j\) and \(m_j\) and the squared exponential kernel \(\phi(x) = \exp(-x^2/2)\). When the plume is in steady-state we have \(t = \infty\) in \Cref{eq:kernel_parameterisation} and there is no dependence on time. Note that if we ignore the difference between a Gaussian kernel and a voxel indicator, then \Cref{eq:kernel_parameterisation} differs from a voxel based parameterisation only in that the basis functions are not fixed in place or size and instead are allowed to vary in these quantities. The non-uniqueness caused by allowing this additional variability is tackled by applying dispersion based constraints to their weights, sizes and position so that the associated inverse problem doesn't become unsolvable.
\paragraph{Entropy of atmospheric dispersion} It is worth noticing that the above model is most useful and accurate in scenarios where the plume is observed over longer periods of time which obviously isn't possible when the release is instant or inhomogeneous. The above model will therefore inevitably have some error and not necessarily represent the true dispersion process accurately. Instead of viewing the above dispersion model as an approximate truth to the real-world particle transport one might also consider reconstructing certain aspects \(F_0(u), \dots, F_k(u)\) of the image \(u\), where \(F_j\) can in principle evaluate any real-valued feature of the image \(u\), regardless of any knowledge about underlying phenomena responsible for the motion of airborne particles. 
If \(F_j\) evaluate localised integrals over equal patches of the domain, i.e. voxels, then \(F_j(u)\) correspond to local averages. Most standard methods seek to reconstruct these aspects of the image and assume that \(u\) is constant within each voxel which coincides with the maximum entropy distribution given the constraints \(F_j(u)\) while the number of constraints used corresponds to the image resolution in the usual sense. 

The concept of filling in the missing information by means of finding the distribution the with maximal entropy that satisfies a set of constraints can be generalised in a straight-forward way to functionals other than voxel averages and represents in a fairly strong sense the least biased and optimal way of choosing distributions given partial or incomplete knowledge \cite{Jaynes:1957}. Similarly we can interpret the number of real-value constraints as a more general form of image resolution (informally one may think of this as coefficients in a different basis even though the \(F_j\) can be non-linear) where the number of parameters needed to represent the function to a certain degree of accuracy heavily depends on the selected features. As an alternative to local averages over voxels we can consider \(F_j\) that correspond to the amount of gas, position and width of the plume which can be expressed through moments of the function \(u\). In domains that are mostly unconstrained, i.e. open space or uncluttered environments, the maximum entropy solution given the first three moment constraints will be (approximately) Gaussian and take a form similar to \Cref{eq:simpler_plume_solution}. As such this can be considered as our best guess for the gas distribution given only the amount, the location and a (homogeneous) dispersion rate independently of any knowledge about the underlying atmospheric transport. 

The entropy based derivation of Equations \cref{eq:general_plume_solution} or \cref{eq:simpler_plume_solution} is arguably more general in that it doesn't assume the existence of a hypothetical average, i.e. ergodicity. The empirical observations that support such assumptions merely suggest that the chosen parameterisation is sensible and turbulence induced deviations behave much like irreversible, entropy increasing operations resulting in errors that in a way resemble random noise. Although it may seem unnecessary at this point to derive the dispersion model in two different ways, the above arguments will be helpful in understanding the reasoning behind our treatment of missing information for the optical transport problem where averaging becomes essentially meaningless (much like for instantaneous gas-releases). In particular, the mechanism responsible for errors due to regularisation of turbulence is essentially identical to what will be used to handle missing information regarding optical parameters due to our inability to measure light to a degree sufficient for reconstruction of all involved parameters.

\subsection{Radiative transfer}
As of now we haven't discussed how the optical measurements are influenced by a plume. In principle the forward model is indifferent to any other optical imaging or tomography setup where it's given by the time-dependent Radiative Transfer Equation (RTE), as studied in \cite{GkioulekasLevinZickler2016}, or some simplification thereof like the steady-state RTE \cite{LevisSchechnerDavis:2015,HolodovskyEtAl:2016,LevisSchechnerDavis:2017} or a single-scatter approximation. The time-dependent RTE describes the intensity of light at location \(x \in X\) travelling in direction \(v \in \mathbb{S}^{2}\) at any time \(t \geq 0\) and is given by the integro-differential equation
\begin{align}\label{eq:time_dependent_RTE}
    \left(\frac{\partial}{\partial t} + v \cdot \nabla_{x} + \sigma_{a+s}(x) \right) H(x,v,t) = \sigma_s(x) \int_{\mathbb{S}^{2}} H(x,v',t) f_p(x, v\cdot v') dv' 
\end{align}
for \((x,v,t) \in X \times \mathbb{S}^{2} \times [0, \infty)\), where the spatial domain of interest \(X \subseteq \mathbb{R}^3\) is a bounded open set with sufficiently regular boundary \(\partial X\). In \Cref{eq:time_dependent_RTE} we have normalised the time component such that the speed of light is equal to \(1\). The functions \(\sigma_a, \sigma_s : X \to [0, \infty) \) are the absorption and scattering coefficients of the material, \(\sigma_{a+s}:= \sigma_a + \sigma_s\) and \(f_p : X \times \mathbb{S}^{2}\times \mathbb{S}^{2} \to [0, \infty)\) is called phase function and describes the distribution of angles for scattered light. Since \(f_p\) is a probability density on the sphere it satisfies
\begin{align*}
    \int_{\mathbb{S}^{2}} f_p(x, v\cdot v') dv' = 1 \qquad \forall x \in X, v \in \mathbb{S}^{2}.
\end{align*}
In order to state the boundary conditions for \Cref{eq:time_dependent_RTE} in the spatial boundary \(\partial X\) we introduce for any \(x \in \partial X\) the hemispheres
\begin{align*}
    \partial V_{(-)}(x) = \left\{v \in \mathbb{S}^{2}: v \cdot \Vec{n}(x) \leq 0\right\} \qquad
    \partial V_{(+)}(x) = \left\{v \in \mathbb{S}^{2}: v \cdot \Vec{n}(x) > 0\right\}
\end{align*}
where \(\Vec{n}(x)\) is the unit outer normal at \(x \in \partial X\). In other words, at a point \(x \in \partial X\), \(\partial V_{(-)}(x)\) and \(\partial V_{(+)}(x)\) correspond to the inward and outward pointing directions respectively. The impact of time-homogeneous ambient illumination will be discussed later and we can assume that at time \(t=0\) there is no light inside the domain and consider a known source at the boundary given by a function defined on \(X \times \mathbb{S}^2 \times \mathbb{R}\) which is inward pointing, i.e. \(\supp g \subseteq \{(x,v,t) \in \partial X \times \mathbb{S}^2 \times (0, \infty): v \in \partial V_{(-)}(x)\}\). In the case of an instantaneous Lidar pulse at time \(t_0\) released in direction \(v_0\) from a source located at \(x_D\) we may think of \(g\) as smooth and compactly supported such that \(g(x,v,t) \approx \delta(x-x_D)\delta(v-v_0)\delta(t-t_0)\). If we ignore reflections at the boundary this translates to
\begin{subequations}
\begin{align}
    H(x,v,0) &= 0  &\forall x \in X, v \in \mathbb{S}^{2}\label{eq:temporal_boundary_interior_RTE}\\
    H(x,v,t) &= g(x,v,t)  &\forall x \in \partial X, v \in \partial{V}_{(-)}(x) \label{eq:temporal_boundary_boundary_RTE}
\end{align}
\end{subequations}
More generally, reflections at the boundary could be expressed through
\begin{align}\label{eq:spatial_boundary_conditions_RTE}
    H(x,v,t) &= \int_{\partial V_{(+)}(x)} H(x,v',t) f_s(x, v \to v') \lvert \Vec{n}(x) \cdot v' \rvert dv'
\end{align}
for \(x \in \partial X, v \in \partial V_{(-)}(x)\) and \(t>0\). Although it plays a similar role as \(f_p\) we don't require \(f_s\) to be a probability density bur rather that it integrates to at most \(1\). For \Cref{eq:spatial_boundary_conditions_RTE}, also known as the \emph{rendering equation}, can be understood in such a way that it describes the outgoing intensity as the reflected intensity of accumulated incoming light at the boundary represented by the integral. The factor \(\lvert \Vec{n}(x) \cdot v' \rvert\) accounts for the fact that a piece of solid angle from a light source that illuminates a surface will spread over a smaller area when it's perpendicular to the surface than when it comes in almost parallel. Henceforth we will assume that \(f_s = 0\) (the primary purpose of which is to keep the notation simple), i.e. that light is absorbed by the boundary. This isn't always necessary and we will explain how our main theoretical results can be extended to more general non-absorbing boundaries. More details on the derivation and intuition behind the RTE can be found in \cite{Arvo93} and the supplementary material of \cite{GkioulekasLevinZickler2016}. 

\paragraph{Integral representation and Neumann series}
\Cref{eq:time_dependent_RTE} can be recast into an equivalent equation that only has integrals which will make it easier to analyse, at least for our purposes, and will also be used for solving the RTE. This is also known as the \emph{volume rendering equation} due to its similarity with \Cref{eq:spatial_boundary_conditions_RTE}.
\begin{align}
    H(x,v,t) &= \mathcal{I}[g](x,v,t) + \mathcal{K}[H](x,v,t)
\end{align}
where the operators \(\mathcal{I}\) and \(\mathcal{K}\) are defined by
\begin{subequations}
\begin{align}
    \mathcal{I}[g](x,v,t) &:= e^{-\int_{0}^{\ell_{v}(x)}\sigma(x-sv)ds}g\left(x-\ell_{v}(x)v,v,t-\ell_{v}(x)\right)\label{eq:ballistic_operator}\\
    \mathcal{K}[H](x,v,t) &:= \int_{0}^{t} e^{-\int_{0}^{s}\sigma(x-rv)dr}\label{eq:scattering_operator} 1_{(0, \ell_v(x))}(s)\\
    &\phantom{:=} \int_{\mathbb{S}^{2}}H\left(x-sv,v',t-s\right) \sigma_{s}(x-sv) f_p(x-sv, v\cdot v')dv' ds.\nonumber
\end{align}
\end{subequations}
with \(\ell_{v}(x):=\inf\{t > 0:x-tv \not\in X\}\) and \(1_A(x)\) the indicator function on the set \(A\). Note that \(\mathcal{I}\) takes arguments that are functions defined on the boundary while \(\mathcal{K}\) acts on functions that depend on interior points of the domain. Under suitable conditions, which are true for most materials that are of interest to us, we can invert \(\mathrm{Id}-\mathcal{K}\) and obtain the following \emph{Neumann series} for \(t > 0\) 
\begin{align}\label{eq:Neumann_series}
    H(x,v,t) & = (\mathrm{Id}-\mathcal{K})^{-1}\mathcal{I}[g](x,v,t) = \sum_{j = 0}^{\infty} \mathcal{K}^{j}\mathcal{I}[g](x,v,t) = \sum_{j=0}^{\infty} H_j(x,v,t)
\end{align}
which decomposes the contribution into the orders of scattering. Note that the expression for the series in case of non-trivial reflecting boundaries, i.e. when \(f_s \neq 0\), is similar but ends up being slightly more involved than the one from \cite{Bal09} due to the mixed nature of boundary conditions. 

\paragraph{Albedo operators and their Schwartz kernels}
Since the measurements are typically taken at the boundary of the domain it is essential for their description to consider the so called albedo operators. If we define 
\begin{align*}
    \partial D_{(\pm)} := \{(x,v) \in \partial X \times \mathbb{S}^2 : v \in \partial V_{(\pm)}(x)\}
\end{align*}
then these are given by the restriction of the solution from \Cref{eq:time_dependent_RTE} to the spatial boundary in the outward direction \(D_{(+)} \times (0,\infty)\), i.e.
\begin{align}
    \mathcal{A}[g] = (\mathrm{Id}-\mathcal{K})^{-1}\mathcal{I}[g]\big\vert_{\partial D_{(+)} \times (0,\infty)} = H\big\vert_{\partial D_{(+)} \times (0,\infty)}\label{eq:albedo_operator}
\end{align}
For more precise trace results and well-definedness statements for \(\mathcal{A}\) we refer to \cite{ChoulliStefanov1999, Bal:2010}. We will also consider the decomposition into the Neumann series components, i.e.
\begin{align}
    \mathcal{A}_j[g] = H_j\big\vert_{\partial D_{(+)} \times (0,\infty)} 
\end{align}
so that the albedo operator can be expressed equivalently as
\begin{align*}
    \mathcal{A}[g] = \sum_{j=0}^{\infty} H_j\big\vert_{\partial D_{(+)}} = \sum_{j=0}^{\infty} \mathcal{A}_j[g]
\end{align*}
which is very similar to the developments in \cite{Bal09b}. A more detailed treatment of the transport equation is given in \cite{DautrayLions:2012} to which we refer for further details. Finally we want to consider the Schwartz kernel \(\mathfrak{a}(x,v,y,w,s)\) of \(\mathcal{A}\) which has the property that for smooth test functions \(g:D_{(-)}\times(0,\infty)\to \mathbb{R}\) and \(b:D_{(+)}\times(0,\infty) \to \mathbb{R}\) we have
\begin{align}
\begin{split}
    &\int_{\partial D_{(+)}\times (0,\infty)}\mathcal{A}[g](x,v,t)b(x,v,t)dxdvdt \\
    &~= \int_{\partial D_{(+)} \times \partial D_{(-)}\times (0,\infty)^2}\mathfrak{a}(x,v,y,v',s)  g(y,w,t-s)b(x,v,t) dxdv dydv' dtds \label{eq:albedo_kernel}.
\end{split}
\end{align}
Schwartz kernels \(\mathfrak{a}_j\) for the operators \(\mathcal{A}_j\) may be defined in a similar way as in \Cref{eq:albedo_kernel} and will be used in the following developments.
\paragraph{Stochastic process for photon transport} A notable similarity between the solutions of Equations \cref{eq:time_dependent_RTE} and \cref{eq:general_plume_equation} is that, after suitable normalisation, both may be interpreted as probability densities. The solution \(u(\argdot,t)\) of the dispersion model at time \(t\) describes the probability distribution of airborne particles within the domain at time while \(H(\argdot,\argdot,t)\) relates to the position and direction of photons. In a similar way the Schwartz kernels from \Cref{eq:albedo_kernel} can be thought of as densities conditional on a photon being detected. 
More generally these equations describe a distribution over trajectories of the associated particles. In the case of atmospheric dispersion we don't need to concern ourselves with the associated Brownian paths because the way a particle reaches a certain position doesn't influence its optical properties and therefore has no effect on our measurements. In the case of (differential) absorption, which is described by Beer-Lambert's law, the intensity is a function of the path integral and therefore dependent on a photon's trajectory. 

It shouldn't come as much of a surprise that in general we won't be able to reconstruct all optical parameters from two Lidar measurements unless they are significantly constrained. We argue that in a way this is a similar problem to the issues related to turbulence. The wind velocity at every point in time and space alongside the rate of gas release would in theory be sufficient to reconstruct atmospheric transport to an arbitrary degree of detail \cite{Stockie:2011}. Of course such measurements are unrealistic and to the best of our knowledge there is no reliable model for predicting turbulent transport based on crude meteorological data such as wind velocity, pressure, humidity, temperature etc. measured at a single point inside the domain. Similarly unrealistic optical data, i.e. noise free measurements covering the whole domain at every point in time, would be sufficient for perfect reconstruction of (time dependent) images as well. Consequently our image parameterisation is chosen such that the signal strength in the optical measurement is likely to be sufficient in order to reconstruct point estimates of the parameters necessary for the evaluation of an approximation to the atmospheric forward problem in a lower resolution. 

Unlike in the case of atmospheric dispersion, there is no intuitive notion of averages since the source is pulsed and the scattering/absorption parameters depend not only on the number of particles but also their type (e.g. Mie-scattering for spheres \cite{Hulst:1957}) which typically won't change over time. Nevertheless we can compensate for missing information by means of entropy maximisation but before we are ready to explain the details of our method we must develop a notion of entropy on the optical paths as well as a better understanding of the information content within our measurement which will allow us to find suitable parameters/constraints that will be the target of our inverse problem.

\section{DIAL type measurements with wide FOVs}\label{sec:measurement_theory}
In the presence of a gas plume the perturbations in the scattering and attenuation coefficients of the ambient medium that describe the optical properties of the image should be proportional to the gas concentration \(u\). In other words, there are wavelength-dependent constants \(C_{a} > 0\) and \(C_{s} > 0\) such that \(\forall x \in X\)
\begin{align} \label{eq:full_knowledge}
\begin{split}
   \sigma_{a,\mathrm{ambient}}(x) + C_{a} u(x) &= \sigma_{a}(x) \\
   \sigma_{s,\mathrm{ambient}}(x) + C_{s} u(x) &= \sigma_{s}(x)
\end{split}
\end{align}
where \(u(x)\) is the gas concentration at \(x \in X\). If the plume is a mixture of unknown gases, then having access to \(C_s\) and \(C_a\) is not realistic. In the case of DIAL, where we seek to map the concentration of a particular species of gas, it is only assumed that \(C_a > 0\) is known for two carefully selected wavelengths, i.e. we know the absorption behaviour of the gas of interest, \(\sigma_{a,\mathrm{ambient}}(x)\) is known at these wavelengths while all other quantities in \Cref{eq:time_dependent_RTE} are assumed unknown but \emph{identical} in both instances. If we denote the chosen wavelengths by ''\(\mathrm{on}\)'' and ''\(\mathrm{off}\)'' then this can be expressed as
\begin{align} \label{eq:equal_scattering}
\begin{split}
    \sigma_{s(\mathrm{on})}(x) &= \sigma_{s(\mathrm{off})}(x) \\
    f_{p(\mathrm{on})}(x, \argdot) &= f_{p(\mathrm{off})}(x, \argdot)    
\end{split}
\end{align}
and the same holds true for any possibly occurring surface reflections due to \(f_s\) in \Cref{eq:spatial_boundary_conditions_RTE}. Further we assume that the absorption difference \(\sigma_{a(\mathrm{on})} - \sigma_{a(\mathrm{off})}\) for these two two wavelengths satisfies 
\begin{align}\label{eq:differential_absorption_parameterisation}
    C_{\mathrm{ambient}}[\theta] + C_{\mathrm{DIAL}}u[\theta](x) &= \sigma_{a(\mathrm{on})}(x) - \sigma_{a(\mathrm{off})}(x).
\end{align}
for some \emph{known} \(C_{\mathrm{DIAL}} = C_{\mathrm{on}}-C_{\mathrm{off}} > 0\) and possibly unknown and \(\theta\)-dependent constant \(C_{\mathrm{ambient}} = C_{\mathrm{DIAL}}u_0\geq 0\) which corresponds to a not necessarily known amount of homogeneously distributed gas in the ambient atmosphere. In cases where the concentration is much bigger than what's typically present in the atmosphere we may assume \(C_{\mathrm{ambient}} \approx 0\). Note that although \(u\) in general depends on time, such as in \Cref{eq:simpler_plume_solution} and \cref{eq:kernel_parameterisation}, the quantity in \Cref{eq:differential_absorption_parameterisation} decidedly doesn't depend on \(t\) because the time scales relevant for \Cref{eq:time_dependent_RTE} and our optical measurements are many orders of magnitude smaller than those on which any meaningful change in the gas concentration \(u\) might occur. 
\subsection{Wide FOVs \& angularly averaged measurements}
It is worth noticing that our approach in principle generalises to scenarios where the ambient gas concentration is inhomogeneous but known, which might be the case if an area is being actively monitored and accurate measurements that have been averaged over long time periods are available, while for problems where a plume is released into an unknown inhomogeneous atmosphere our dispersion based prior knowledge cannot account for most of the variability in the absorption difference \(\sigma_{a(\mathrm{on})} - \sigma_{a(\mathrm{off})}\) and our method becomes unsuitable. These assumptions are quite different than assuming that the ambient quantities in \Cref{eq:full_knowledge} are fully known and all unknown perturbations are caused by the plume since we only assume that the difference in absorption is caused by the plume in a known way while all other quantities are free, although we will restrict them to some extent later on, as long as they are equal on both wavelengths. For the single scattering approximation and narrow FOVs this of course results in the usual DIAL method. When more complicated measurements are to be used there is no closed form solution for the gas concentration anymore and, depending on the measurement setup and assumptions on the optical parameters it might not be possible to reconstruct all unknown quantities of \Cref{eq:time_dependent_RTE} and careful analysis of the required assumptions and expected errors must be performed.

\begin{figure}
\centering
    \scalebox{0.5}{\input{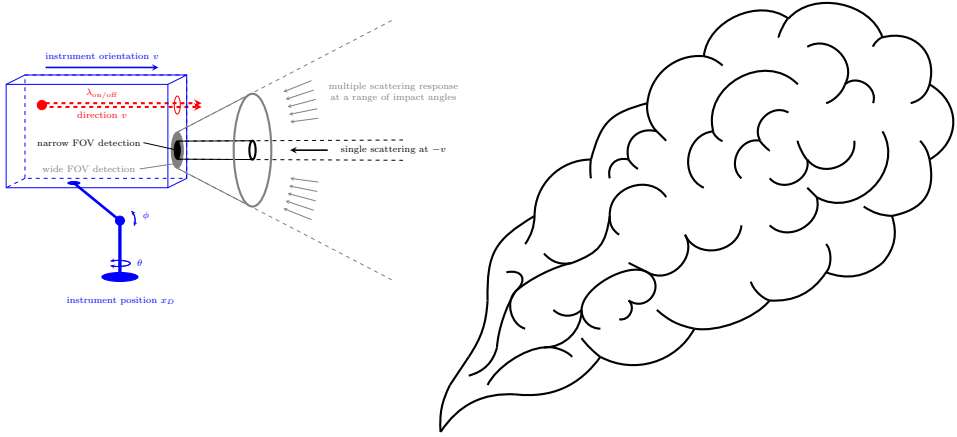}}
    \caption{Laser pulses of two different wavelengths, \(\lambda_{\mathrm{on/off}}\) (red), are released at the source, located at \(x_D\), in direction \(v\). A detector, also located at \(x_D,\) with a narrow narrow FOV (black) captures the single scattering response incident at \(-v\) whereas a wide FOV captures light from a range of directions which have scattered multiple times. After a measurement in direction \(v\) has been taken the instrument is re-oriented by adjusting \(\theta\) or \(\phi\) (blue). This procedure is carried out for azimuthal angle and polar angles at a fixed instrument location \(x_D\).}
    \label{fig:measurement_modality}
\end{figure}
In practice we measure light that reaches a detector that is orders of magnitude smaller than \(\partial X\) and can be essentially thought of as a point measurement at the boundary (see \cref{fig:measurement_modality}). However, the solutions for \Cref{eq:time_dependent_RTE} as well as the traces are developed in such a way that, in order to have a mathematically consistent description of the measurement, we can a priori only restrict the albedo operator to some (small) region around a fixed point and think of a point measurement as a limit when the area of that region tends to zero. For an instantaneous pulse released at time \(t_0\) in direction \(v_0\) from a source and detector located at \(x_D\) that takes an angularly averaged measurement at time \(s\) within a field of view modelled by a continuous function \(b_\mathfrak{m}\) that means taking \(g \to \delta(x-x_D)\delta(v-v_0)\delta(t-t_0)\) and \(b \to \delta(x-x_D) \delta(t-s) b_{\mathfrak{m}}(x,v)\) in the situation of \Cref{eq:albedo_kernel}. The following result ensures that the resulting limit makes sense.
\begin{lemma}[Angularly averaged distribution kernels]\label{lem:path_continuity}
    Consider the Schwartz kernels \(\mathfrak{a}_k\) for the summands of the albedo operator as in \Cref{eq:albedo_kernel} and assume that the optical parameters \(\sigma_{a}, \sigma_{s}\) and \(f_p\) are bounded and continuous. Let
    \begin{align*}
        M:=\big\{(x,y,v,s) \in \partial X \times \partial X \times \mathbb{S}^2 \times (0,\infty):\|x-y\| < s ~\mathrm{and}~ v \in \partial V_{(-)}(y)\big\}
    \end{align*}
    Then for any \(k \geq 1\) and smooth function \(b_\mathfrak{m}:\partial X \times \mathbb{S}^2 \to \mathbb{R}\) 
    \begin{align}
        (x,y,v,s) \mapsto \mathfrak{m}_k(x,y,v,s):= \int_{\partial V_{(+)}(x)}\mathfrak{a}_k(x,v',y,v,s) b_\mathfrak{m}(x,v)dv' 
    \end{align}
    as well as 
    \begin{align}
        (x,y,v,s) \mapsto \mathfrak{m}(x,y,v,s):= \sum_{k=0}^{\infty} \mathfrak{m}_k(x,y,v,s)
    \end{align}
    are continuous functions when restricted to \(M\). In particular, for any \(x_D \in \partial X\) the expression \(\mathfrak{m}(\argdot \mid x_D): \partial V_{(-)}(x_D) \times (0,\infty) \to \mathbb{R}\)
    \begin{align*}
        (v,s) \mapsto \mathfrak{m}(x_D,x_D,v,s)
    \end{align*}
    is well defined and continuous.
    
\end{lemma}
\subsection{Semi-parametric form of the RTE-based inverse problem}\label{sec:uniqueness}
We need to rigorously address the question whether the collected data is sufficient for a reliable estimation of the dispersion, or equivalently differential absorption, parameters. Standard RTE related uniqueness conditions such as those given in \cite{ChoulliStefanov:1996,Bal09} for a variety of measurements types, require a detector almost everywhere on the boundary which may be realistic in some cases but clearly not when working with Lidar systems on such a large scale. The inverse transport problem with angularly resolved sources and angularly averaged measurements which is similar to our setting, apart from lacking temporal resolution, was considered in \cite{Langmore:2008} where uniqueness and stability results for the optical parameters are given. In \cite{GkioulekasLevinZickler2016} the authors present similarity relations, i.e. negative results, for the case of arbitrary heterogeneous materials. In our case, we don't seek to recover all optical parameters and additionally have a strong constraint based on the dispersion process which means that the right question to ask is whether, and if so to what extent, the dispersion constraint regularises the inverse problem and whether this is enough to recover the dispersion from single-detector data distorted by noise. In this section we focus on the uniqueness of the reconstruction for which we start with the following definition.
\begin{definition}\label{def:kernel_space}
For any function \(\phi \in C_c^{\infty}(\mathbb{R})\), i.e. \(\phi\) is a smooth function with compact support, we define the set of functions \(K(\phi) \subseteq C^\infty(X)\) as
\begin{align}\label{eq:kernel_space}
    K(\phi):= \left\{w_0+\sum_{j=1}^N w_j \phi\left(\frac{\|\argdot - m_j\|}{h_j}\right):N \in \mathbb{N}, m_j \in X, h_j > 0 \mathrm{~and~} w_j\in \mathbb{R} \right\}.
\end{align}
We may, without loss of generality, assume that \(\supp \phi \subseteq [-1,+1]\) and also define for any point \(x \in \partial X\) the set \(K(\phi \mid x)\) of functions that don't extend beyond \(x\) for such kernels \(\phi\). We may express this more formally as
\begin{align}\label{eq:kernel_space_x}
    K(\phi \mid x):= \left\{w_0+\sum_{j=1}^N w_j \phi\left(\frac{\|\argdot - m_j\|}{h_j}\right) \in K(\phi): h_j < \|x-m_j\|~ \forall j=1,\dots,N\right\}.
\end{align}
A function \(\phi \in C_c^{\infty}(\mathbb{R})\setminus \{0\}\) is called a good kernel, if its non-negative and decreases on \([0,\infty)\), i.e. for all \(z \geq 0\) we require \(\phi(z) \geq 0\) as well as  \(\phi'(z) \leq 0\).

\end{definition}
The following theorem shows that, under the condition from \cref{eq:equal_scattering} and fixed scattering parameters, a time resolved differential absorption measurement taken at the point \(x_D\) of the light source is good enough to recover the absorption difference \(\sigma_{\mathrm{on}(a)} - \sigma_{\mathrm{off}(a)} = \alpha \in K(\phi \mid x_D)\) for a good kernel \(\phi\), i.e. \(\alpha\) that admits a representation akin to that of \Cref{eq:kernel_parameterisation}. Apart from not being compactly supported, although numerically there is arguably no real difference between a Gaussian and an approximation thereof by an appropriately chosen compactly supported kernel \(\phi\), the function \(\exp(-x^2/2)\) is essentially a good kernel. Note that  \(K(\phi) \not\subseteq C_c^\infty(X)\), i.e. functions in \(K(\phi)\) may extend over the boundaries of \(X\), which is the case in situations such as when the gas is near the ground, and the condition \(\alpha \in K(\phi \mid x_D)\) states that \(x_D\) is outside the support of \(\alpha\) which is weaker than requiring \(\alpha \in C_c^\infty(X) \cap K(\phi)\). Intuitively, if \(\alpha \in K(\phi \mid x_D)\) and \(X\) is convex, then an observer standing at \(x_D\) can ``see'' the support of \(\alpha\)  without having to fully turn around. 
\paragraph{Simple uniqueness results} 
As \(K(\phi)\), and thus \(K(\phi \mid x_D)\), only contains functions that can be represented by a finite tuple of real-valued parameters, \cref{lem:fixed_scatter_uniqueness} can be seen as a uniqueness result for discretised parameters. However, we don't require explicit knowledge of \(\phi\), the number of kernel functions used or any other, possibly dispersion related, constraints that would restrict/regularise the problem further. Indeed, \(K(\phi \mid x_D)\) is infinite dimensional and looking for \(\alpha \in K(\phi \mid x_D)\) without further regularisation is still an ill-posed problem.

\begin{lemma}(Uniqueness for given scattering)\label{lem:fixed_scatter_uniqueness}
    Assume that \(X\) is convex as well as otherwise sufficiently regular, \(x_D \in \partial X\) and that \(\mathfrak{m}_{\mathrm{off}}\) is as in \cref{lem:path_continuity} for continuous optical parameters \(\sigma_{a}, \sigma_{s}\) and \(f_p\). If \(\alpha_1, \alpha_2  \in K(\phi \mid x_D)\) for a good kernel \(\phi\), then 
    \begin{align*}
        \alpha_1 = \alpha_2 \iff \mathfrak{m}^{(1)}_{\mathrm{on}}(\argdot \mid x_D) = \mathfrak{m}^{(2)}_{\mathrm{on}}(\argdot \mid x_D)
    \end{align*}
    where \(\mathfrak{m}^{(1/2)}_{\mathrm{on}}\) correspond to optical parameters \(\sigma_{a}+\alpha_{1/2}, \sigma_{s}\) and \(f_p\). In other words, the forward map with respect to differential absorption is injective on \(K(\phi \mid x_D)\) for any good kernel \(\phi\).
    \begin{proof}
    See \cref{proof:uniqueness}
    \end{proof}
\end{lemma}
In general the scattering parameters won't be fixed a priori and it is necessary to consider an inverse problem where \(\sigma_a, \sigma_s\) and \(f_p\) are unknown as well. 

In the case of narrow (or separately measured) FOVs the validity of a single scattering approximation means that the terms necessary to prove \cref{lem:fixed_scatter_uniqueness} are always accessible from the measured data which therefore provides sufficient information about the optical parameters in order to reconstruct the absorption field (which is the well known foundation for standard DIAL). The following result shows that if only a wide FOV is available the issue could essentially be further reduced to the influence of a subset of the optical parameters.
\begin{lemma}\label{lem:fixed_phase_uniqueness}
    Assume that \(X\) is convex as well as otherwise sufficiently regular, \(x_D \in \partial X\) and that \(\mathfrak{m}_{\mathrm{off}}\) is as in \cref{lem:path_continuity} for continuous optical parameters \(\sigma_{a}, \sigma_{s}\) and \(f_p\). If \(\sigma_{a}, \sigma_{s} \in K(\phi \mid x_D)\) for a good kernel \(\phi\) with a known proportionality constant between them while \(f_p\) is known, then \(\mathfrak{m}_{\mathrm{off}}\) uniquely determines \(\sigma_a\) and \(\sigma_s\).
    \begin{proof}
    See \cref{proof:uniqueness}
    \end{proof}
\end{lemma}

Combining \cref{lem:fixed_scatter_uniqueness} and \cref{lem:fixed_phase_uniqueness} we obtain that the measured data is arguably good enough to recover an additional scattering related parameter even if it's assumed to have a virtually unconstrained and very general form. However, since we measure light from where it was released (see \cref{fig:measurement_modality}) the measurements will be greatly impacted by the phase function or, more precisely, its back scattering coefficient. Although this could potentially be compensated by means of additional kernel functions, it unfortunately becomes infeasible to use this naturally occurring parameterisation when solving the inverse problem computationally given how it requires evaluations of the RTE for potentially complex and high-dimensional parameters. In the next paragraph we will present a way around this issue that allows us to use a high-dimensional parameterisation that (partially) accounts for the dominant scattering behaviour but doesn't present an impractical computational challenge.  

\paragraph{A semi-parametric relaxation}
From a practical perspective, one of the key features associated with DIAL is the ability to recover the absorption field without having to concern oneself with the the scattering parameters. On a more technical level however, the are the differential absorption field \(\alpha\) as well as \(\mathfrak{m}_{\mathrm{off}}\) for which, assuming a free line of sight and thus ignoring the boundary of the domain, we have
\begin{subequations}
\begin{align}
    \mathfrak{m}_{\mathrm{off}}(v,t \mid x_D) &\propto \frac{1}{t^2} e^{-2 \int_{0}^{t/2}\sigma_{a+s}(x_D+rv)dr} \sigma_{s}(x_D+tv) f_p(x_D+tv, -v\cdot v) \label{eq:off_data}\\
    \mathfrak{m}_{\mathrm{on}}(v,t \mid x_D) &\propto \frac{1}{t^2} e^{-2 \int_{0}^{t/2}\sigma_{a+s}(x_D+rv) + \alpha(x_D+rv)dr} \sigma_{s}(x_D+tv) f_p(x_D+tv, -v\cdot v) \label{eq:on_data}
\end{align}
\end{subequations}

The reconstruction of this parameter is of course trivial as it's measured directly but at the same time also all that is required in order to account for arbitrary scattering behaviour when the detector uses a narrow FOV. 

Note that the measurements, regardless of the FOV (and model), can be written as \((\tilde{\rho}_{\mathrm{on}}(v,t),\tilde{\rho}_{\mathrm{off}}(v,t))\) subject to the constraints
\begin{subequations}
\begin{align}
    \tilde{\rho}_{\mathrm{off}}(v,t) &= \mathfrak{m}_{\mathrm{off}}(v,t) \tilde{\psi}(v,t) \label{eq:off_constraint} \\
    \tilde{\rho}_{\mathrm{on}}(v,t) &= \mathfrak{m}_{\mathrm{on}}(v,t) \tilde{\psi}(v,t) \label{eq:on_constraint} \\
    \tilde{\psi}(v,t) &= 1 \label{eq:psi_constraint} 
\end{align}
\end{subequations}
where \(\mathfrak{m}_{\mathrm{on/off}}\) have the same underlying optical parameters and \(\tilde{\psi}\) is (for now) a redundant constant. In classical DIAL the constraint in \Cref{eq:psi_constraint} is dropped which results in what is essentially a modified forward model that, due to the multiplicative relationship between \Cref{eq:on_data} and \Cref{eq:off_data}, can capture all effects resulting from different \(\sigma_a, \sigma_s\) and \(f_p\). It is easily seen that the value for the differential absorption \(\alpha\) is still uniquely determined after this relaxation. The same is not true for the other optical parameters or \(\tilde{\psi}\).

Since \(\tilde{\psi}\) is a free parameter after the relaxation, we may divide by \(\mathfrak{m}_{\mathrm{off}}\) and obtain
\begin{subequations}
\begin{align}
    \rho_{\mathrm{off}}(v,t) &=  \psi(v,t) \label{eq:off_constraint_modified} \\
    \rho_{\mathrm{on}}(v,t) &= \frac{\mathfrak{m}_{\mathrm{on}}(v,t)}{\mathfrak{m}_{\mathrm{off}}(v,t)} \psi(v,t) \label{eq:on_constraint_modified} \\
    &= e^{-2 \int_{0}^{t/2}\alpha(x_D+rv)dr} \psi(v,t) \label{eq:on_constraint_ssct}
\end{align}
\end{subequations}
with \cref{eq:on_constraint_ssct} only being true for single scattering. Comparing Equations \cref{eq:off_data} and \cref{eq:on_data} with \cref{eq:on_constraint_ssct} and \cref{eq:off_constraint_modified} it becomes apparent that \(\psi\) essentially acts as a scaled back-scattering coefficient in the case of narrow FOVs and a single scattering model. 

Formally we may think of \(\psi\) as a density corresponding to photon detection. More specifically, if we denote by \(\Delta t\) the bin width, \(H_D\) the energy of the light source (expressed as expected number of photons per pulse) and \(A_D\) the detector size, then we have  
\[\mathbb{P}(\text{Photon from direction} ~ v ~ \text{observed in bin centred at} ~ t) \approx H_D^{-1}  \Delta t A_D \psi(v,t) .\]
Furthermore we also have
\[\mathbb{P}(\text{Differential absorption} \mid \text{Detection from direction} ~ v ~ \text{at time} ~ t) \approx 1-\frac{\mathfrak{m}_{\mathrm{on}}(v,t)}{\mathfrak{m}_{\mathrm{off}}(v,t)} .\]

Before we give an explanation for the roles of the quantities involved in this re-parameterised for of the model we give a uniqueness statement similar to those in \cref{lem:fixed_scatter_uniqueness} and \cref{lem:fixed_phase_uniqueness}. 

In view of \cref{eq:on_constraint_ssct} it becomes apparent that \(\psi\) is proportional to the local back-scattering rate \(\sigma_s(x_D+tv)f_p(x_D+tv, -v\cdot v)\). We may (loosely) associate \(\psi\) with a similar quantity also in a more general RTE setting. Indeed, the proof of \cref{lem:path_continuity} relied on the fact that a photon cannot be observed unless there was at least one scattering event resulting in a sufficiently large change of direction, the magnitude of which depends on the order of scattering under consideration and is maximal for single scattering. Assuming that the phase function \(f_p\) is (approximately) homogeneous in the spatial coordinate, which will be true if the same particles are responsible for the majority of scattering events throughout the domain, then the intensity along each path that corresponds to a low order detection will contain a factor that resembles the back scattering coefficient \(f_p(x_D+tv, -v\cdot v)\). In such situations our high-dimensional nuisance parameters \(\psi\) can therefore be viewed as a parameter that accounts for the back- or, more generally, large angle scattering rate.

\begin{theorem}\label{thm:relaxed_uniqueness}
    Assume that \(X\) is convex as well as otherwise sufficiently regular, \(x_D \in \partial X\) and that \(\mathfrak{m}_{\mathrm{off}}\) is as in \cref{lem:path_continuity} for continuous optical parameters \(\sigma_{a}, \sigma_{s}\) and \(f_p\). Assume that \(\sigma_{a}, \sigma_{s} \in K(\phi_\mathfrak{s} \mid x_D)\) for a good kernel \(\phi_\mathfrak{s}\) and that \(f_p\) is known. Further assume that \(\alpha \in K(\phi_\mathfrak{a} \mid x_D)\) for a good kernel \(\phi_\mathfrak{a}\) such that for any \(h > 0\)
    \begin{align}
       \frac{ \phi_\mathfrak{s}\left(1-\frac{\varepsilon}{h}\right) dz}{\int_{0}^{\varepsilon} \phi_\mathfrak{a}\left(1-z\right) dz} \xrightarrow{\varepsilon \to 0} 0 \label{eq:kernel_tails}.
    \end{align}
    Further assume that \(\sigma_{a}, \sigma_{s}\) and \(\alpha\) have common mid-points as well as widths while the weights are equal up to proportionality (with a constant shared between all summands) and the constant terms, which correspond to the homogeneous ambient quantities in the dispersion and are denoted by \(w_0\) in \cref{def:kernel_space}, are known. Then the absorption \(\frac{\mathfrak{m}_{\mathrm{off}}}{\mathfrak{m}_{\mathrm{on}}}\) uniquely determines \(\sigma_a\), \(\sigma_s\) and \(\alpha\).
    \begin{proof}
    Similar to \cref{lem:fixed_scatter_uniqueness}. Condition \cref{eq:kernel_tails} ensures that the scattering parameters do not ``interfere`` with single scattering as the uniqueness is based on the observation that single scattering ``precedes`` multiple scattering at the edge of a puff. It effectively requires that the scattering kernel has diminishing mass in its tails where its contribution to the signal will primarily be single/back-scattering and is meant to be compensated by \(\psi\) and ensures that we can treat reconstructions ``within`` the plume in the same way as those at the ``edges``. In fact, if the plume was given by a single kernel function this condition would not be necessary. Since the regularity imposed by the dispersion model make the plume behave much more like a single structre as opposed to many independent puffs, we believe that this technical condition is best thought of as an artifact of our proof.
    Note that regardless of \(\phi_\mathfrak{a}\) one can always find a kernel \(\phi_\mathfrak{s}\) that is arbitrarily similar to \(\phi_\mathfrak{a}\) in the \(L_1\) sense but satisfies \cref{eq:kernel_tails}. The details can be found in \cref{sec:proofs}.
    \end{proof}
\end{theorem}

The essence of \cref{thm:relaxed_uniqueness} is that the parameters associated with the relaxed model as given in \cref{eq:off_constraint_modified} and \cref{eq:on_constraint_modified} can be identified from the data if the scattering and absorption particles follow the same atmospheric dispersion process. Knowledge of the ambient quantities is in principle not necessary but arguably more realistic than an accurate reconstruction from noisy measurements. It should be noted that conditions such as the required proportionality of parameters are somewhat necessary since we cannot reconstruct anything from differential absorption data in the absence of the trace gas. In particular, any obstacles before the plume can only be modeled as a \emph{known} baseline concentration of scattering particles. We should also note that the conditions imposed as \cref{eq:kernel_tails} are such that the resulting function spaces can produce a more accurate approximation in those scenarios where the application of wide FOV measurements in more appropriate but also highlight an inherent limitation of wider FOVs and the RTE model. If exact reconstruction is sought then we rely on separation of single and multiple scattering which is possible only locally at times where we measure signals scattered at plume boundaries. Such data is fundamentally unstable as not only the tails of a distribution are can be heavily perturbed by even small changes in the parameters but the signal will also be particularly weak, noisy and dominated by ambient photons. 

Nonetheless \cref{thm:relaxed_uniqueness} can be used as a best-case reference as to which parameters we can hope to resolve by wide FOVs and although exact reconstruction is problematic, it turns out that other aspects of the measurement are more stable under perturbations of the scattering parameters than the tails and wide FOV measurements can still bear useful information with regards to the differential absorption. The necessary results are developed in \cref{sec:stability}.

\paragraph{Path spaces and entropy}
Recall that in the case of single scattering, i.e. sufficiently narrow FOVs, the differential absorption is independent of \(f_p\) and not being able to reconstruct it has no effect on the reconstructed gas concentration. This is of course not true in a more general RTE setting and yet \cref{thm:relaxed_uniqueness} requires knowledge of \(f_p\). As argued earlier, the nuisance parameter \(\psi\) can account for variability in the back scattering behaviour, which is arguably the most important aspect, but even if \(f_p\) is assumed spatially homogeneous, i.e. \(f_p(x, v \cdot v')=g_p(v \cdot v')=f_p(y, v \cdot v')\) for any \(x,y \in X\) and some spatially homogeneous \(g\), it may still have unaccounted irregularities on \([-1+\gamma,1]\) where \(\gamma > 0\) serves (informally) as cut-off for large angles. In practice, the most commonly modeled component of \(f_p\) is the forward peak, e.g. by means of the Henyey-Greenstein phase function \cite{HenyeyGreenstein:1941}. If, as we assume, the distribution of direction depends only on the cosine of inward and outward direction \(v \cdot v'\), this is somewhat equivalently to modelling the first angular moments. As mentioned earlier, the introduction of \(\psi\) implies that the modelled intensity of measured light is independent of the RTE scattering parameters, particularly \(\sigma_s\). Note that up to first order we can approximate \(g \approx \lambda \delta_{\{v \cdot v' = 1\}} + (1-\lambda)g_0 =: g_\delta\) by putting
\[
\lambda  = \int_{\mathbb{S}^2} v \cdot v' g(v\cdot v') dv'
\]
and picking \(g_0\) such that the analogous integral is 0. It is not hard to see that solving the RTE \cref{eq:time_dependent_RTE} with \(f_p \to g_\delta\) is the same as setting \(f_p \to g_0\) as well as \(\sigma_s \to (1-\lambda) \sigma_s\). Note that this observation is the foundation for other commonly used phase function approximations that use \(\delta\)-distributions for forward peaks, e.g. the \(\delta\)-Eddington and related methods \cite{Wiscombe:1976,Wiscombe:1977}. The inverse idea is also true, i.e. one may also select a baseline \(g_0\) with a forward peak, which effectively results in an implicit choice of a baseline which gets homogenised as \(\sigma_s\) increases. 

\begin{figure}[h]%
    \centering
    \subfloat[\(\sigma_s\) high]{\scalebox{0.25}{\begin{tikzpicture}[ray/.style={decoration={markings,mark=at position .5 with {\arrow[>=latex]{>}}},postaction=decorate}]

		\node (0) at (-1.75, 3) {};
		\node (1) at (-0.25, 5) {};
		\node (2) at (2.75, 5.75) {};
		\node (3) at (6, 5) {};
		\node (5) at (-1.25, -5) {};
		\node (6) at (2.75, -7) {};
		\node (7) at (6, -6.75) {};
		\node (8) at (8, -5.25) {};
		\node (9) at (-2.75, -1.25) {};
		\node (10) at (7.75, 3.75) {};
		\node (11) at (8.75, 1.25) {};
		\node (12) at (9.75, -2.25) {};
		\node (13) at (-0.75, -4.5) {};
		\node (14) at (2.75, -6.5) {};
		\node (15) at (5, -5.75) {};
		\node (16) at (7.75, -4.75) {};
		\node (17) at (9, -2.25) {};
		\node (18) at (8, 1.25) {};
		\node (19) at (2.75, 5) {};
		\node (20) at (0.5, 4.5) {};
		\node (21) at (-1, 2.75) {};
		\node (22) at (-1.75, -1.25) {};
		\node (23) at (5.75, 4) {};
		\node (24) at (6.75, 3) {};
		\node (25) at (-2.25, -3.5) {};
		\node (26) at (-1.75, -3.25) {};
		\draw [semithick, dashed, bend right=90, looseness=1.50] (5.center) to (6.center);
		\draw [semithick, dashed, bend right=60] (6.center) to (7.center);
		\draw [semithick, dashed, bend right=45, looseness=0.75] (7.center) to (8.center);
		\draw [semithick, dashed, bend right=60, looseness=1.25] (1.center) to (0.center);
		\draw [semithick, dashed, bend right=45] (2.center) to (1.center);
		\draw [semithick, dashed, bend left=60] (2.center) to (3.center);
		\draw [semithick, dashed, bend left=60, looseness=1.25] (9.center) to (0.center);
		\draw [semithick, dashed, bend left=60] (12.center) to (8.center);
		\draw [semithick, dashed, bend left=75, looseness=1.25] (3.center) to (10.center);
		\draw [semithick, dashed, bend left=60] (10.center) to (11.center);
		\draw [semithick, dashed, bend left=75] (11.center) to (12.center);
		\draw [semithick, dashed, bend left=15] (1.center) to (20.center);
		\draw [semithick, dashed, bend right=45, looseness=0.75] (2.center) to (19.center);
		\draw [semithick, dashed, bend left=15] (3.center) to (23.center);
		\draw [semithick, dashed, bend right] (10.center) to (24.center);
		\draw [semithick, dashed, bend left=45] (11.center) to (18.center);
		\draw [semithick, dashed, bend left] (12.center) to (17.center);
		\draw [semithick, dashed, bend right] (8.center) to (16.center);
		\draw [semithick, dashed, bend right] (15.center) to (7.center);
		\draw [semithick, dashed, bend left=15, looseness=0.75] (14.center) to (6.center);
		\draw [semithick, dashed, bend left=45, looseness=0.75] (13.center) to (5.center);
		\draw [semithick, dashed, bend right] (9.center) to (22.center);
		\draw [semithick, dashed, bend right] (0.center) to (21.center);
		\draw [semithick, dashed, bend right=75, looseness=1.25] (9.center) to (25.center);
		\draw [semithick, dashed, bend right=60, looseness=1.25] (25.center) to (5.center);
            \draw[thick,fill=white,white] (-3.6,.235) ellipse (.05 and .05);
            \draw[thick,fill=white,white] (-2.84,-1.145) ellipse (.05 and .05);
            \draw[thick,fill=white,white] (-2.81,-1.325) ellipse (.05 and .05);
            \draw[thick,fill=white,white] (-3.4,1.88) ellipse (.15 and .15);
            \draw[thick,fill=white,white] (-1.6,2.88) ellipse (.05 and .05);

    \draw[thick,gray] (-12,-.30) ellipse (0.11 and 0.33);

    \draw[thick,red!35] (-12,.36) -- (-8,.3225);
    \draw[thick,red!35] (-12,.34) -- (-8,.3025);
    \draw[thick,red!35] (-8,.3025) -- (0,.2275);
    \draw[thick,red!35] (-12,.32) -- (-8,.2825);
    \draw[thick,red!35] (-8,.2825) -- (2,0.18875);
    \draw[thick,red!35] (-12,.30) -- (-8,.2625);
    \draw[thick,red!35] (-8,0.2625) -- (-2,0.20625);
    \draw[thick,red!35] (-12,.26) -- (-8,.2225);
    \draw[thick,red!35] (-8,.2225) -- (-6,0.20375);
    \draw[thick,red!35] (-12,.28) -- (-8,.2425);
    \draw[thick,red!35] (-8,.2425) -- (-3,0.205);
    \draw[dashed,thick,red!35,-stealth] (2,0.18875) -- (6,0.105);

    \draw [dashed, thin, -stealth, red](-13.150,.36) -- (-12.1,.36);
    \draw [dashed, thin, -stealth, red](-13.045,.26) -- (-12.1,.26);

    \draw[line width=2.5pt,gray!70,ray] (-8,.3225) -- (-5,3.4025);
    \draw[line width=2.5pt,gray!70,ray] (-5,3.4025) -- (-3,1.4025);
    \draw[line width=2.5pt,gray!70,ray] (-3,1.4025) -- (-2.8,-1.7025);
    \draw[line width=2.5pt,gray!70,ray] (-2.8,-1.7025) -- (-12,-.41);
    \draw[line width=2.5pt,gray!70,-stealth] (-2.8,-1.7025) -- (-12,-.41);

    \draw[line width=2.5pt,gray!70,ray] (0,.2275) -- (-0.1,-.525);
    \draw[line width=2.5pt,gray!70,ray] (-0.1,-.525) -- (0.1,-1.3825);
    \draw[line width=2.5pt,gray!70,ray] (0.1,-1.3825) -- (-0.3,-2.825);
    \draw[line width=2.5pt,gray!70,ray] (-0.3,-2.825) -- (-12,-.48);
    \draw[line width=2.5pt,gray!70,-stealth] (-0.3,-2.825) -- (-12,-.48);

    \draw[line width=1.2pt,gray!70, dashed, ray] (2,0.18875) -- (3.1,0.75);
    \draw[line width=1.2pt,gray!70, dashed,ray] (3.1,0.75) -- (-12,-.27);
    \draw[line width=1.2pt,gray!70, dashed,-stealth] (3.1,0.75) -- (-12,-.27);

    \draw[line width=1.2pt,gray!70, dashed,ray] (-2,0.20625) -- (1.5,1.25);
    \draw[line width=1.2pt,gray!70, dashed,ray] (1.5,1.25) -- (0.75,2);
    \draw[line width=1.2pt,gray!70, dashed,ray] (0.75,2) -- (-12,-.20);
    \draw[line width=1.2pt,gray!70, dashed,-stealth] (0.75,2) -- (-12,-.20);

    \draw[line width=1.2pt,gray!70, dashed,ray] (-3,0.205) -- (2.5,-.62425);
    \draw[line width=1.2pt,gray!70, dashed,ray] (2.5,-.62425) -- (-12,-.34);
    \draw[line width=1.2pt,gray!70, dashed,-stealth] (2.5,-.62425) -- (-12,-.34);

    \draw[line width=2.5pt,gray!70,ray] (-6,0.20375) -- (-1.45,1.8825);
    \draw[line width=2.5pt,gray!70,ray] (-1.45,1.8825) -- (-0.15,2.9325);
    \draw[line width=2.5pt,gray!70,ray] (-0.15,2.9325) -- (-12,-.13);
    \draw[line width=2.5pt,gray!70,-stealth] (-0.15,2.9325) -- (-12,-.13);

\end{tikzpicture}}}%
    \subfloat[\(\sigma_s\) low]{\scalebox{0.25}{\begin{tikzpicture}[ray/.style={decoration={markings,mark=at position .5 with {\arrow[>=latex]{>}}},postaction=decorate}]

		\node (0) at (-1.75, 3) {};
		\node (1) at (-0.25, 5) {};
		\node (2) at (2.75, 5.75) {};
		\node (3) at (6, 5) {};
		\node (5) at (-1.25, -5) {};
		\node (6) at (2.75, -7) {};
		\node (7) at (6, -6.75) {};
		\node (8) at (8, -5.25) {};
		\node (9) at (-2.75, -1.25) {};
		\node (10) at (7.75, 3.75) {};
		\node (11) at (8.75, 1.25) {};
		\node (12) at (9.75, -2.25) {};
		\node (13) at (-0.75, -4.5) {};
		\node (14) at (2.75, -6.5) {};
		\node (15) at (5, -5.75) {};
		\node (16) at (7.75, -4.75) {};
		\node (17) at (9, -2.25) {};
		\node (18) at (8, 1.25) {};
		\node (19) at (2.75, 5) {};
		\node (20) at (0.5, 4.5) {};
		\node (21) at (-1, 2.75) {};
		\node (22) at (-1.75, -1.25) {};
		\node (23) at (5.75, 4) {};
		\node (24) at (6.75, 3) {};
		\node (25) at (-2.25, -3.5) {};
		\node (26) at (-1.75, -3.25) {};
		\draw [semithick, dashed, bend right=90, looseness=1.50] (5.center) to (6.center);
		\draw [semithick, dashed, bend right=60] (6.center) to (7.center);
		\draw [semithick, dashed, bend right=45, looseness=0.75] (7.center) to (8.center);
		\draw [semithick, dashed, bend right=60, looseness=1.25] (1.center) to (0.center);
		\draw [semithick, dashed, bend right=45] (2.center) to (1.center);
		\draw [semithick, dashed, bend left=60] (2.center) to (3.center);
		\draw [semithick, dashed, bend left=60, looseness=1.25] (9.center) to (0.center);
		\draw [semithick, dashed, bend left=60] (12.center) to (8.center);
		\draw [semithick, dashed, bend left=75, looseness=1.25] (3.center) to (10.center);
		\draw [semithick, dashed, bend left=60] (10.center) to (11.center);
		\draw [semithick, dashed, bend left=75] (11.center) to (12.center);
		\draw [semithick, dashed, bend left=15] (1.center) to (20.center);
		\draw [semithick, dashed, bend right=45, looseness=0.75] (2.center) to (19.center);
		\draw [semithick, dashed, bend left=15] (3.center) to (23.center);
		\draw [semithick, dashed, bend right] (10.center) to (24.center);
		\draw [semithick, dashed, bend left=45] (11.center) to (18.center);
		\draw [semithick, dashed, bend left] (12.center) to (17.center);
		\draw [semithick, dashed, bend right] (8.center) to (16.center);
		\draw [semithick, dashed, bend right] (15.center) to (7.center);
		\draw [semithick, dashed, bend left=15, looseness=0.75] (14.center) to (6.center);
		\draw [semithick, dashed, bend left=45, looseness=0.75] (13.center) to (5.center);
		\draw [semithick, dashed, bend right] (9.center) to (22.center);
		\draw [semithick, dashed, bend right] (0.center) to (21.center);
		\draw [semithick, dashed, bend right=75, looseness=1.25] (9.center) to (25.center);
		\draw [semithick, dashed, bend right=60, looseness=1.25] (25.center) to (5.center);
            \draw[thick,fill=white,white] (-3.6,.235) ellipse (.05 and .05);
            \draw[thick,fill=white,white] (-2.84,-1.145) ellipse (.05 and .05);
            \draw[thick,fill=white,white] (-2.81,-1.325) ellipse (.05 and .05);
            \draw[thick,fill=white,white] (-3.4,1.88) ellipse (.15 and .15);
            \draw[thick,fill=white,white] (-1.6,2.88) ellipse (.05 and .05);

    \draw[thick,gray] (-12,-.30) ellipse (0.11 and 0.33);

    \draw[thick,red!35] (-12,.36) -- (-8,.3225);
    \draw[thick,red!35] (-12,.34) -- (-8,.3025);
    \draw[thick,red!35] (-8,.3025) -- (0,.2275);
    \draw[thick,red!35] (-12,.32) -- (-8,.2825);
    \draw[thick,red!35] (-8,.2825) -- (2,0.18875);
    \draw[thick,red!35] (-12,.30) -- (-8,.2625);
    \draw[thick,red!35] (-8,0.2625) -- (-2,0.20625);
    \draw[thick,red!35] (-12,.26) -- (-8,.2225);
    \draw[thick,red!35] (-8,.2225) -- (-6,0.20375);
    \draw[thick,red!35] (-12,.28) -- (-8,.2425);
    \draw[thick,red!35] (-8,.2425) -- (-3,0.205);
    \draw[dashed,thick,red!35,-stealth] (2,0.18875) -- (6,0.105);

    \draw [dashed, thin, -stealth, red](-13.150,.36) -- (-12.1,.36);
    \draw [dashed, thin, -stealth, red](-13.045,.26) -- (-12.1,.26);

    \draw[line width=1.2pt,gray!70,dashed,ray] (-8,.3225) -- (-5,3.4025);
    \draw[line width=1.2pt,gray!70,dashed,ray] (-5,3.4025) -- (-3,1.4025);
    \draw[line width=1.2pt,gray!70,dashed,ray] (-3,1.4025) -- (-2.8,-1.7025);
    \draw[line width=1.2pt,gray!70,dashed,ray] (-2.8,-1.7025) -- (-12,-.41);
    \draw[line width=1.2pt,gray!70,dashed,-stealth] (-2.8,-1.7025) -- (-12,-.41);

    \draw[line width=1.2pt,gray!70,dashed,ray] (0,.2275) -- (-0.1,-.525);
    \draw[line width=1.2pt,gray!70,dashed,ray] (-0.1,-.525) -- (0.1,-1.3825);
    \draw[line width=1.2pt,gray!70,dashed,ray] (0.1,-1.3825) -- (-0.3,-2.825);
    \draw[line width=1.2pt,gray!70,dashed,ray] (-0.3,-2.825) -- (-12,-.48);
    \draw[line width=1.2pt,gray!70,dashed,-stealth] (-0.3,-2.825) -- (-12,-.48);

    \draw[line width=2.5pt,gray!70, ray] (2,0.18875) -- (3.1,0.75);
    \draw[line width=2.5pt,gray!70,ray] (3.1,0.75) -- (-12,-.27);
    \draw[line width=2.5pt,gray!70,-stealth] (3.1,0.75) -- (-12,-.27);

    \draw[line width=2.5pt,gray!70,ray] (-2,0.20625) -- (1.5,1.25);
    \draw[line width=2.5pt,gray!70,ray] (1.5,1.25) -- (0.75,2);
    \draw[line width=2.5pt,gray!70,ray] (0.75,2) -- (-12,-.20);
    \draw[line width=2.5pt,gray!70,-stealth] (0.75,2) -- (-12,-.20);

    \draw[line width=2.5pt,gray!70,ray] (-3,0.205) -- (2.5,-.62425);
    \draw[line width=2.5pt,gray!70,ray] (2.5,-.62425) -- (-12,-.34);
    \draw[line width=2.5pt,gray!70,-stealth] (2.5,-.62425) -- (-12,-.34);

    \draw[line width=1.2pt,gray!70,dashed,ray] (-6,0.20375) -- (-1.45,1.8825);
    \draw[line width=1.2pt,gray!70,dashed,ray] (-1.45,1.8825) -- (-0.15,2.9325);
    \draw[line width=1.2pt,gray!70,dashed,ray] (-0.15,2.9325) -- (-12,-.13);
    \draw[line width=1.2pt,gray!70,dashed,-stealth] (-0.15,2.9325) -- (-12,-.13);

\end{tikzpicture}}}%
    \caption{An increase in the scattering function \(\sigma_s\) means light is more likely to scatter resulting in less low order scattering at distance corresponding to deep inside the plume. As the non-parametric component controls the response strength it acts as a normalising constant while \(\sigma_s\) controls the diffusivity and thus the likely trajectories of the measured photons. A high amount of scattering particles as shown in image (a) favours high order scattering (solid line) over forward peaked trajectories that reach deeper inside the cloud (dashed lines) whereas the opposite can be observed in figure (b).}%
    \label{fig:wFOV_weighting}%
\end{figure}
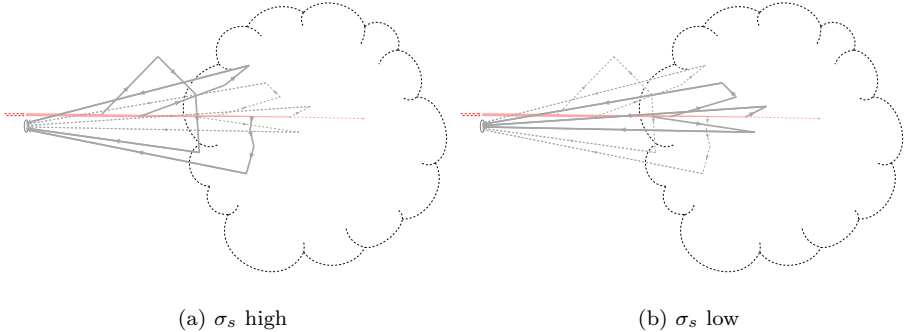

Since the above developments indicate that the optical parameters that control the absorption \(\frac{\mathfrak{m}_{\mathrm{on}}(v,t)}{\mathfrak{m}_{\mathrm{off}}(v,t)}\) can be used to model forward peaked behaviour, \cref{thm:relaxed_uniqueness} can be interpreted as a set of restrictions on \(\sigma_{a,s}\) which allow partial reconstruction of the phase function or at least certain aspects thereof (see \cref{fig:wFOV_weighting}).

Given that we cannot evaluate the forward model without a full set of optical parameters, we must find a way to deal with the lack of information regarding \(f_p\) or at least find a suitable candidate for \(g_0\). In fact, it only makes sense to use wider FOVs, or even multiple FOVs where the narrow component can be separated, if we handle the resulting ill-posedness without falling victim to regularisation induced errors/biases that are larger than a simple narrow FOV reconstruction. In theory we can select \(g_0\) that maximises the entropy on the sphere \(\mathbb{S}^2\) subject to given constraints which should be based on prior information/beliefs. It is worth noticing however that we found the choice of \(g_0\) to be of diminishing importance when it comes to the actual performance of our algorithms. An indication as to why this phenomenon will likely be observed more generally, at least in the situations of relevance for our measurements and methods, is given in the next section.

\section{Signal quality and noise considerations}
\label{sec:stability}
The results that were shown in \cref{lem:fixed_scatter_uniqueness}, \cref{lem:fixed_phase_uniqueness} as well as \cref{thm:relaxed_uniqueness} are similar to what is available for the classical DIAL approaches based on single scattering but weaker in certain ways. Combining \cref{lem:fixed_scatter_uniqueness} with \cref{lem:fixed_phase_uniqueness} essentially yields the result of classical DIAL with the additional, unrealistic assumption such as that the phase function be known. The more practical result in \cref{thm:relaxed_uniqueness} uses the same data as classical dial but only provides what is perhaps best thought of as regularity conditions under which a subset of the optical parameters is uniquely determined by differential absorption only.
Notably each of our proofs relied on extracting the effect of single scattering in order to identify the perturbation \(\alpha\) in the absorption function and we must ask ourselves whether there is any benefit in taking wider FOVs and multiple scattering into account. If we can separate first and higher order scattering contributions, e.g. through partially resolving the incident angles, then a good reconstruction algorithm won't make things worse since we could just ignore all but the narrow FOV with the single scattering measurement. In situations where such a separation is not an option, things get more complicated and angularly averaged point measurements from wide FOVs can be better or worse than narrow ones, depending on the situation. 
\subsection{Poisson noise model}
In order to make more formal statements, we must first introduce a noise model for our optical measurement. If the light source releases two instantaneous impulses at time \(t=t_{\mathrm{on/off}}\) in direction \(v\) consisting of \(\mathsf{Poisson}(H_D)\) many photons then the fully time resolved observation at the detector in the interval \((t_{\mathrm{on/off}},t_{\mathrm{on/off}}+T)\) will take the form of a Poisson point process with intensity measure \(H_D \mathfrak{m}_{\mathrm{on/off}}(v, \argdot \mid x_D)\) respectively. In practice it is arguably more realistic to have measurements that are aggregated into \(N_t\) bins, rather than infinitely time resolved, which would mean that for each direction \(v_i\) for \(i = 1, \dots, N_v\), in which a light pulse is released, we observe two independent random arrays \(\bm{m}_{v_i}, \bm{n}_{v_i}\) with independent entries such that, at least approximately,
\begin{align}\label{eq:data_noise_model}
\begin{split}
    \bm{m}_{v_i,t_j} &\sim \mathsf{Poisson}\left(\Delta t H_D A_D \mathfrak{m}_{\mathrm{on}} (v_i, t_j\mid x_D)\tilde{\psi}_{i,j}\right) \\
    \bm{n}_{v_i,t_j} &\sim \mathsf{Poisson}\left(\Delta t H_D A_D\mathfrak{m}_{\mathrm{off}} (v_i, t_j\mid x_D)\tilde{\psi}_{i,j}\right)    
\end{split}
\end{align}
where \(\Delta t\) is the bin width, \(A_D\) the detector size, \(t_j\) is the centre of bin \(j\) for \(j = 1, \dots, N_t\) and a mid-point quadrature rule was used to integrate \(\mathfrak{m}_{\mathrm{off}} (v_i, \argdot \mid x_D)\) over the bin and detector surface. The values \(\tilde{\psi}_{i,j}\) are a discrete version of the semi parametric component introduced in \cref{sec:uniqueness}. It should be noted that \Cref{eq:data_noise_model} ignores contributions of ambient light which change the mean but not the fact that the data has a \(\mathsf{Poisson}\) distribution. The effect this has on image reconstruction quality is discussed later in this section. Realistically, there are likely other detector related sources of distortion in our data apart from the \(\mathsf{Poisson}\) nature of the optical measurement but for weak signals, i.e. very low photon counts per bin, or, perhaps more realistically, low differential absorption we can expect that the optical noise dominates and our assumption in \Cref{eq:data_noise_model} becomes a good approximation for the measurement. Recall that we want to use such data \((\bm{m},\bm{n})\) to find the differential absorption 
\begin{align*}
    \alpha(x) &= \sigma_{a(\mathrm{on})}(x) - \sigma_{a(\mathrm{off})}(x)
\end{align*}
where the spatial absorption difference \(\alpha=\alpha[\theta]\) is parameterised by \(\theta \in \Theta\) and has the same form as in \Cref{eq:differential_absorption_parameterisation}. Consequently \(\mathfrak{m}_{\mathrm{on/off}}=\mathfrak{m}_{\mathrm{on/off}}[\theta]\) are functions of \(\theta\) and we may form the negative \(\log\)-likelihood \(\mathsf{L}(\tilde{\psi},\theta \mid \bm{m},\bm{n})\) for data \((\bm{m},\bm{n})\) as in \Cref{eq:data_noise_model}
\begin{align}\label{eq:DIAL_llik}
\begin{split}
    \mathsf{L}(\tilde{\psi},\theta \mid \bm{m},\bm{n}) &= \sum_{i=1}^{N_v}\sum_{j = 1}^{N_t} H_D \Delta t A_D (\mathfrak{m}_{\mathrm{on}}(v_i,t_j\mid x_D)+\mathfrak{m}_{\mathrm{off}}(v_i,t_j\mid x_D))\tilde{\psi}_{i,j} \\
    &-\bm{m}_{v_i,t_j} \log(\tilde{\psi}_{i,j}\mathfrak{m}_{\mathrm{on}}(v_i,t_j\mid x_D))  - \bm{n}_{v_i,t_j} \log(\tilde{\psi}_{i,j}\mathfrak{m}_{\mathrm{off}}(v_i,t_j\mid x_D)).
\end{split}
\end{align}
For that parameter we may compute the (efficient) Fisher information matrix 
\begin{align}\label{eq:DIAL_fisher_information_total}
    \mathsf{J}(\theta) = \mathbb{E}(\partial_{\theta\theta}\max_{\tilde{\psi} \in \Psi}\mathsf{L}(\theta \mid \bm{m},\bm{n}))
\end{align}
where \(\Psi\) is a (typically non-trivial) set of feasible values for the array \(\tilde{\psi}\) which can be used to enforce a set of constraints such as \cref{eq:psi_constraint}. Assuming that there is \(\theta^* \in \Theta\) such that \(\alpha[\theta^*]\) is the true difference in absorption, then we know that under suitable regularity conditions \cite{Murphy:2000} the maximum likelihood estimator \(\hat{\theta}_{\mathrm{MLE}}\) satisfies
\begin{align}\label{eq:MLE_AN}
    \hat{\theta}_{\mathrm{MLE}} \sim \mathsf{AN}(\theta^*, \mathsf{J}(\theta^*)^{-1}),
\end{align}
where \(\mathsf{AN}\) denotes asymptotic normality. In certain situations \cref{eq:MLE_AN} can be a reasonable approximation for the distribution of \(\hat{\theta}_{\mathrm{MLE}}\) but even if we ignore that fact we easily see that for simple situation, such as the one considered in \cref{ssec:trivial}, \Cref{eq:MLE_AN} is essentially the square of the sensitivity of the noiseless measurement, i.e. the mean, scaled by its standard deviation which a measure for the amount of noise which is arguably a reasonable measure of signal quality in \(\mathfrak{m}_{\mathrm{on/off}}\) for a particular set of optical parameters.
\subsection{Trivial edge cases}\label{ssec:trivial}
Before considering more general scenarios we can develop a basic intuition by looking at the situation of known scattering as in \cref{lem:fixed_scatter_uniqueness} and in particular \(\tilde{\psi}_{i,j}=1\) fixed. In that case we have
\begin{align*}
    \mathsf{J}(\theta) = \sum_{i=1}^{N_v}\sum_{j = 1}^{N_t}  J_{\mathrm{on}}(v_i,t_j\mid x_D)[\theta]+J_{\mathrm{off}}(v_i,t_j\mid x_D)[\theta]
\end{align*}
where we abbreviated
\begin{align}\label{eq:DIAL_fisher_information_single}
    J_{\mathrm{on/off}}(v_i,t_j\mid x_D)[\theta] = H_D \Delta t A_D \frac{\partial_{\theta}\mathfrak{m}_{\mathrm{on/off}}(v_i,t_j\mid x_D) \partial_{\theta}\mathfrak{m}_{\mathrm{on/off}}(v_i,t_j\mid x_D)^{\top}}{\mathfrak{m}_{\mathrm{on/off}}(v_i,t_j\mid x_D)}
\end{align}
and suppressed the \(\theta\)-dependence in \(\mathfrak{m}_{\mathrm{on/off}}\). 
\paragraph{Narrow vs. wide vs. multiple FOVs}
Note that we can think of averaged measurements from wider FOVs, corresponding to choices of the function \(b\) in \cref{lem:path_continuity} with a larger support, as the sum of single and multiple scattering photons which are independent \(\mathsf{Poisson}\) distributed random variables with intensity \(\Delta t H_D \mathfrak{m}_{\mathrm{on/off},1}\) and \(\Delta t H_D \sum_{j=2}^\infty \mathfrak{m}_{\mathrm{on/off},j}\) respectively. First consider a constant absorption difference, i.e. we have 
\begin{align*}
    C_{\mathrm{ambient}}[\theta] = \sigma_{a(\mathrm{on})}(x) - \sigma_{a(\mathrm{off})}(x)
\end{align*}
alongside some known \(\sigma_s,\sigma_a > 0\) on \(X\). We easily see that 
\begin{align*}
    \mathfrak{m}_{\mathrm{on}}(v_i,t_j\mid x_D)=\mathfrak{m}_{\mathrm{off}}(v_i,t_j\mid x_D)\exp(-t_j C_{\mathrm{ambient}}).
\end{align*} 
for any \(i=1,\dots,N_v\), \(j=1,\dots,N_t\) and \(\theta \in \Theta\). If the scattering parameters, and thus \(\mathfrak{m}_{\mathrm{off}}\), are fixed then \(\mathfrak{m}_{\mathrm{off}}\) is independent of \(\theta\) and we have \(J_{\mathrm{off}} = 0\) while
\begin{align}\label{eq:fully_homogeneous_fisher_information}
    J_{\mathrm{on,wide}}(v_i,t_j\mid x_D)[\theta] =  J_{\mathrm{on,narrow}}(v_i,t_j\mid x_D)[\theta] \frac{\mathfrak{m}_{\mathrm{off,wide}}(v_i,t_j\mid x_D)}{\mathfrak{m}_{\mathrm{off,narrow}}(v_i,t_j\mid x_D)}
\end{align}
where \(J_{\mathrm{on,wide}}\) and \(J_{\mathrm{on,narrow}}\) are as in \Cref{eq:DIAL_fisher_information_single} and correspond to wide and narrow FOVs respectively. As \(\mathfrak{m}_{\mathrm{off,wide}} \geq \mathfrak{m}_{\mathrm{off,narrow}}\) the information in each measurement point increases, the overall quality of the signal would increase considerably in optically thick environments. 

The case where the gas is spread evenly in the domain is in a way ideal because we don't have to worry about how light came to the detector and averaging over a wider FOV has no downsides. This situation can be interpreted as the extreme case of a very large gas plume where its size is taken to \(\infty\). Its counterpart can be seen as the situation where the same, unknown amount \(C_{\mathrm{ambient}}[\theta] \mathrm{vol}(X)\) of gas is taken and accumulated in a very small area around a point \(\frac{t_j}{2}v_i + x_D\). As the area is taken to \(0\) we end up with \(\partial_\theta \mathfrak{m}_{\mathrm{on},k}(v_i,t_j\mid x_D) \to 0\) for \(k \geq 2\) since most of the time multiply scattered light won't pass through that small area where the gas is located. Single scattering on the other hand is much more localised and remains sensitive to arbitrarily small objects. As an immediate consequence we obtain
\begin{align}\label{eq:fully_localised_fisher_information}
    J_{\mathrm{on,wide}}(v_i,t_j\mid x_D)[\theta] \approx J_{\mathrm{on,narrow}}(v_i,t_j\mid x_D)[\theta] \frac{\mathfrak{m}_{\mathrm{off,narrow}}(v_i,t_j\mid x_D)}{\mathfrak{m}_{\mathrm{off,wide}}(v_i,t_j\mid x_D)}
\end{align}
which is exactly the reverse scaling compared to \Cref{eq:fully_homogeneous_fisher_information}. For times larger than \(t_j\) the multiply scattered measurement component does retain some sensitivity but the overall situation is considerably worse than in the case of what can be thought of as an infinitely spread out gas plume to the point where averaging over wider FOVs will likely make things worse rather than yield an improvement in any meaningful way.
\begin{figure}[h]%
    \centering
    \subfloat[large absorbing region]{\scalebox{0.25}{\begin{tikzpicture}[ray/.style={decoration={markings,mark=at position .5 with {\arrow[>=latex]{>}}},postaction=decorate}]

            \node (0) at (-1.75, 3) {};
		\node (1) at (-0.25, 5) {};
		\node (2) at (2.75, 5.75) {};
		\node (3) at (6, 5) {};
		\node (5) at (-1.25, -5) {};
		\node (6) at (2.75, -7) {};
		\node (7) at (6, -6.75) {};
		\node (8) at (8, -5.25) {};
		\node (9) at (-2.75, -1.25) {};
		\node (10) at (7.75, 3.75) {};
		\node (11) at (8.75, 1.25) {};
		\node (12) at (9.75, -2.25) {};
		\node (13) at (-0.75, -4.5) {};
		\node (14) at (2.75, -6.5) {};
		\node (15) at (5, -5.75) {};
		\node (16) at (7.75, -4.75) {};
		\node (17) at (9, -2.25) {};
		\node (18) at (8, 1.25) {};
		\node (19) at (2.75, 5) {};
		\node (20) at (0.5, 4.5) {};
		\node (21) at (-1, 2.75) {};
		\node (22) at (-1.75, -1.25) {};
		\node (23) at (5.75, 4) {};
		\node (24) at (6.75, 3) {};
		\node (25) at (-2.25, -3.5) {};
		\node (26) at (-1.75, -3.25) {};
		\draw [white, bend right=90, looseness=1.50] (5.center) to (6.center);
		\draw [white, bend right=60] (6.center) to (7.center);
		\draw [white, bend right=45, looseness=0.75] (7.center) to (8.center);
		\draw [white, bend right=60, looseness=1.25] (1.center) to (0.center);
		\draw [white, bend right=45] (2.center) to (1.center);
		\draw [white, bend left=60] (2.center) to (3.center);
		\draw [white, bend left=60, looseness=1.25] (9.center) to (0.center);
		\draw [white, bend left=60] (12.center) to (8.center);
		\draw [white, bend left=75, looseness=1.25] (3.center) to (10.center);
		\draw [white, bend left=60] (10.center) to (11.center);
		\draw [white, bend left=75] (11.center) to (12.center);
		\draw [white, bend left=15] (1.center) to (20.center);
		\draw [white, bend right=45, looseness=0.75] (2.center) to (19.center);
		\draw [white, bend left=15] (3.center) to (23.center);
		\draw [white, bend right] (10.center) to (24.center);
		\draw [white, bend left=45] (11.center) to (18.center);
		\draw [white, bend left] (12.center) to (17.center);
		\draw [white, bend right] (8.center) to (16.center);
		\draw [white, bend right] (15.center) to (7.center);
		\draw [white, bend left=15, looseness=0.75] (14.center) to (6.center);
		\draw [white, bend left=45, looseness=0.75] (13.center) to (5.center);
		\draw [white, bend right] (9.center) to (22.center);
		\draw [white, bend right] (0.center) to (21.center);
		\draw [white, bend right=75, looseness=1.25] (9.center) to (25.center);
		\draw [white, bend right=60, looseness=1.25] (25.center) to (5.center);

    \draw[thick,gray] (-12,-.30) ellipse (0.11 and 0.33);
    \draw[thick,color=blue!5,fill=blue!3] (1.85,-.5) ellipse (4.2 and 4.2);

    \coordinate (A) at (-12,0);
    \coordinate (A') at (3,0);
    \coordinate (A'') at (-5.6,3);

    \draw[thick,red!35] (-12,.36) -- (-8,.3225);
    \draw[thick,red!35] (-12,.34) -- (-8,.3025);
    \draw[thick,red!35] (-8,.3025) -- (0,.2275);
    \draw[thick,red!35] (-12,.32) -- (-8,.2825);
    \draw[thick,red!35] (-8,.2825) -- (2,0.18875);
    \draw[thick,red!35] (-12,.30) -- (-8,.2625);
    \draw[thick,red!35] (-8,0.2625) -- (-2,0.20625);
    \draw[thick,red!35] (-12,.26) -- (-8,.2225);
    \draw[thick,red!35] (-8,.2225) -- (-6,0.20375);
    \draw[thick,red!35] (-12,.28) -- (-8,.2425);
    \draw[thick,red!35] (-8,.2425) -- (-3,0.205);
    \draw[dashed,thick,red!35,-stealth] (2,0.18875) -- (6,0.105);

    \draw [dashed, thin, -stealth, red](-13.150,.36) -- (-12.1,.36);
    \draw [dashed, thin, -stealth, red](-13.045,.26) -- (-12.1,.26);

    \draw[ultra thick,gray!50,ray] (-8,.3225) -- (-5,3.4025);
    \draw[ultra thick,gray!50,ray] (-5,3.4025) -- (-3,1.4025);
    \draw[ultra thick,gray!50,ray] (-3,1.4025) -- (-2.8,-1.7025);
    \draw[ultra thick,gray!50,ray] (-2.8,-1.7025) -- (-12,-.41);
    \draw[ultra thick,gray!50,-stealth] (-2.8,-1.7025) -- (-12,-.41);

    \draw[ultra thick,blue!30,ray] (0,.2275) -- (-0.1,-.525);
    \draw[ultra thick,blue!30,ray] (-0.1,-.525) -- (0.1,-1.3825);
    \draw[ultra thick,blue!30,ray] (0.1,-1.3825) -- (-0.3,-2.825);
    \draw[ultra thick,blue!30,ray] (-0.3,-2.825) -- (-12,-.48);
    \draw[ultra thick,blue!30,-stealth] (-0.3,-2.825) -- (-12,-.48);

    \draw[ultra thick,blue!70,ray] (2,0.18875) -- (3.1,0.75);
    \draw[ultra thick,blue!70,ray] (3.1,0.75) -- (-12,-.27);
    \draw[ultra thick,blue!70,-stealth] (3.1,0.75) -- (-12,-.27);

    \draw[ultra thick,blue!30,ray] (-2,0.20625) -- (1.5,1.25);
    \draw[ultra thick,blue!30,ray] (1.5,1.25) -- (0.75,2);
    \draw[ultra thick,blue!30,ray] (0.75,2) -- (-12,-.20);
    \draw[ultra thick,blue!30,-stealth] (0.75,2) -- (-12,-.20);

    \draw[ultra thick,blue!70,ray] (-3,0.205) -- (2.5,-.62425);
    \draw[ultra thick,blue!70,ray] (2.5,-.62425) -- (-12,-.34);
    \draw[ultra thick,blue!70,-stealth] (2.5,-.62425) -- (-12,-.34);

    \draw[ultra thick,blue!15,ray] (-6,0.20375) -- (-1.45,1.8825);
    \draw[ultra thick,blue!15,ray] (-1.45,1.8825) -- (-0.15,2.9325);
    \draw[ultra thick,blue!15,ray] (-0.15,2.9325) -- (-12,-.13);
    \draw[ultra thick,blue!15,-stealth] (-0.15,2.9325) -- (-12,-.13);


\end{tikzpicture}}}%
    \subfloat[small absorbing region]{\scalebox{0.25}{\begin{tikzpicture}[ray/.style={decoration={markings,mark=at position .5 with {\arrow[>=latex]{>}}},postaction=decorate}]
            \node (0) at (-1.75, 3) {};
		\node (1) at (-0.25, 5) {};
		\node (2) at (2.75, 5.75) {};
		\node (3) at (6, 5) {};
		\node (5) at (-1.25, -5) {};
		\node (6) at (2.75, -7) {};
		\node (7) at (6, -6.75) {};
		\node (8) at (8, -5.25) {};
		\node (9) at (-2.75, -1.25) {};
		\node (10) at (7.75, 3.75) {};
		\node (11) at (8.75, 1.25) {};
		\node (12) at (9.75, -2.25) {};
		\node (13) at (-0.75, -4.5) {};
		\node (14) at (2.75, -6.5) {};
		\node (15) at (5, -5.75) {};
		\node (16) at (7.75, -4.75) {};
		\node (17) at (9, -2.25) {};
		\node (18) at (8, 1.25) {};
		\node (19) at (2.75, 5) {};
		\node (20) at (0.5, 4.5) {};
		\node (21) at (-1, 2.75) {};
		\node (22) at (-1.75, -1.25) {};
		\node (23) at (5.75, 4) {};
		\node (24) at (6.75, 3) {};
		\node (25) at (-2.25, -3.5) {};
		\node (26) at (-1.75, -3.25) {};
		\draw [white, bend right=90, looseness=1.50] (5.center) to (6.center);
		\draw [white, bend right=60] (6.center) to (7.center);
		\draw [white, bend right=45, looseness=0.75] (7.center) to (8.center);
		\draw [white, bend right=60, looseness=1.25] (1.center) to (0.center);
		\draw [white, bend right=45] (2.center) to (1.center);
		\draw [white, bend left=60] (2.center) to (3.center);
		\draw [white, bend left=60, looseness=1.25] (9.center) to (0.center);
		\draw [white, bend left=60] (12.center) to (8.center);
		\draw [white, bend left=75, looseness=1.25] (3.center) to (10.center);
		\draw [white, bend left=60] (10.center) to (11.center);
		\draw [white, bend left=75] (11.center) to (12.center);
		\draw [white, bend left=15] (1.center) to (20.center);
		\draw [white, bend right=45, looseness=0.75] (2.center) to (19.center);
		\draw [white, bend left=15] (3.center) to (23.center);
		\draw [white, bend right] (10.center) to (24.center);
		\draw [white, bend left=45] (11.center) to (18.center);
		\draw [white, bend left] (12.center) to (17.center);
		\draw [white, bend right] (8.center) to (16.center);
		\draw [white, bend right] (15.center) to (7.center);
		\draw [white, bend left=15, looseness=0.75] (14.center) to (6.center);
		\draw [white, bend left=45, looseness=0.75] (13.center) to (5.center);
		\draw [white, bend right] (9.center) to (22.center);
		\draw [white, bend right] (0.center) to (21.center);
		\draw [white, bend right=75, looseness=1.25] (9.center) to (25.center);
		\draw [white, bend right=60, looseness=1.25] (25.center) to (5.center);

    \draw[thick,gray] (-12,-.30) ellipse (0.11 and 0.33);
    \draw[thick,color=blue!5,fill=blue!15] (1.85,-.5) ellipse (1.15 and 1.15);

    \coordinate (A) at (-12,0);
    \coordinate (A') at (3,0);
    \coordinate (A'') at (-5.6,3);

    \draw[thick,red!35] (-12,.36) -- (-8,.3225);
    \draw[thick,red!35] (-12,.34) -- (-8,.3025);
    \draw[thick,red!35] (-8,.3025) -- (0,.2275);
    \draw[thick,red!35] (-12,.32) -- (-8,.2825);
    \draw[thick,red!35] (-8,.2825) -- (2,0.18875);
    \draw[thick,red!35] (-12,.30) -- (-8,.2625);
    \draw[thick,red!35] (-8,0.2625) -- (-2,0.20625);
    \draw[thick,red!35] (-12,.26) -- (-8,.2225);
    \draw[thick,red!35] (-8,.2225) -- (-6,0.20375);
    \draw[thick,red!35] (-12,.28) -- (-8,.2425);
    \draw[thick,red!35] (-8,.2425) -- (-3,0.205);
    \draw[dashed,thick,red!35,-stealth] (2,0.18875) -- (6,0.105);

    \draw [dashed, thin, -stealth, red](-13.150,.36) -- (-12.1,.36);
    \draw [dashed, thin, -stealth, red](-13.045,.26) -- (-12.1,.26);

    \draw[ultra thick,gray!50,ray] (-8,.3225) -- (-5,3.4025);
    \draw[ultra thick,gray!50,ray] (-5,3.4025) -- (-3,1.4025);
    \draw[ultra thick,gray!50,ray] (-3,1.4025) -- (-2.8,-1.7025);
    \draw[ultra thick,gray!50,ray] (-2.8,-1.7025) -- (-12,-.41);
    \draw[ultra thick,gray!50,-stealth] (-2.8,-1.7025) -- (-12,-.41);

    \draw[ultra thick,gray!50,ray] (0,.2275) -- (-0.1,-.525);
    \draw[ultra thick,gray!50,ray] (-0.1,-.525) -- (0.1,-1.3825);
    \draw[ultra thick,gray!50,ray] (0.1,-1.3825) -- (-0.3,-2.825);
    \draw[ultra thick,gray!50,ray] (-0.3,-2.825) -- (-12,-.48);
    \draw[ultra thick,gray!50,-stealth] (-0.3,-2.825) -- (-12,-.48);

    \draw[ultra thick,blue!30,ray] (2,0.18875) -- (3.1,0.75);
    \draw[ultra thick,blue!30,ray] (3.1,0.75) -- (-12,-.27);
    \draw[ultra thick,blue!30,-stealth] (3.1,0.75) -- (-12,-.27);

    \draw[ultra thick,gray!50,ray] (-2,0.20625) -- (1.5,1.25);
    \draw[ultra thick,gray!50,ray] (1.5,1.25) -- (0.75,2);
    \draw[ultra thick,gray!50,ray] (0.75,2) -- (-12,-.20);
    \draw[ultra thick,gray!50,-stealth] (0.75,2) -- (-12,-.20);

    \draw[ultra thick,blue!70,ray] (-3,0.205) -- (2.5,-.62425);
    \draw[ultra thick,blue!70,ray] (2.5,-.62425) -- (-12,-.34);
    \draw[ultra thick,blue!70,-stealth] (2.5,-.62425) -- (-12,-.34);

    \draw[ultra thick,gray!50,ray] (-6,0.20375) -- (-1.45,1.8825);
    \draw[ultra thick,gray!50,ray] (-1.45,1.8825) -- (-0.15,2.9325);
    \draw[ultra thick,gray!50,ray] (-0.15,2.9325) -- (-12,-.13);
    \draw[ultra thick,gray!50,-stealth] (-0.15,2.9325) -- (-12,-.13);


\end{tikzpicture}}}%
    \caption{Figure (a) and (b) show why the absorption in the wide FOV is heavily dependent on ambient scattering. Only blue trajectories are sensitive to the patch of interest. Dark blue paths are strongly affected whereas photons along grey paths have a noise-like effect on the measurement.}%
    \label{fig:wFOV_feature}%
\end{figure}
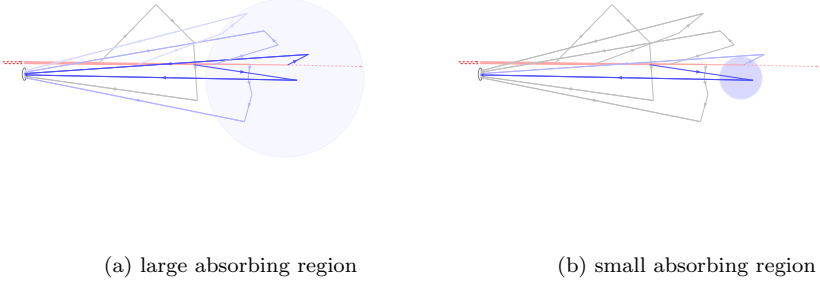

\paragraph{Effect of ambient light} In practice we may also have to consider photons that reach the detector from light sources that are unrelated to our instrument and are out of our control, such as sunlight. Assuming that these sources are approximately constant over a sufficiently large period of time we can think of them as an independent \(\mathsf{Poisson}\) distributed quantity with known intensity \(H_0\) that is added on top of \(\bm{m}\) and \(\bm{n}\). The amount of ambient photons scales with bin width, detector size and width of aperture so this changes the mean of the data from \(H_D \Delta t A_D \mathfrak{m}_{\mathrm{on/off}}\) to \(H_D \Delta t A_D \mathfrak{m}_{\mathrm{on/off}} + \Delta t A_D H_0\) for some \(H_0\) that depends on the FOV, but it doesn't change its underpinning distribution. In that situation the quantities from \Cref{eq:DIAL_fisher_information_single} decrease to
\begin{align}\label{eq:ambient_scaling}
    \left(1+\frac{H_0}{H_D \mathfrak{m}_{\mathrm{on/off}}(v_i,t_j\mid x_D)}\right)^{-1}J_{\mathrm{on/off}}(v_i,t_j\mid x_D)[\theta]
\end{align}
which can be substantial when \(H_0\) is of the same order as the detector response \(H_D\mathfrak{m}_{\mathrm{on/off}}\), or even larger. In reality the factor in \Cref{eq:ambient_scaling} depends on a lot of variables other than the FOV such as the spectral filter or the impulse strength \(H_D\). \(\mathfrak{m}_{\mathrm{on/off}}\) can also be highly inhomogeneous and its magnitude in the relevant data points, i.e. where the plume is located, may scale with the square of its distance from the detector. Quantifying the adversarial effect of an ambient light source is therefore not something that can be done in general but will depend on the use case. Nonetheless, the effects are always worse for larger FOVs because the aperture acts as a filter not only for multiple scattering but also for ambient illumination. Apart from having better resolution for fine image features, this observation is perhaps the primary reason why it is preferable to use multiple FOVs whenever possible. Indeed, if instead of averaging the data we are able to separate single and multiple scattering and the Fisher information of the measurement is uniformly, i.e. for any set of optical parameters and any \(\alpha[\theta]\), larger than that corresponding to just a narrow FOV and single scattering data. This should come as no surprise since with multiple FOV data we can perform any reconstruction procedure that would be possible with a single FOV alone.

\subsection{Weak absorption errors}
In general, neither of the above scenarios are realistic and the true concentration field lies most likely somewhere in between a fully homogeneous and strictly local distribution of gas and for approximately homogeneous \(\sigma_s\) one should only expect improvements by including photons from wider FOVs for fairly large plumes or, equivalently, ones close to the detector. However, if the \(\sigma_s\) is not homogeneous and scattering is caused mostly by particles in the plume, e.g. smoke or water vapour, we end up with a \emph{very} different behaviour in the multiple scattering component. Indeed, assuming that ambient scatterers are rare and concentrated mostly around the plume of shape \(\alpha = \sigma_{a(\mathrm{on})} - \sigma_{a(\mathrm{off})}\), then most photons that reach the detector will have travelled through the absorbing medium of interest and we end up with a situation much more like \Cref{eq:fully_homogeneous_fisher_information}, even in the case of highly localised plumes. Unlike with the two trivial edge cases considered in \cref{ssec:trivial}, for non-trivial images it is not possible to eliminate the scattering parameters from the quantities that can be used in order to quantify the uncertainty/errors of the reconstructed absorption parameters. In the following we will study a more relevant scenario where the plume shape is known and the concentration is low. Although not entirely unrealistic this may seem rather restrictive at first but we will later discuss how this transfers to considerably more general plumes of unknown shape. We start with the following simple result which relates \(\mathsf{Poisson}\) and \(\mathsf{Normal}\) random variables when the \(\mathsf{Poisson}\) intensities are relatively large. It further establishes a natural connection with our semi-parametric formulation as we explicitly consider the differential absorption in the form of \(\frac{\bm{n}_{v_i, t_j}}{\bm{m}_{v_i, t_j}}\).

\begin{lemma}\label{lem:normal_approximation}
    Let \(\bm{m}_{v_i, t_j}\) and \(\bm{n}_{v_i, t_j}\) be as in \cref{eq:data_noise_model} and define for each direction \(v_i\) and time bin corresponding to mid-point \(t_j\) the random variables \(\bm{\bm{y}_{v_i, t_j}}\) and \(\bm{z}_{v_i, t_j}\) via
    \begin{align}\label{eq:quotient_transform}
    \begin{split}
        \bm{y}_{v_i, t_j} &:= \log \left(\frac{\bm{n}_{v_i, t_j}}{\bm{m}_{v_i, t_j}} \right) \\
        \bm{z}_{v_i, t_j} &:= \frac{\bm{n}_{v_i, t_j} + \bm{m}_{v_i, t_j}}{2}.
    \end{split}
    \end{align}
    Further assume that the differential absorption is considerably smaller than the signal, i.e. that we have
    \begin{align*}
        \frac{\mathfrak{m}_{\mathrm{off}} (v_i, t_j\mid x_D) - \mathfrak{m}_{\mathrm{on}} (v_i, t_j\mid x_D)}{\mathfrak{m}_{\mathrm{off}} (v_i, t_j\mid x_D)} \ll 1 
    \end{align*}
    as well as
    \begin{align*}
        \Delta t H_D A_D \mathfrak{m}_{\mathrm{on}} (v_i, t_j\mid x_D) \gg 1.
    \end{align*}
    Then \(\bm{y}_{v_i, t_j}\) as well as \(\bm{z}_{v_i, t_j}\) are mutually independent and their distribution can be approximated as 
    \begin{align}\label{eq:normal_approximation}
    \begin{split}
        \bm{y}_{v_i, t_j} &\sim \mathsf{Normal}\left(\frac{\mathfrak{m}_{\mathrm{off}} (v_i, t_j\mid x_D)}{\mathfrak{m}_{\mathrm{on}} (v_i, t_j\mid x_D)} - 1,\frac{2}{\Delta t H_D A_D \mathfrak{m}_{\mathrm{off}} (v_i, t_j\mid x_D)}\right)\\
        \bm{z}_{v_i, t_j} &\sim \mathsf{Normal}\left(\Delta t H_D A_D \mathfrak{m}_{\mathrm{off}} (v_i, t_j\mid x_D),\frac{\Delta t H_D A_D \mathfrak{m}_{\mathrm{off}} (v_i, t_j\mid x_D)}{2}\right).
    \end{split}
    \end{align}
    \begin{proof}
Note that \(\bm{y}_{v_i, t_j}\) and \(\bm{z}_{v_i, t_j}\) are not independent but the arrays \(\bm{y}\) and \(\bm{z}\) created from objects as in \cref{eq:quotient_transform} have independent entries due to the assumed independence of the data in each bin and direction. The normal approximation for the \(\bm{z}_{v_i, t_j}\) is standard while the result for \(\bm{y}_{v_i, t_j}\) requires \cite{DazFrancs:2013} after a Taylor approximation of the logarithm \(\log(x) \approx x-1\) for \(x \approx 1\). Note that the effect of (differential) absorption is neglected in the variance terms.
\end{proof}
\end{lemma}

It is easily seen that \(\bm{y}\) is essentially a discretised and noisy version of the same quantity (subject to a bijective transform) as in \cref{thm:relaxed_uniqueness} which considered its dependence on the dispersion parameters. Recall that due to \cref{eq:on_data} and \cref{eq:off_data} for single scattering we have 
\begin{align*}
    \frac{\mathfrak{m}_{\mathrm{off}} (v_i, t_j\mid x_D)}{\mathfrak{m}_{\mathrm{on}} (v_i, t_j\mid x_D)} - 1 = e^{2 \int_{0}^{t/2}\alpha[\theta](x_D+rv)dr} - 1 \approx 2 \int_{0}^{t/2}\alpha[\theta](x_D+rv)dr
\end{align*}
where the final approximation is due to \(\alpha\) being small. For the remainder of this section we assume that \(\alpha\) is known up to a constant, i.e. we have access to \(\alpha_0\) such that \(\alpha[\theta](x) = \theta_0 \alpha_0(x)\), where \(\theta_0\) denotes the first component of \(\theta \in \Theta \subseteq [0,\infty) \times \mathbb{R}^{d}\), for each \(x \in X\). This is identical to the cases discussed in the previous section except that the plume may now take an arbitrary but known shape given by \(\alpha_0\) and that we consider a \(d\)-dimensional nuisance component which controls the scattering behaviour. In this case we can write
\begin{align*}
    \frac{\mathfrak{m}_{\mathrm{off}} (v_i, t_j\mid x_D)}{\mathfrak{m}_{\mathrm{on}} (v_i, t_j\mid x_D)} - 1 &= \int_{\ell_{00}(X)} e^{\int_{\Gamma(\xi)}  \alpha[\theta](z)dz} Q_{v_i, t_j}[\theta](d \xi) - 1 \\
    & \approx \theta_0 \int_{\ell_{00}(X)}\int_{\Gamma(\xi)}   \alpha_0(z)dz Q_{v_i, t_j}[\theta] (d\xi)
\end{align*}
where \(Q_{v_i, t_j}\) is a probability distribution over the set of finite sequences \(\ell_{00}(X)\) of points in \(X\). In particular, each \(\xi\) can be represented as an ordered set of points \(\xi_0, \dots, \xi_{k+1}\) meaning that \(Q_{v_i, t_j}\) is effectively a probability distribution over piece-wise linear paths which can be described as defined in \cref{eq:path_definition} of the appendix. Note that in the case of a narrow FOV the distribution is a delta concentrated on the unique path that corresponds to a single scattering event for a fixed direction and time. As such the single scattering absorption is independent of the nuisance component in \(\theta\) whereas in the case of wider FOVs the resulting \(Q_{v_i, t_j}\) depends on the nuisance components through the scattering parameters. This naturally has an effect on our ability to recover the parameter of interest \(\theta_0\). For brevity we will introduce
\begin{align}
    \bar{L}_{i,j}[\theta] &= \frac{\Delta t H_D A_D \mathfrak{m}_{\mathrm{off}} (v_i, t_j\mid x_D)}{2} \label{eq:Z_IJ}\\
    \bar{q}_{i,j}[\theta] &= \int_{\ell_{00}(X)}\int_{\Gamma(\xi)}   \alpha_0(z)dz Q_{v_i, t_j}[\theta] (d\xi) \label{eq:Q_IJ} .
\end{align}
and note that even though \(\theta\) is unknown, the quantity \(\bar{L}_{i,j}[\theta]\) can be estimated accurately from \(\bm{z}_{v_i,t_j}\) and thus (under the conditions of \cref{lem:normal_approximation}) be assumed as given/observed. We have the following.
\begin{theorem}[Optimal detection of low concentrations]\label{thm:UMP_test}
    Let the random variables \(\bm{\bm{y}}_{v_i, t_j}\) be as in \cref{lem:normal_approximation} and \(\bar{L}_{i,j}[\theta]\) and \(\bar{q}_{i,j}[\theta]\) as in \cref{eq:Z_IJ} and \cref{eq:Q_IJ}. Further let \(\hat{q}_{i,j} > 0\) be an arbitrary collection of positive scalars and that \(\bar{L}_{i,j}[\theta]\) is known for any \(\theta \in \Theta\). If we define the set
    \begin{align*}
        \Theta_0 = \left\{\theta \in \Theta : \theta_0 = 0\right\}
    \end{align*}
    then we have for all \(\theta \in \Theta\)
    \begin{align}\label{eq:test_statistics}
        \frac{\sum_{i,j} \hat{q}_{i,j} \bar{L}_{i,j}[\theta] \bm{y}_{v_i,t_j} - \theta_0 \hat{q}_{i,j} \bar{q}_{i,j}[\theta] \bar{L}_{i,j}[\theta]}{\sqrt{\sum_{i,j} \hat{q}_{i,j}^2 \bar{L}_{i,j}[\theta]}} \sim \mathsf{Normal}(0,1).
    \end{align}
    In particular, \cref{eq:test_statistics} can be used to test the hypothesis \(H_0: \theta \in \Theta_0\) vs. \(H_1 : \theta \in \Theta \setminus \Theta_0\) where the rejection regions take the form
    \begin{align}\label{eq:rejection_region}
        \frac{\sum_{i,j} \hat{q}_{i,j} \bar{L}_{i,j}[\theta] \bm{y}_{v_i,t_j}}{ \sqrt{\sum_{i,j} \hat{q}_{i,j}^2 \bar{L}_{i,j}[\theta]}} \geq R > 0
    \end{align}
    with \(R\) is a quantile of the standard normal distribution that depends only on the significance level. The power of the test is given by 
    \begin{align}\label{eq:test_power}
        \mathbb{P}\left(\mathsf{Normal}(0,1) > R - \theta_0 \frac{\sum_{i,j} \hat{q}_{i,j} \bar{q}_{i,j}[\theta] \bar{L}_{i,j}[\theta]}{ \sqrt{\sum_{i,j} \hat{q}_{i,j}^2 \bar{L}_{i,j}[\theta]}} \right).
    \end{align}
    and is maximised for \(\bar{q}_{i,j}[\theta] = \hat{q}_{i,j}\) which is independent of \(\theta\) in the case of narrow FOV measurements and yields the uniformly most powerful test for the above hypothesis. 
    \begin{proof}
    The only non-trivial claim is that the test for narrow FOV data is uniformly most powerful, the rest is essentially an immediate consequence of \cref{lem:normal_approximation} following simple rearrangements. Note that in the case of single scattering the quantities in \Cref{eq:Z_IJ} as well as \cref{eq:Q_IJ} are assumed known (as \(q_{i,j}\) doesn't depend on the scattering). The probability density for the sample \(\bm{y}_{v_i,t_j}\) is given by
    \begin{align*}
        \prod_{i,j} \sqrt{\frac{L_{i,j}}{2 \pi}} \exp\left(-\frac{1}{2} \sum_{i,j} \bm{y}_{v_i,t_j} L_{i,j} \bm{y}_{v_i,t_j}+\theta_0^2q_{i,j} L_{i,j} q_{i,j} \right) 
        \exp\left(\theta_0 \sum_{i,j} \bm{y}_{v_i,t_j} L_{i,j} q_{i,j} \right)
    \end{align*}
    from which the factorisation theorem implies that \(T(\bm{y}):=\sum_{i,j} \bm{y}_{v_i,t_j} L_{i,j} q_{i,j}\) is a sufficient statistic for \(\theta_0\) (see Theorem 6.2.6 in \cite{CasellaBerger:2001}) and the likelihood ratio is an increasing function of \(T(\bm{y})\).
    The optimality of the rejection regions as in \Cref{eq:rejection_region} is a direct consequence of the Karlin-Rubin Theorem (see e.g. Theorem 8.3.17 in \cite{CasellaBerger:2001}) because \(\theta_0 \geq 0\).
    \end{proof}
\end{theorem}
It is worth noticing that the procedure presented in \cref{thm:UMP_test} is in a rather strong sense optimal when it comes to detecting low concentration plumes with narrow FOV data. As such, an improvement with respect to that particular task is achieved as soon as we can construct a test from wide FOV data that outperforms its narrow FOV equivalent. It is worth noticing that the power is essentially a function of 
\begin{align*}
    \Pi_L(\bar{q}[\theta] \to  \hat{q} ):=\frac{\sum_{i,j} \hat{q}_{i,j} \bar{q}_{i,j}[\theta] \bar{L}_{i,j}[\theta]}{ \sqrt{\sum_{i,j} \hat{q}_{i,j}^2 \bar{L}_{i,j}[\theta]}}
\end{align*}
which is nothing but a scaled inner product of \(\hat{q}\) with \(\bar{q}[\theta]\) weighted by \(\bar{L}_{i,j}\) or, more precisely, the scalar projection of \(\bar{q}[\theta]\) onto \(\hat{q}\) w.r.t. that inner product. Typically we won't have access to the scattering parameters but we may select an arbitrary \(\hat{\theta} \in \Theta\) and set \(\hat{q}_{i,j} = \bar{q}_{i,j}[\hat{\theta}]\) resulting in a power no worse than had we chosen
\begin{align}\label{eq:worst_case_power}
    \hat{\Pi}_{\Theta, L}[\theta] := \inf_{\hat{q} \in \bar{q}[\Theta]} \Pi_L(\bar{q}[\theta] \to  \hat{q} ).
\end{align}
Note that \(\theta\) can be seen as a ground truth parameter whereas \(\Theta\) is a region determined by prior knowledge about the true parameter. The quantity subject to minimisation in \Cref{eq:worst_case_power} depends on the FOV through \(\bar{q}\) and becomes trivial when the FOV is assumed to collect only single scattering, thus light with known trajectories. To make the distinction more clear we will refer to \(\hat{\Pi}_{\Theta, L}[\theta]\) only for a given/fixed wide FOV and define
\begin{align*}
    \Pi^{0}_L[\theta] := \sqrt{\sum_{i,j} q^0_{i,j} L^{0}_{i,j}[\theta] q^0_{i,j} }
\end{align*}
where \(q^0_{i,j}\) and \(L^{0}_{i,j}[\theta]\) are as in \Cref{eq:Z_IJ} and \cref{eq:Q_IJ} with only \(L^0\) depending on \(\theta\). By studying the ratio
\begin{align}
    \Phi_\Theta[\theta] := \frac{\hat{\Pi}_{\Theta, L}[\theta]}{\Pi^{0}_L[\theta]} - 1
\end{align}
for varying ground truth parameters \(\theta\) and constraints on \(\Theta\) we can gain insight into the relative performance of wide and narrow FOV by means of a simple scalar indicator without having to assume known scattering parameters. It should be noted that the insight is primarily of qualitative nature as the signal strength \(L\) for wide and narrow FOVs, in particular the ratio nFOV:wFOV, has a major impact on the behaviour of \(\Phi\) whose precise quantitative characteristics therefore depend on unknown quantities. 

In general we will have parameters that result in the minimisation in \Cref{eq:worst_case_power} having intractable integrals that cannot be dealt with analytically. The following plots show a simple case where \(\alpha_0\) is determined by a single kernel with radius approximately 20m which placed at a distance of 100m away from the detector. We consider \(\theta = (\theta_0, g[\theta], \sigma_0[\theta], \sigma_s[\theta])\), i.e. the nuisance component of \(\theta\) has three components. An ambient scattering intensity (uniformly distributed scattering particles) denoted by \(\sigma_0[\theta]\) measured from detector to the mid point of the plume kernel, a kernel scattering intensity (scattering particles aligned with the plume) denoted by \(\sigma_s[\theta]\) measured as its total thickness and a Henyey-Greenstein parameter \(g[\theta]\) which controls the phase function and is constrained to be between \(0\) and \(0.7\). Unlike the case depicted in \cref{fig:wFOV_feature} we kept the size of the kernel fixed. However, an increased presence of ambient scattering particles shifts the distribution of trajectories along which photons are likely to be observed towards the scenario in \cref{fig:wFOV_feature} (b).  
\begin{figure}[ht]%
    \centering
    \subfloat[\centering unconstrained]{{\includegraphics[scale=0.45,trim={0cm 0cm 1.5cm 0cm},clip]{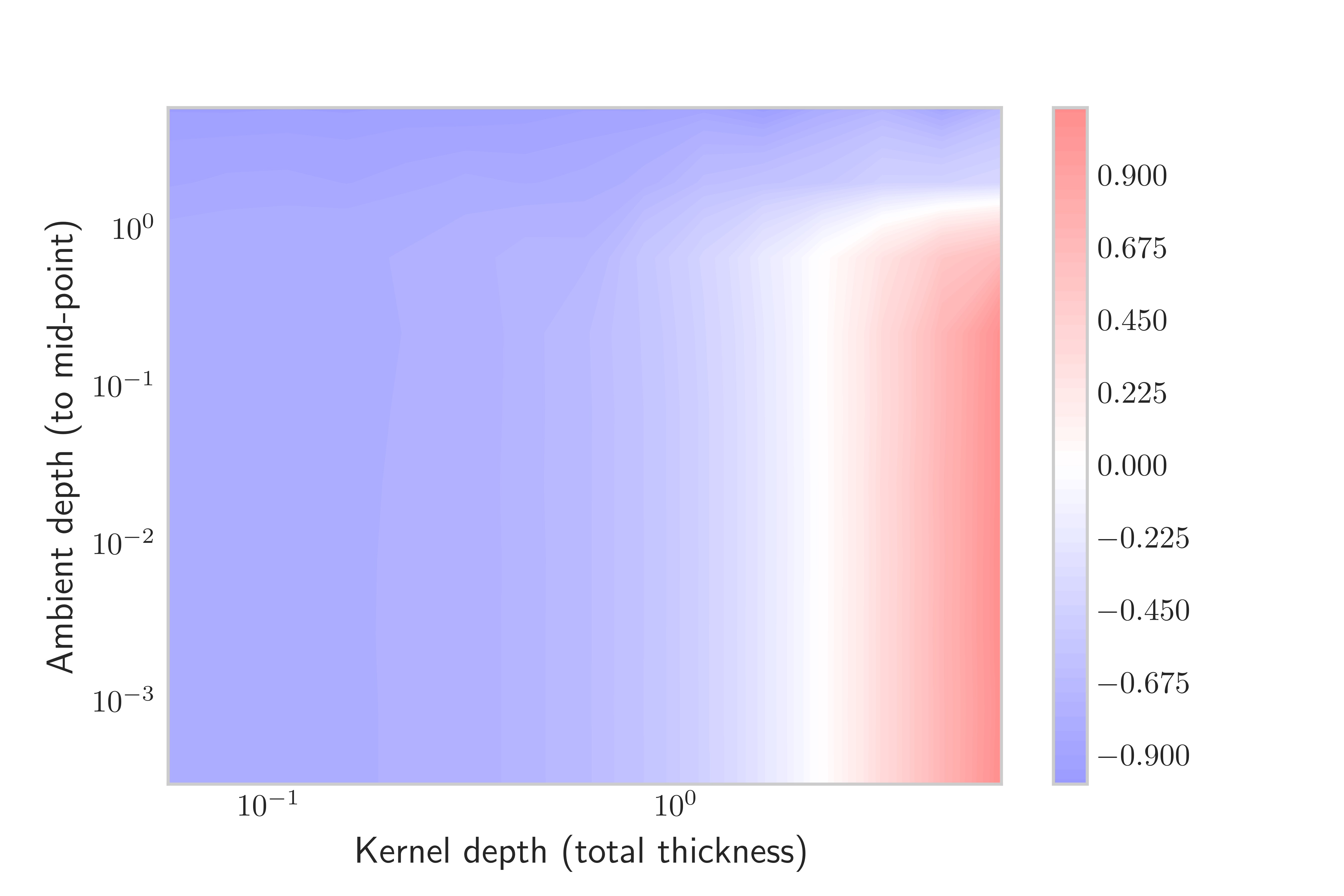} }}%
    \subfloat[\centering known phase function]{{\includegraphics[scale=0.45,trim={0cm 0cm 1.5cm 0cm},clip]{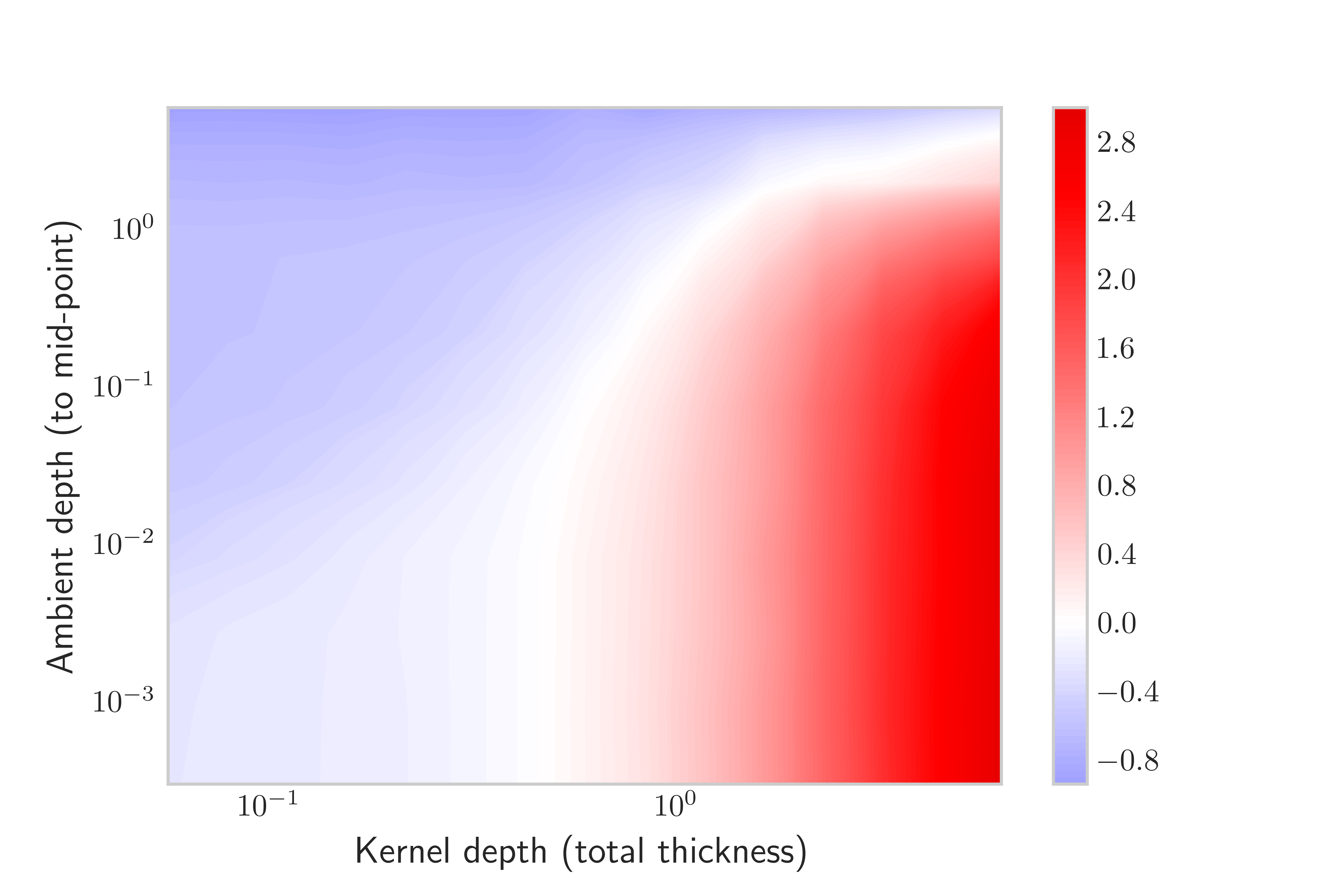} }}\\%
    \subfloat[\centering known optical depth]{{\includegraphics[scale=0.45,trim={0cm 0cm 1.5cm 0cm},clip]{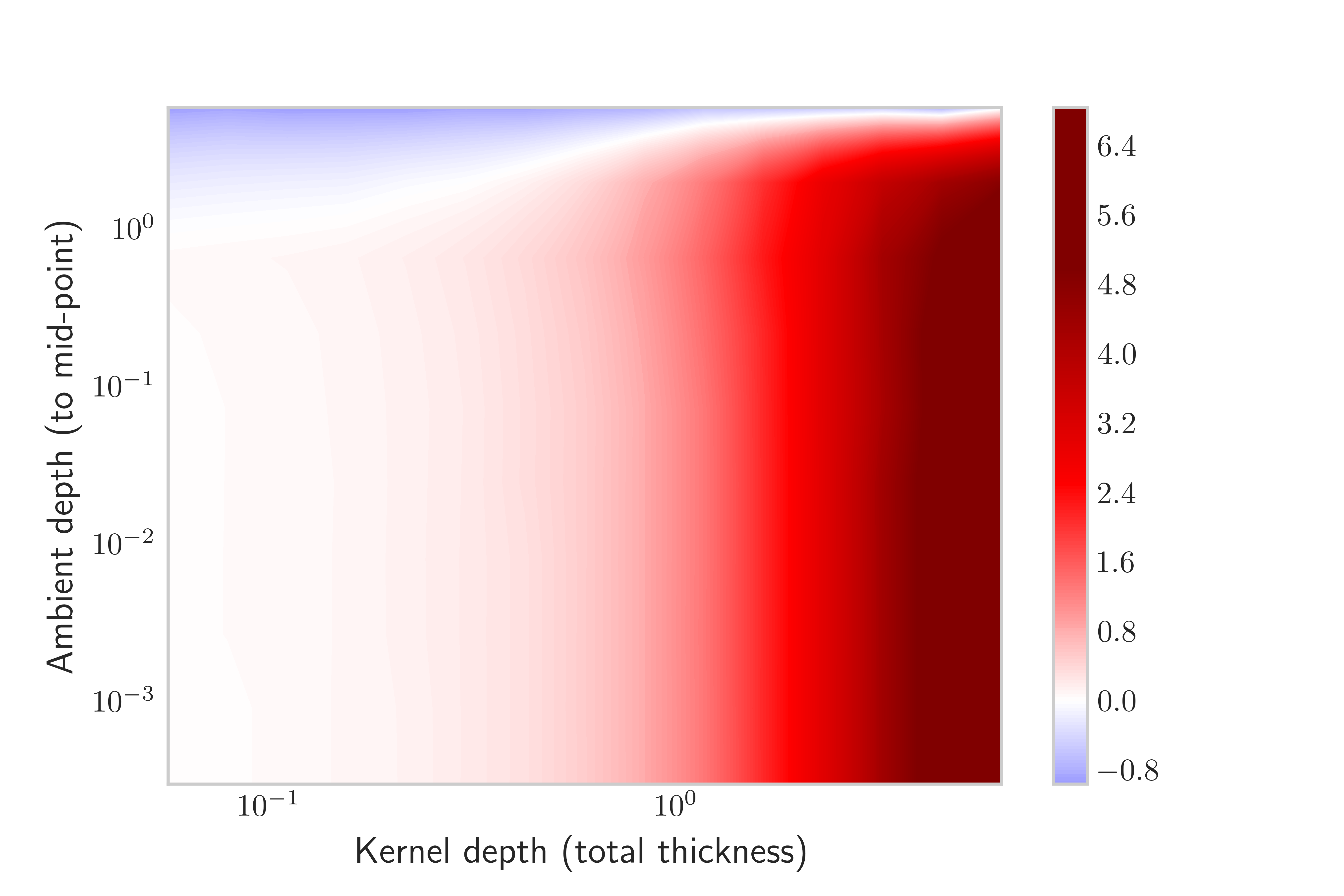} }}%
    \subfloat[\centering full knowledge]{{\includegraphics[scale=0.45,trim={0cm 0cm 1.5cm 0cm},clip]{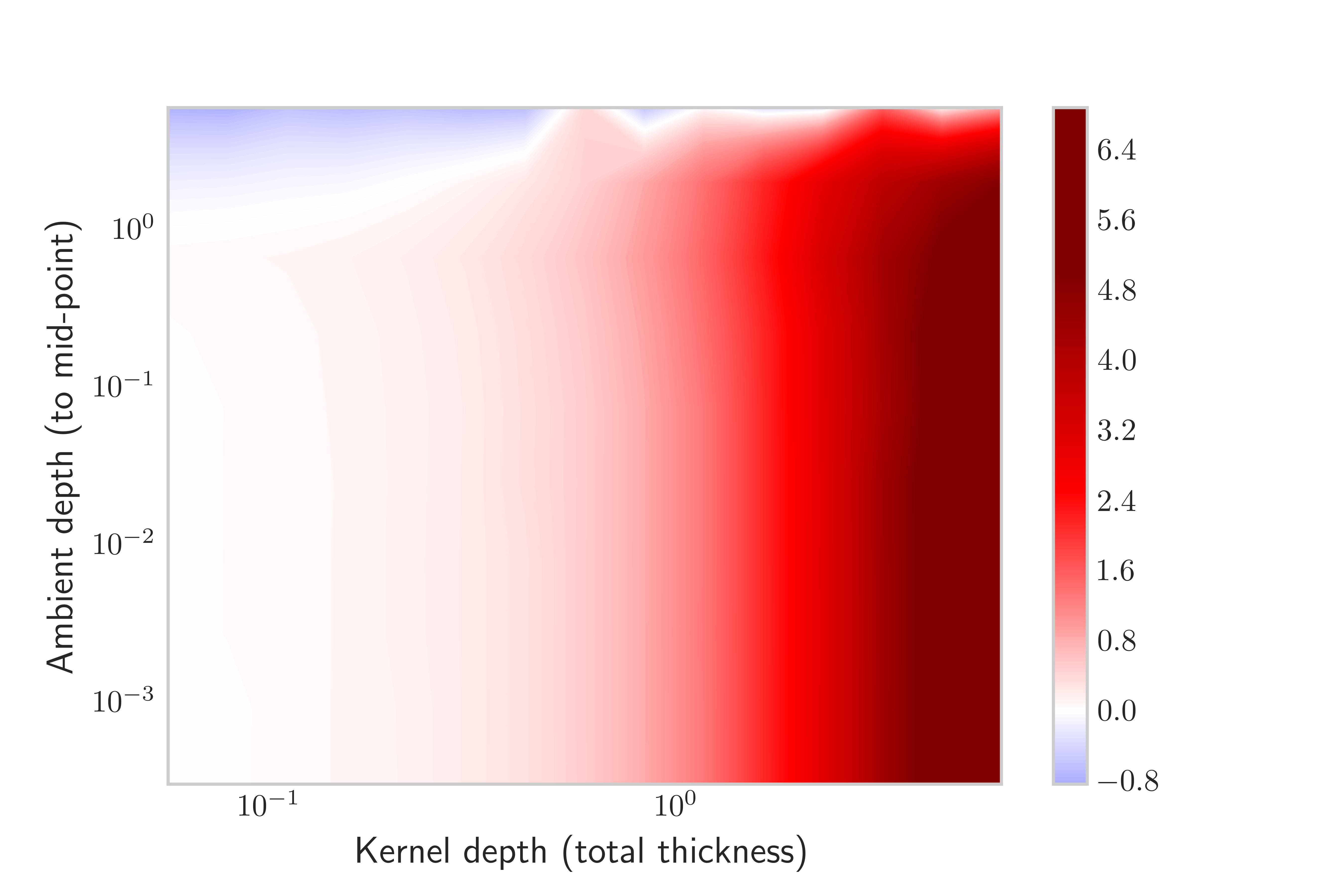} }}%
    \caption{The red regions indicate better detection wide wide FOV data. Figure (a) shows that if nothing other than loose the bounds for the distribution of ambient particles is known, then wide FOVs will improve the reconstruction only for optically thick plumes. Figures (b), where the optical depth is assumed unknown while \(g[\theta]\) is fixed to the correct value, and (c), where only the phase function is unknown while \(\sigma_0[\theta]\) and \(\sigma_s[\theta]\) are fixed, show that knowledge regarding the distribution of scattering particles is considerably more valuable than knowledge of the phase function. The former is in this example virtually equivalent to full knowledge of the scattering particles, i.e. \(\Theta = \{\theta\}\), which is shown in Figure (d).}%
    \label{fig:free_test}%
\end{figure}

The graphic in \cref{fig:free_test}, which shows cross sections of \(\theta \mapsto \Phi_{\Theta,L}[\theta]\) at \(g[\theta] = 0\) for different choices \(\Theta\). Note that we fixed the phase function in the ground truth so the plots become 2-dimensional, homogeneous scattering was chose as it maximises the single scattering component in our case, and the behaviour is similar for other cross sections. The phase function was left variable in \(\Theta\) unless stated otherwise. The plots, suggest that the quantities that matter for good ``worst-case alignment`` are precisely those that are considered as reconstructable in \cref{thm:relaxed_uniqueness}. 
\begin{figure}[ht]%
    \centering
    \subfloat[\centering limited ambient scattering]{{\includegraphics[scale=0.45,trim={0cm 0cm 1.5cm 0cm},clip]{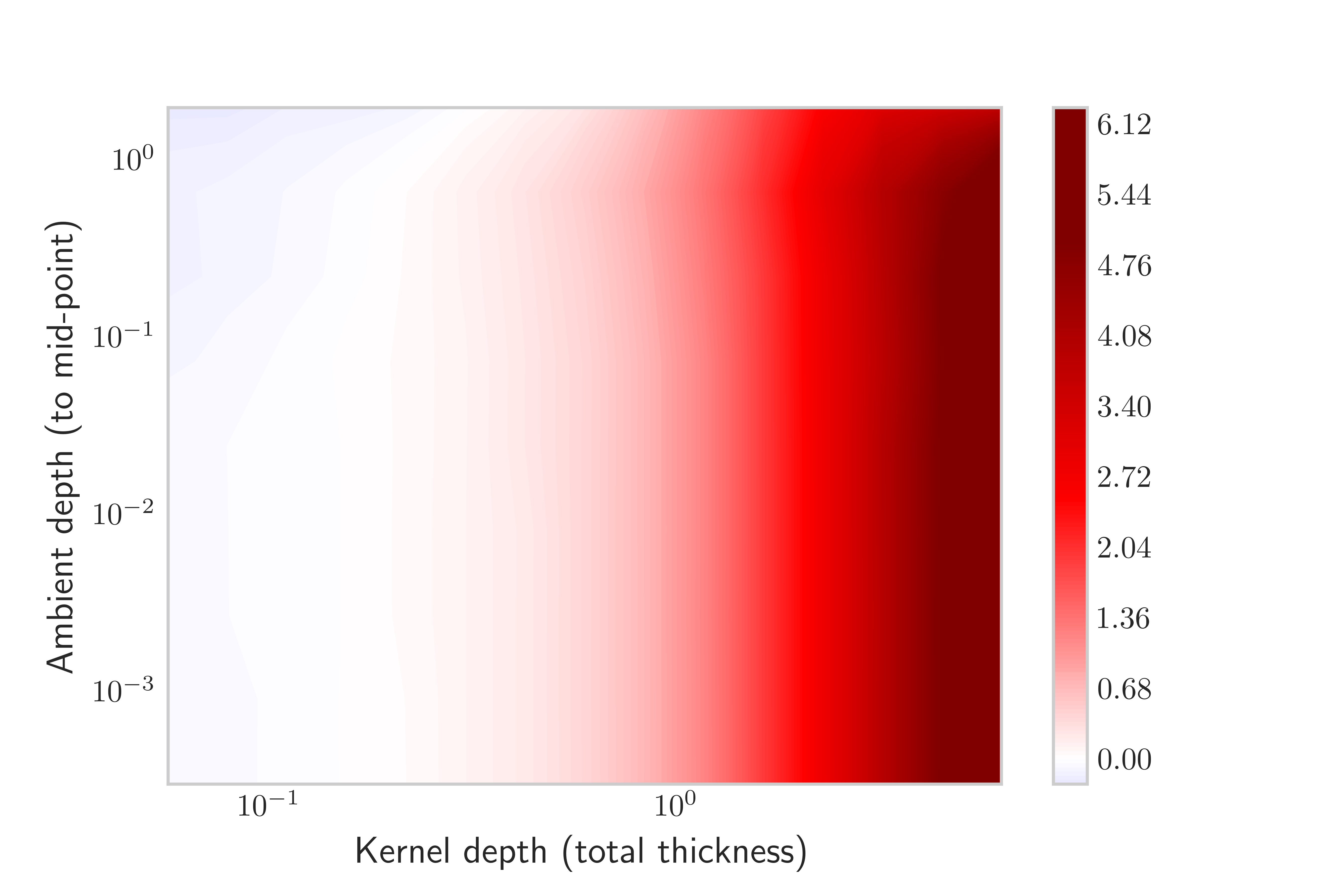} }}%
    \subfloat[\centering full knowledge]{{\includegraphics[scale=0.45,trim={0cm 0cm 1.5cm 0cm},clip]{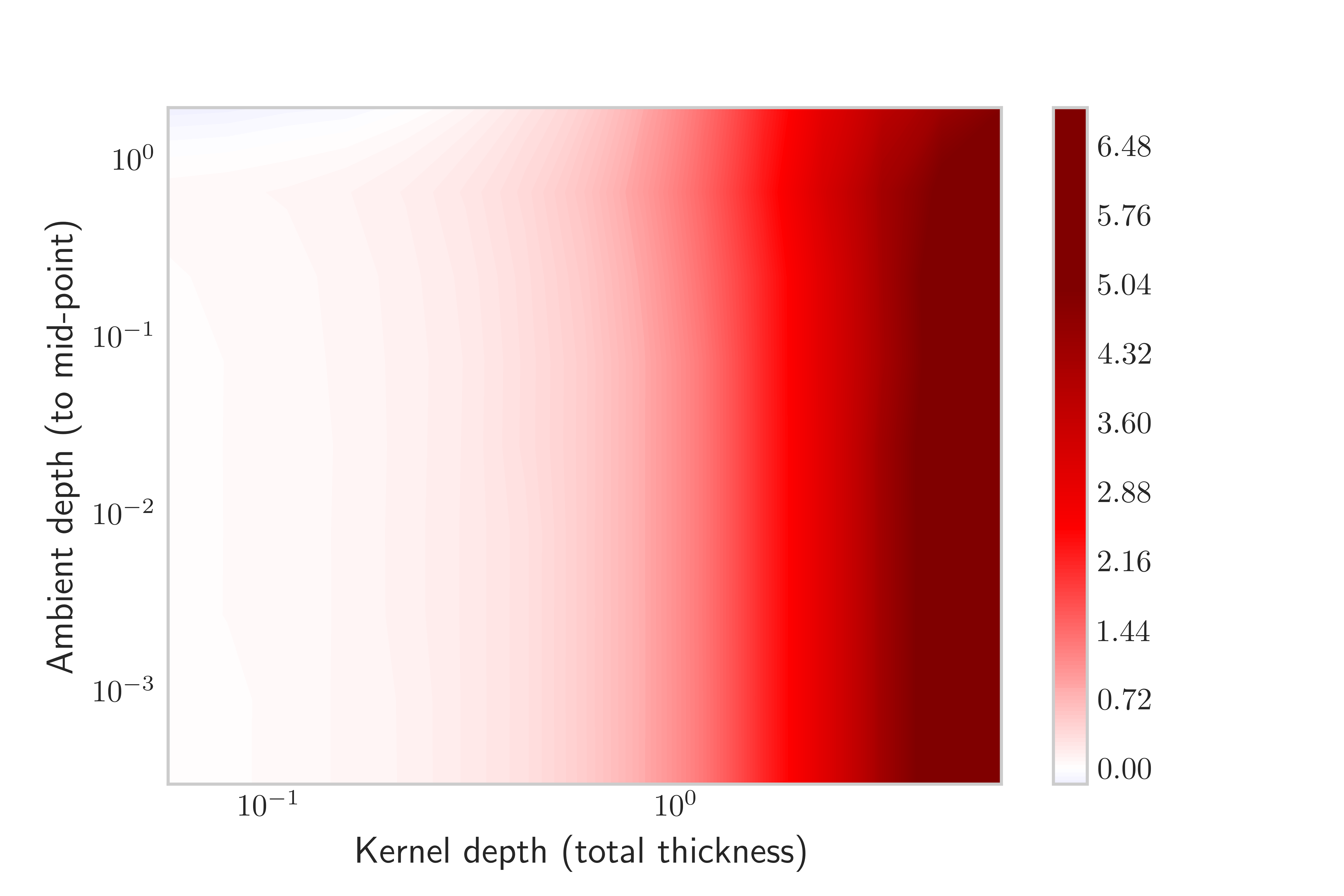} }}%
    \caption{The red regions indicate better detection wide FOV data. Figure (a) and (b) show very similar results indicating that full knowledge of the scattering parameters does not yield considerably better results than a merely a limit on ambient scattering corresponding to approximately \(\sigma_0[\theta] \lesssim 2\).}%
    \label{fig:known_test}%
\end{figure}

\Cref{fig:free_test} and \cref{fig:known_test} clearly indicate that amongst the optical nuisance parameters the distribution of scattering particles has the biggest effect on the quality of our measurement. Indeed, as \cref{fig:known_test} shows even knowledge of the complete set of scattering parameters does not provide much better alignment w.r.t. \(\Pi_L\) as a simple constraint that limits the amount of ambient scattering particles. Although we only considered a single parameter and a simplified scenario, the same idea extends to more general \(\theta\) which is to say that photons from wider FOVs are most useful when we are interested in aspects of the differential concentration field \(\alpha\) that are relatively well aligned with \(\sigma_s\), i.e not too small or far away from most of the scattering particles. Indeed, the simple scenario can more generally be thought of as an indication of how much information wider FOVs preserve in the presence of unknown scattering parameters with regard to a single component of an element from a good kernel space. 
\begin{figure}[h]%
    \centering
    \subfloat[large feature]{\scalebox{0.25}{\begin{tikzpicture}[ray/.style={decoration={markings,mark=at position .5 with {\arrow[>=latex]{>}}},postaction=decorate}]

		\node (0) at (-1.75, 3) {};
		\node (1) at (-0.25, 5) {};
		\node (2) at (2.75, 5.75) {};
		\node (3) at (6, 5) {};
		\node (5) at (-1.25, -5) {};
		\node (6) at (2.75, -7) {};
		\node (7) at (6, -6.75) {};
		\node (8) at (8, -5.25) {};
		\node (9) at (-2.75, -1.25) {};
		\node (10) at (7.75, 3.75) {};
		\node (11) at (8.75, 1.25) {};
		\node (12) at (9.75, -2.25) {};
		\node (13) at (-0.75, -4.5) {};
		\node (14) at (2.75, -6.5) {};
		\node (15) at (5, -5.75) {};
		\node (16) at (7.75, -4.75) {};
		\node (17) at (9, -2.25) {};
		\node (18) at (8, 1.25) {};
		\node (19) at (2.75, 5) {};
		\node (20) at (0.5, 4.5) {};
		\node (21) at (-1, 2.75) {};
		\node (22) at (-1.75, -1.25) {};
		\node (23) at (5.75, 4) {};
		\node (24) at (6.75, 3) {};
		\node (25) at (-2.25, -3.5) {};
		\node (26) at (-1.75, -3.25) {};
		\draw [semithick, dashed, bend right=90, looseness=1.50] (5.center) to (6.center);
		\draw [semithick, dashed, bend right=60] (6.center) to (7.center);
		\draw [semithick, dashed, bend right=45, looseness=0.75] (7.center) to (8.center);
		\draw [semithick, dashed, bend right=60, looseness=1.25] (1.center) to (0.center);
		\draw [semithick, dashed, bend right=45] (2.center) to (1.center);
		\draw [semithick, dashed, bend left=60] (2.center) to (3.center);
		\draw [semithick, dashed, bend left=60, looseness=1.25] (9.center) to (0.center);
		\draw [semithick, dashed, bend left=60] (12.center) to (8.center);
		\draw [semithick, dashed, bend left=75, looseness=1.25] (3.center) to (10.center);
		\draw [semithick, dashed, bend left=60] (10.center) to (11.center);
		\draw [semithick, dashed, bend left=75] (11.center) to (12.center);
		\draw [semithick, dashed, bend left=15] (1.center) to (20.center);
		\draw [semithick, dashed, bend right=45, looseness=0.75] (2.center) to (19.center);
		\draw [semithick, dashed, bend left=15] (3.center) to (23.center);
		\draw [semithick, dashed, bend right] (10.center) to (24.center);
		\draw [semithick, dashed, bend left=45] (11.center) to (18.center);
		\draw [semithick, dashed, bend left] (12.center) to (17.center);
		\draw [semithick, dashed, bend right] (8.center) to (16.center);
		\draw [semithick, dashed, bend right] (15.center) to (7.center);
		\draw [semithick, dashed, bend left=15, looseness=0.75] (14.center) to (6.center);
		\draw [semithick, dashed, bend left=45, looseness=0.75] (13.center) to (5.center);
		\draw [semithick, dashed, bend right] (9.center) to (22.center);
		\draw [semithick, dashed, bend right] (0.center) to (21.center);
		\draw [semithick, dashed, bend right=75, looseness=1.25] (9.center) to (25.center);
		\draw [semithick, dashed, bend right=60, looseness=1.25] (25.center) to (5.center);

    \draw[thick,gray] (-12,-.30) ellipse (0.11 and 0.33);
    \draw[thick,color=blue!5,fill=blue!3] (1.85,-.5) ellipse (4.2 and 4.2);

    \coordinate (A) at (-12,0);
    \coordinate (A') at (3,0);
    \coordinate (A'') at (-5.6,3);

    \draw[thick,red!35] (-12,.36) -- (-8,.3225);
    \draw[thick,red!35] (-12,.34) -- (-8,.3025);
    \draw[thick,red!35] (-8,.3025) -- (0,.2275);
    \draw[thick,red!35] (-12,.32) -- (-8,.2825);
    \draw[thick,red!35] (-8,.2825) -- (2,0.18875);
    \draw[thick,red!35] (-12,.30) -- (-8,.2625);
    \draw[thick,red!35] (-8,0.2625) -- (-2,0.20625);
    \draw[thick,red!35] (-12,.26) -- (-8,.2225);
    \draw[thick,red!35] (-8,.2225) -- (-6,0.20375);
    \draw[thick,red!35] (-12,.28) -- (-8,.2425);
    \draw[thick,red!35] (-8,.2425) -- (-3,0.205);
    \draw[dashed,thick,red!35,-stealth] (2,0.18875) -- (6,0.105);

    \draw [dashed, thin, -stealth, red](-13.150,.36) -- (-12.1,.36);
    \draw [dashed, thin, -stealth, red](-13.045,.26) -- (-12.1,.26);

    \draw[ultra thick,gray!50,ray] (-8,.3225) -- (-5,3.4025);
    \draw[ultra thick,gray!50,ray] (-5,3.4025) -- (-3,1.4025);
    \draw[ultra thick,gray!50,ray] (-3,1.4025) -- (-2.8,-1.7025);
    \draw[ultra thick,gray!50,ray] (-2.8,-1.7025) -- (-12,-.41);
    \draw[ultra thick,gray!50,-stealth] (-2.8,-1.7025) -- (-12,-.41);

    \draw[ultra thick,blue!30,ray] (0,.2275) -- (-0.1,-.525);
    \draw[ultra thick,blue!30,ray] (-0.1,-.525) -- (0.1,-1.3825);
    \draw[ultra thick,blue!30,ray] (0.1,-1.3825) -- (-0.3,-2.825);
    \draw[ultra thick,blue!30,ray] (-0.3,-2.825) -- (-12,-.48);
    \draw[ultra thick,blue!30,-stealth] (-0.3,-2.825) -- (-12,-.48);

    \draw[ultra thick,blue!70,ray] (2,0.18875) -- (3.1,0.75);
    \draw[ultra thick,blue!70,ray] (3.1,0.75) -- (-12,-.27);
    \draw[ultra thick,blue!70,-stealth] (3.1,0.75) -- (-12,-.27);

    \draw[ultra thick,blue!30,ray] (-2,0.20625) -- (1.5,1.25);
    \draw[ultra thick,blue!30,ray] (1.5,1.25) -- (0.75,2);
    \draw[ultra thick,blue!30,ray] (0.75,2) -- (-12,-.20);
    \draw[ultra thick,blue!30,-stealth] (0.75,2) -- (-12,-.20);

    \draw[ultra thick,blue!70,ray] (-3,0.205) -- (2.5,-.62425);
    \draw[ultra thick,blue!70,ray] (2.5,-.62425) -- (-12,-.34);
    \draw[ultra thick,blue!70,-stealth] (2.5,-.62425) -- (-12,-.34);

    \draw[ultra thick,blue!15,ray] (-6,0.20375) -- (-1.45,1.8825);
    \draw[ultra thick,blue!15,ray] (-1.45,1.8825) -- (-0.15,2.9325);
    \draw[ultra thick,blue!15,ray] (-0.15,2.9325) -- (-12,-.13);
    \draw[ultra thick,blue!15,-stealth] (-0.15,2.9325) -- (-12,-.13);


\end{tikzpicture}}}%
    \subfloat[small feature]{\scalebox{0.25}{\begin{tikzpicture}[ray/.style={decoration={markings,mark=at position .5 with {\arrow[>=latex]{>}}},postaction=decorate}]

		\node (0) at (-1.75, 3) {};
		\node (1) at (-0.25, 5) {};
		\node (2) at (2.75, 5.75) {};
		\node (3) at (6, 5) {};
		\node (5) at (-1.25, -5) {};
		\node (6) at (2.75, -7) {};
		\node (7) at (6, -6.75) {};
		\node (8) at (8, -5.25) {};
		\node (9) at (-2.75, -1.25) {};
		\node (10) at (7.75, 3.75) {};
		\node (11) at (8.75, 1.25) {};
		\node (12) at (9.75, -2.25) {};
		\node (13) at (-0.75, -4.5) {};
		\node (14) at (2.75, -6.5) {};
		\node (15) at (5, -5.75) {};
		\node (16) at (7.75, -4.75) {};
		\node (17) at (9, -2.25) {};
		\node (18) at (8, 1.25) {};
		\node (19) at (2.75, 5) {};
		\node (20) at (0.5, 4.5) {};
		\node (21) at (-1, 2.75) {};
		\node (22) at (-1.75, -1.25) {};
		\node (23) at (5.75, 4) {};
		\node (24) at (6.75, 3) {};
		\node (25) at (-2.25, -3.5) {};
		\node (26) at (-1.75, -3.25) {};
		\draw [semithick, dashed, bend right=90, looseness=1.50] (5.center) to (6.center);
		\draw [semithick, dashed, bend right=60] (6.center) to (7.center);
		\draw [semithick, dashed, bend right=45, looseness=0.75] (7.center) to (8.center);
		\draw [semithick, dashed, bend right=60, looseness=1.25] (1.center) to (0.center);
		\draw [semithick, dashed, bend right=45] (2.center) to (1.center);
		\draw [semithick, dashed, bend left=60] (2.center) to (3.center);
		\draw [semithick, dashed, bend left=60, looseness=1.25] (9.center) to (0.center);
		\draw [semithick, dashed, bend left=60] (12.center) to (8.center);
		\draw [semithick, dashed, bend left=75, looseness=1.25] (3.center) to (10.center);
		\draw [semithick, dashed, bend left=60] (10.center) to (11.center);
		\draw [semithick, dashed, bend left=75] (11.center) to (12.center);
		\draw [semithick, dashed, bend left=15] (1.center) to (20.center);
		\draw [semithick, dashed, bend right=45, looseness=0.75] (2.center) to (19.center);
		\draw [semithick, dashed, bend left=15] (3.center) to (23.center);
		\draw [semithick, dashed, bend right] (10.center) to (24.center);
		\draw [semithick, dashed, bend left=45] (11.center) to (18.center);
		\draw [semithick, dashed, bend left] (12.center) to (17.center);
		\draw [semithick, dashed, bend right] (8.center) to (16.center);
		\draw [semithick, dashed, bend right] (15.center) to (7.center);
		\draw [semithick, dashed, bend left=15, looseness=0.75] (14.center) to (6.center);
		\draw [semithick, dashed, bend left=45, looseness=0.75] (13.center) to (5.center);
		\draw [semithick, dashed, bend right] (9.center) to (22.center);
		\draw [semithick, dashed, bend right] (0.center) to (21.center);
		\draw [semithick, dashed, bend right=75, looseness=1.25] (9.center) to (25.center);
		\draw [semithick, dashed, bend right=60, looseness=1.25] (25.center) to (5.center);
		\draw [thick, dashed, bend left=15] (25.center) to (26.center);
            \draw[thick,fill=white,white] (-3.6,.235) ellipse (.05 and .05);
            \draw[thick,fill=white,white] (-2.84,-1.145) ellipse (.05 and .05);
            \draw[thick,fill=white,white] (-2.81,-1.325) ellipse (.05 and .05);
            \draw[thick,fill=white,white] (-3.4,1.88) ellipse (.15 and .15);
            \draw[thick,fill=white,white] (-1.6,2.88) ellipse (.05 and .05);

    \draw[thick,gray] (-12,-.30) ellipse (0.11 and 0.33);
    \draw[thick,color=blue!5,fill=blue!15] (1.85,-.5) ellipse (1.15 and 1.15);

    \coordinate (A) at (-12,0);
    \coordinate (A') at (3,0);
    \coordinate (A'') at (-5.6,3);

    \draw[thick,red!35] (-12,.36) -- (-8,.3225);
    \draw[thick,red!35] (-12,.34) -- (-8,.3025);
    \draw[thick,red!35] (-8,.3025) -- (0,.2275);
    \draw[thick,red!35] (-12,.32) -- (-8,.2825);
    \draw[thick,red!35] (-8,.2825) -- (2,0.18875);
    \draw[thick,red!35] (-12,.30) -- (-8,.2625);
    \draw[thick,red!35] (-8,0.2625) -- (-2,0.20625);
    \draw[thick,red!35] (-12,.26) -- (-8,.2225);
    \draw[thick,red!35] (-8,.2225) -- (-6,0.20375);
    \draw[thick,red!35] (-12,.28) -- (-8,.2425);
    \draw[thick,red!35] (-8,.2425) -- (-3,0.205);
    \draw[dashed,thick,red!35,-stealth] (2,0.18875) -- (6,0.105);

    \draw [dashed, thin, -stealth, red](-13.150,.36) -- (-12.1,.36);
    \draw [dashed, thin, -stealth, red](-13.045,.26) -- (-12.1,.26);

    \draw[ultra thick,gray!50,ray] (-8,.3225) -- (-5,3.4025);
    \draw[ultra thick,gray!50,ray] (-5,3.4025) -- (-3,1.4025);
    \draw[ultra thick,gray!50,ray] (-3,1.4025) -- (-2.8,-1.7025);
    \draw[ultra thick,gray!50,ray] (-2.8,-1.7025) -- (-12,-.41);
    \draw[ultra thick,gray!50,-stealth] (-2.8,-1.7025) -- (-12,-.41);

    \draw[ultra thick,gray!50,ray] (0,.2275) -- (-0.1,-.525);
    \draw[ultra thick,gray!50,ray] (-0.1,-.525) -- (0.1,-1.3825);
    \draw[ultra thick,gray!50,ray] (0.1,-1.3825) -- (-0.3,-2.825);
    \draw[ultra thick,gray!50,ray] (-0.3,-2.825) -- (-12,-.48);
    \draw[ultra thick,gray!50,-stealth] (-0.3,-2.825) -- (-12,-.48);

    \draw[ultra thick,blue!30,ray] (2,0.18875) -- (3.1,0.75);
    \draw[ultra thick,blue!30,ray] (3.1,0.75) -- (-12,-.27);
    \draw[ultra thick,blue!30,-stealth] (3.1,0.75) -- (-12,-.27);

    \draw[ultra thick,gray!50,ray] (-2,0.20625) -- (1.5,1.25);
    \draw[ultra thick,gray!50,ray] (1.5,1.25) -- (0.75,2);
    \draw[ultra thick,gray!50,ray] (0.75,2) -- (-12,-.20);
    \draw[ultra thick,gray!50,-stealth] (0.75,2) -- (-12,-.20);

    \draw[ultra thick,blue!70,ray] (-3,0.205) -- (2.5,-.62425);
    \draw[ultra thick,blue!70,ray] (2.5,-.62425) -- (-12,-.34);
    \draw[ultra thick,blue!70,-stealth] (2.5,-.62425) -- (-12,-.34);

    \draw[ultra thick,gray!50,ray] (-6,0.20375) -- (-1.45,1.8825);
    \draw[ultra thick,gray!50,ray] (-1.45,1.8825) -- (-0.15,2.9325);
    \draw[ultra thick,gray!50,ray] (-0.15,2.9325) -- (-12,-.13);
    \draw[ultra thick,gray!50,-stealth] (-0.15,2.9325) -- (-12,-.13);


\end{tikzpicture}}}%
    \caption{Despite most of the scattering happening inside the plume (dashed region), most trajectories don't cross the small kernel in figure (b). As such the marginal distributions of parameters reconstructed from wide FOV data corresponding to small features will suffer from essentially the same issues as the shown in \cref{fig:free_test} for high ambient scattering.}%
    \label{fig:wFOV_feature_cloud}%
\end{figure}
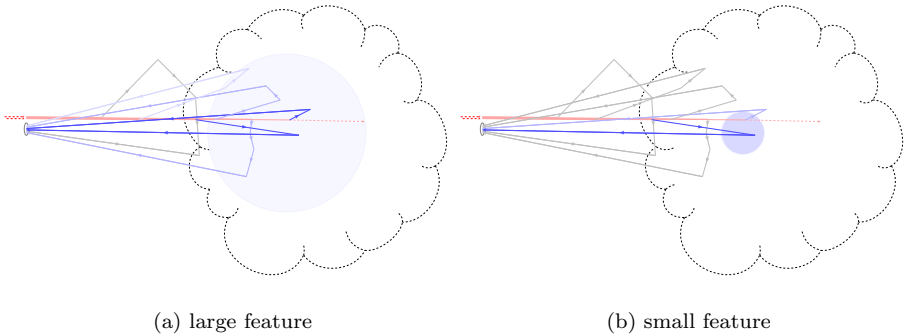
Similar results are to be expected for position and size parameters for which an alignment of the scattering parameters with respect to the corresponding gradients is required. Furthermore the ambient scattering may be replaced by scattering particles from the plume surrounding a smaller feature. In other words, the information within wide FOVs declines as soon as the plume feature can only be resolved with a higher resolution than the bulk of the plume, even if the scattering parameters are driven by the same dispersion model as the absorption profile. This in effect limits the usefulness of wide FOVs when a high resolution image is sought (compare also \cref{fig:wFOV_feature} and \cref{fig:wFOV_feature_cloud}). We note that the behaviour for high ambient scattering is in a way similar, although for a different type of measurement, to the findings of decreased sensitivity at larger optical depths from \cite{LoveridgeEtAl:2023a, LoveridgeEtAl:2023b}. The approach taken here differs insofar that image features aren't strictly local but rather we consider more smoothly varying aspects of the gas concentration profile which as a whole will ca be rather stable even for plumes of large optical depths.

Although ambient noise has been ignored in this comparison it is not hard to see that its effect will be very similar to that of ambient scattering particles, i.e. it doesn't matter whether we observe photons that didn't originate from the controlled source or light that never reached the region of interest (gray trajectories in \cref{fig:wFOV_feature} and \cref{fig:wFOV_feature_cloud}). In much the same as with ambient scattering it becomes an issue once the number of peak signal photons observed is of a similar magnitude as the ambient noise. 

\section{Computational modelling and reconstruction}
\label{sec:computing}
As it turns out, the relaxed semi-parametric form is, when paired with an expression as in \Cref{eq:DIAL_llik}, becomes very convenient for optimisation purposes. Consider the loss function
\begin{align}\label{eq:loss_functional}
\begin{split}
    \mathsf{Loss}(\tilde{\psi},\theta \mid \bm{m},\bm{n}) &= \sum_{i=1}^{N_v}\sum_{j = 1}^{N_t} H_D \Delta t A_D (\mathfrak{m}_{\mathrm{on}}[\theta](v_i,t_j\mid x_D)+\mathfrak{m}_{\mathrm{off}}[\theta](v_i,t_j\mid x_D))\tilde{\psi}_{i,j} \\
    &-\bm{m}_{v_i,t_j} \log(\mathfrak{m}_{\mathrm{on}}[\theta](v_i,t_j\mid x_D)) - \bm{n}_{v_i,t_j} \log(\mathfrak{m}_{\mathrm{off}}[\theta](v_i,t_j\mid x_D))\\
    &-(\bm{m}_{i,j}+ \bm{n}_{i,j})\log(\tilde{\psi}_{i,j})-\mathsf{R}(\theta)
\end{split}
\end{align}
for \(\tilde{\psi}\) free. It is formed as the sum of \Cref{eq:DIAL_llik} and a prior term \(\mathsf{R}\) that enforces dispersion based constraints such as alignment according to the wind direction or continuity and smoothness constraints (i.e. no ``holes``) of the kernel components. 

Optimisation of \Cref{eq:loss_functional} can now be carried out sequentially. Note that \(\mathsf{R}\) by assumption only depends on \(\theta\) so optimising in the nuisance parameter \(\tilde{\psi}\) for fixed \(\theta\), i.e. taking
\begin{align*}
    \tilde{\psi}^*(\theta) = \arg \max_{\tilde{\psi}} \mathsf{Loss}(\theta, \tilde{\psi}\mid \bm{m}),
\end{align*}
it is easily seen that we have
\begin{align}\label{eq:profile_nuisance}
    \tilde{\psi}_{i,j}^*(\theta) = \frac{\bm{m}_{i,j}+\bm{n}_{i,j}}{H_D \Delta t A_D (\mathfrak{m}_{\mathrm{on}}[\theta](v_i,t_j\mid x_D)+\mathfrak{m}_{\mathrm{off}}[\theta](v_i,t_j\mid x_D))}.
\end{align}
We may plug \(\tilde{\psi}_{i,j}^*(\theta)\) into \Cref{eq:loss_functional} which yields, after rearranging and removal of additive constants, the expression
\begin{align}\label{eq:profiled_loss}
    \begin{split}
    \mathsf{Loss}(\theta \mid \bm{m},\bm{n}; \tilde{\psi}_{i,j}^*) 
    &= \sum_{i=1}^{N_v}\sum_{j = 1}^{N_t} \bm{m}_{i,j} \log\left(\frac{\mathfrak{m}_{\mathrm{on}}[\theta](v_i,t_j\mid x_D)}{\mathfrak{m}_{\mathrm{on}}[\theta](v_i,t_j\mid x_D)+\mathfrak{m}_{\mathrm{off}}[\theta](v_i,t_j\mid x_D)}\right) \\
    &+\sum_{i=1}^{N_v}\sum_{j = 1}^{N_t} \bm{n}_{i,j} \log\left(\frac{\mathfrak{m}_{\mathrm{off}}[\theta](v_i,t_j\mid x_D)}{\mathfrak{m}_{\mathrm{on}}[\theta](v_i,t_j\mid x_D)+\mathfrak{m}_{\mathrm{off}}[\theta](v_i,t_j\mid x_D)}\right) \\
    & - R(\theta)
\end{split}
\end{align}
which has essentially the same structure as the negative \(\log\)-likelihood of a collection of binomial random variables with trial lengths \(\bm{n}_{i.j}+\bm{m}_{i.j}\) and probabilities
\begin{align*}
    P_{i,j}(\theta) = \frac{\mathfrak{m}_{\mathrm{on}}[\theta](v_i,t_j\mid x_D)}{\mathfrak{m}_{\mathrm{on}}[\theta](v_i,t_j\mid x_D)+\mathfrak{m}_{\mathrm{off}}[\theta](v_i,t_j\mid x_D)}.
\end{align*}
As such \Cref{eq:profiled_loss} may be optimised by what can be regarded as Fisher-Scoring for a binomial distribution which will avoid a computation of second derivatives. In other words, we iterate
\begin{align}\label{eq:theta_iteration}
    \theta^{(q+1)} &= \theta^{(q)} + \zeta \mathsf{H}(\theta^{(q)})^{-1}\partial_{\theta} \mathsf{Loss}(\theta^{(q)} \mid \bm{m},\bm{n}; \tilde{\psi}_{i,j}^*)
\end{align}
where \(\zeta\) is a step size and \(\mathsf{H}(\theta)\) is given by
\begin{align*}
    \mathsf{H}(\theta) &= \partial_{\theta \theta}\mathsf{R}(\theta) + \sum_{i=1}^{N_v}\sum_{j = 1}^{N_t} (\bm{n}_{i.j}+\bm{m}_{i.j}) \frac{\partial_{\theta}P_{i,j}(\theta) \partial_{\theta} P_{i,j}(\theta)^{\top}}{P_{i,j}(\theta)(1-P_{i,j}(\theta))}
\end{align*}

\subsection{Computational complexity \& randomisation induced errors}
Performing the iterations as described earlier is straight-forward once we have a method for computing \(P_{i,j}\) as well as the gradients \(\partial_{\theta}P_{i,j}\). A solution of \Cref{eq:time_dependent_RTE} with source term \(\delta(x-x_D)\delta(v-v_i)\delta(t)\) yields \(P_{i,j}\) for all \(j=1,\dots, N_t\) and a given \(i \in \{1,\dots, N_v\}\) and each partial derivative can be obtained at the same cost. means that \((1+\dim (\theta)) N_v\) solutions of the RTE are necessary to evaluate the loss functional and perform the an update \cref{eq:theta_iteration}. In practice, the cost of all other operations is negligible in comparison and we may focuse our attention to the complexity of an RTE evaluation. 

Given that the present work is primarily of interest when the data is noisy, it is not be necessary to evaluate the forward map and thus the RTE to a high degree of accuracy and we are be able to get away with random approximations \(\hat{P}_{i,j}(\theta)\) and derivatives thereof obtained from MC ray-tracing as described in \cite{GkioulekasLevinZickler2016} which makes in effect results in a stochastic optimisation procedure. The derivatives of \(\hat{P}_{i,j}(\theta)\) may be computed through score function estimators \cite{mohamed:2020, GkioulekasLevinZickler2016} and exhibit an overall similar behaviour. Sampling a photon path of a given length \(t_j\) for scattering order \(k\) requires \(k N\) flops, where \(N\) is the number of kernel functions used for the plume approximation as in \cref{eq:kernel_parameterisation}. Assuming a maximal order of scattering \(k_{\max}\) the cost of evaluating the RTE becomes \(\mathcal{O}(N_p N_t k_{\max}^2 N )\) with \(N_p\) being the number of sampled paths. The values \(N\) and \(N_t\) are determined by the plume and measurement granularity and effectively fixed while \(N_p\) controls the accuracy of an RTE evaluation. 
It should be noted that the sums of random vectors and matrices such as \(\mathsf{H}(\theta)\) or \(\partial_\theta \mathsf{Loss}\) used in out optimisation behave much like approximations used in randomised numerical linear algebra. Indeed our reduction of the image to dispersion related parameters means that \(\dim(\theta) \ll N_v N_t\) and, assuming that we ensure boundedness of the integrand that is to be evaluated via MC path tracing, we may use the matrix Bernstein inequality \cite{Tropp:2012} to obtain concentration estimates for quantities such as the spectral norm error \(\|\hat{H}(\theta) - H(\theta)\|\). For further details we refer to \cite{Tropp:2015, MartinssonTropp:2020}. Estimates for the extreme eigenvalues of the spd matrix \(\hat{H}(\theta)\) in the form of matrix Chernoff inequalities are also available \cite{Tropp:2012, CohenEtAl:2015}.

The sampling used in this work is straight-forward and no attempt to further reduce the variance has been made, nor was any action necessary. In situations with less well-behaved parameters or more complex environments a further reduction in variance is possible through the use of more sophisticated, bi-directional path tracing strategies \cite{Pediredla:2019}.

\section{A numerical example}\label{sec:experiments}
In order to validate our method on simulated data we consider a simulated reconstruction from \(30 \times 10 \times 50\) Lidar scan of a \(14\) parameter dispersion which can be recovered when conventional reconstruction fails due to the low SNR. As suggested by our developments in \cref{sec:stability}, we consider instances where the scattering is caused largely by particles around the gas plume and the effectively required resolution is roughly granularity of absorbing gas within the scatterer. 

In order to avoid contrived scenarios we have fixed the phase function \(f_p\) to different values for simulation and reconstruction respectively.
Fixed system parameters as in \cref{tab:system_parameters} were selected in order to mimic realistic conditions as if a DIAL instrument for methane measurement was used for the experiment. 
\begin{table}[]
\caption{System parameters used in the simulation}
\begin{center}
\small
\begin{tabular}{|p{4cm}|p{6cm}|}  \hline
     Detector & 3cm lens with 4\% detection rate \\ \hline
     Methane amount & 50mol or 0.8kg \\ \hline
     Distance & 100m \\ \hline
     Wavelength (absorbing) & 1645.55nm \\ \hline
     Pulse Energy & 250\(\mu\)J \\ \hline
     Ambient intensity & 0.025W uniformly over hemisphere \\ \hline
\end{tabular}
\end{center}
\label{tab:system_parameters}
\end{table}

In addition to the optical noise we have also considered perturbations of the dispersion model via a super-position of branching jump diffusion processes with the idea to validate the ability of the non-parametric component to capture noise that isn't accounted for by the likelihood. The observed level of noise will in general depend on atmospheric conditions (which affect the amount of turbulence within dispersion process) and the amount of (temporal) averaging done as part of the data acquisition. In particular, the plume will in practice be a dynamic object and for a stationary release rate its temporal average will, if taken over sufficiently long periods, resemble the low dimensional, smooth average model. For the purpose of this simulations we did not simulate a dynamic turbulent dispersion process (mainly due to computational constraints) and instead used a realisation of a random process as midpoints which combined with suitably chosen kernel widths will in expectation satisfy \cref{eq:general_plume_equation}. Note that the resulting method can be thought of as snapshot of a process which aims to mimic the empirical observation that ``Big whorls have little whorls which have lesser whorls...`` in a turbulent dispersion process. For further details regarding the noise structure as well as the hyperparameters and prior distributions used in the fitting of the dispersion parameters we refer to \cref{sec:implementation_details}.

\cref{tab:errors} shows different error metrics for various situations and measurements using a single pulse per direction with nFOV, wFOV and mFOV denoting narrow, wide and multiple FOVs respectively. The ratio (nFOV:wFOV) is a measure for the optical depth of the plume, a value of 7.2 indicates that 7.2 times as much light was measured in the narrow FOV than in the wide FOV. The first column indicates whether the plume was perturbed additionally via the random process described in \cref{sec:implementation_details}. ``Yes`` therefore indicates a higher amount of noise in the measurement and higher errors are to be expected. Note that the level of optical noise depends on the counts ratio in a non trivial way. While the number of photons increases with more scattering particles, i.e. a lower ratio, the increased optical depth means that back scattered light does not penetrate the plume as much which results in a lower differential absorption.  

\begin{table}[h]
\caption{Reconstruction errors (average from 10 simulation runs per row)}
\begin{center}
\small
\begin{tabular}{|p{2.5cm}|p{1.6cm}|p{3.5cm}|p{3.5cm}|}  \hline
 Dispersion noise & Counts \tiny{(nFOV:wFOV)} & $\mathrm{L_1 \; errors}$ \tiny{(nFOV|wFOV|mFOV)}  & Release amount errors \tiny{(nFOV|wFOV|mFOV)}\\ 
 \hline
 Yes & 7.2 & (44\% |\, 45\% |\, 45\%) & (11\% |\, 12\% |\, 14\%)\\ 
 Yes & 3.9 & (44\% |\, 45\% |\, 46\%) & (12\% |\, 12\% |\, 10\%)\\
 Yes & 2.4 & (56\% |\, 39\% |\, 36\%) & (21\% |\, 12\% |\, 8\% )\\  
 Yes & 1.7 & (64\% |\, 36\% |\, 36\%) & (32\% |\, 14\% |\, 14\%)\\  
 No  & 7.2 & (42\% |\, 32\% |\, 33\%) & (18\% |\, 16\% |\, 15\%)\\  
 No  & 3.9 & (30\% |\, 22\% |\, 23\%) & (14\% |\, 9\%  |\, 8\% )\\  
 No  & 2.4 & (30\% |\, 20\% |\, 21\%) & (14\% |\, 8\%  |\, 10\%)\\  
 No  & 1.7 & (41\% |\, 25\% |\, 24\%) & (13\% |\, 14\% |\, 13\%)\\ 
 \hline
\end{tabular}
\end{center}
\label{tab:errors}
\end{table}

\Cref{tab:errors}, which essentially corresponds to a more general point estimation problem in the situation of \cref{fig:known_test}, shows that in those cases the wider FOV tends to outperform a like for like reconstruction that uses narrow FOVs only. In fact, the average improvement with respect to \(L_1\) and release rate is 24.5\% and 23.6\% respectively. Note that the reduction in error from an additional (statistically independent) narrow FOV measurement would be 29.3\%, i.e. an accuracy much closer to what one might expect from two instead of one sets of measurement. Of course this does not take accurately into account how the turbulent component might behave. We emphasise that the absolute errors should not be taken as any indication of what can be expected in practice as this will be highly dependent on the system parameters (see \cref{tab:system_parameters}). As long as the ambient intensity does not attain levels comparable to the signal response, i.e. we have a high enough pulse energy, we may use our findings as an indication of the relative improvement that the proposed measurements can have compared to classical narrow FOV methods.

\begin{figure}[]%
    \begin{center}
    \subfloat[\centering narrow FOV]{{\includegraphics[scale=0.475,trim={4.8cm 1cm 3cm 2cm},clip]{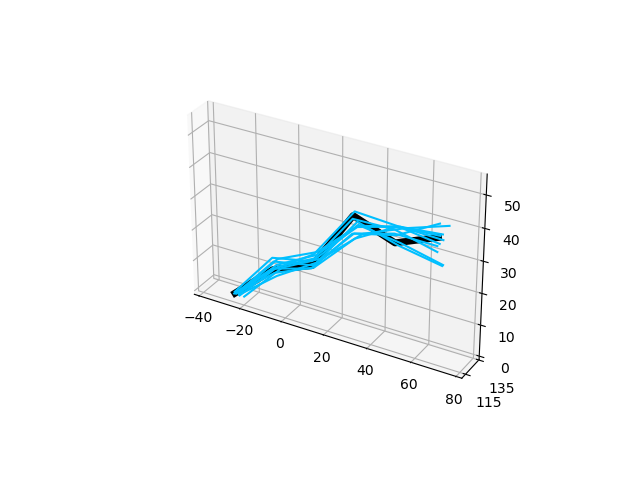} }}%
    \subfloat[\centering wide FOV]{{\includegraphics[scale=0.475,trim={4.8cm 1cm 3cm 2cm},clip]{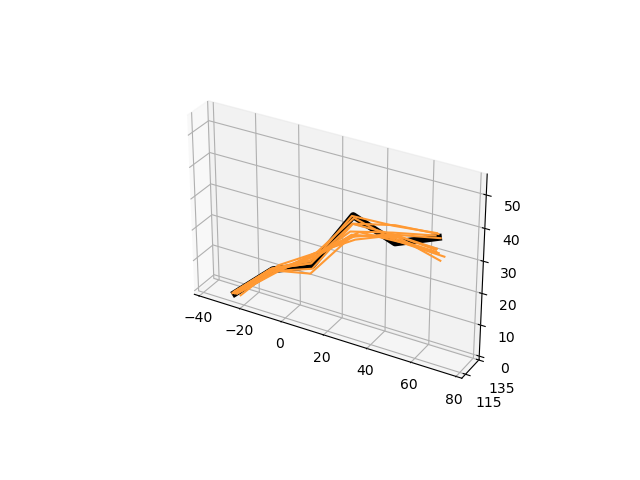} }}%
    \subfloat[\centering multiple FOVs]{{\includegraphics[scale=0.475,trim={4.8cm 1cm 3cm 2cm},clip]{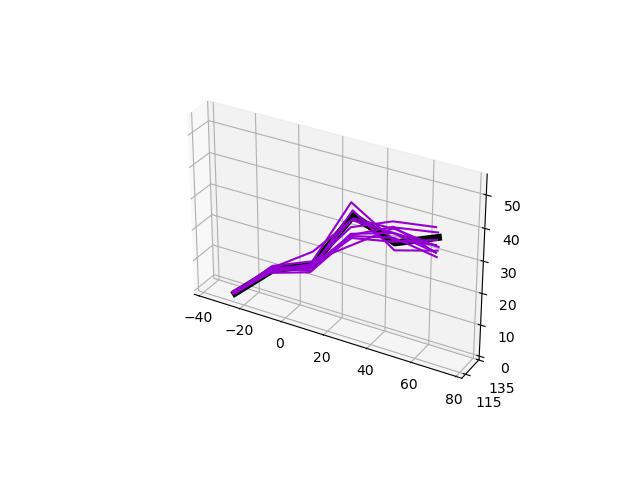} }}%
    \caption{Reconstructed plume centre-lines for a smooth plume and a scattering particle concentration corresponding to a 2.4 nFOV:wFOV ratio in the measured data (line 7 in \cref{tab:errors})}%
    \end{center}
    \label{fig:release_error_mv_smoth}%
\end{figure}

\begin{figure}[]%
    \centering
    \includegraphics[scale=.667,trim={0cm .75cm 0cm 1.45cm},clip]{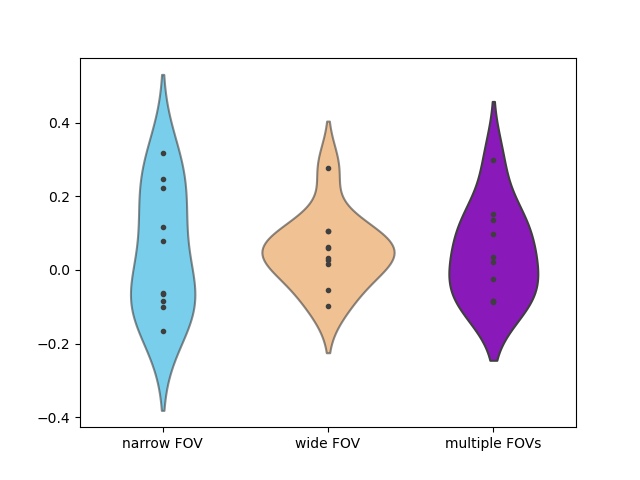}
    \caption{Relative deviation of reconstructed release rates corresponding to the same instance as shown in \cref{fig:release_error_mv_smoth} (line 7 in \cref{tab:errors})}%
    \label{fig:release_error_mv_smoth_violin}%
\end{figure}

\Cref{fig:release_error_mv_smoth} and \Cref{fig:release_error_mv_smoth_violin} show an instance that can be thought of as a typical level of improvement (line 7 in \cref{tab:errors}). In particular, the primary benefits are found in the reconstruction of the release rate as well as close to the source. This is to be expected as the release rate is the ``smoothest`` possible parameter while the source location has a global effect on the image. The errors in the centre-lines increase with the distance from the source as the plume widens which results in a less (multiple) scattering along with a lower sensitivity to the relevant parameters due to increased smoothness. Note that the improvements are observed in the instance with perhaps the lowest absolute errors which suggests that there unknown scattering (which are only partly compensated by the non-parametric nuisance component) do not cause any notable bias in the reconstruction. A similar pattern would be observed for the width of the plume which is not shown explicitly but large errors would have a negative effect on the \(L_1\) errors in \cref{tab:errors}. It should be noted that multiple and wide FOV reconstructions suffer from outliers which can skew the averages towards larger errors. However, we acknowledge that the increased amount of randomness from the RTE solver paired with a lack of information regarding the scattering parameters appears to cause such phenomena with non-negligible regularity and at present we have no way to avoid these situations. As such they should not, and were not, disregarded from performance metrics.

\begin{figure}[]%
    \centering
    \subfloat[\centering narrow FOV]{{\includegraphics[scale=0.45,trim={4.5cm 1cm 2.8cm 2cm},clip]{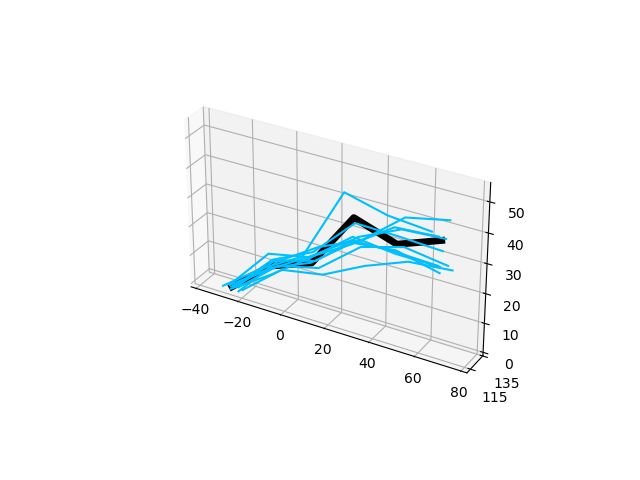} }}%
    \subfloat[\centering wide FOV]{{\includegraphics[scale=0.45,trim={4.5cm 1cm 2.8cm 2cm},clip]{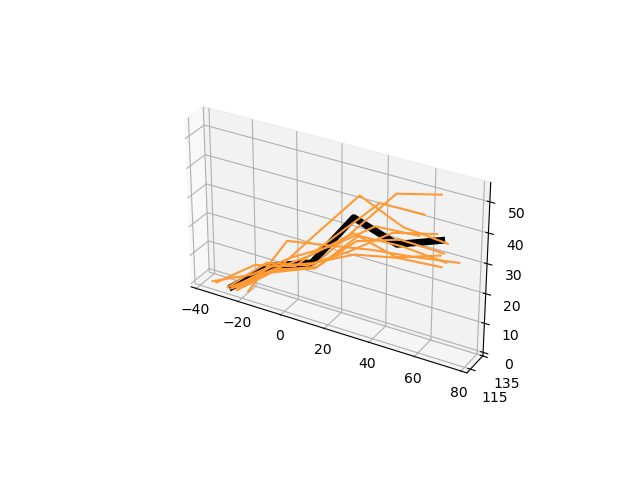} }}%
    \subfloat[\centering multiple FOVs]{{\includegraphics[scale=0.45,trim={4.5cm 1cm 2.8cm 2cm},clip]{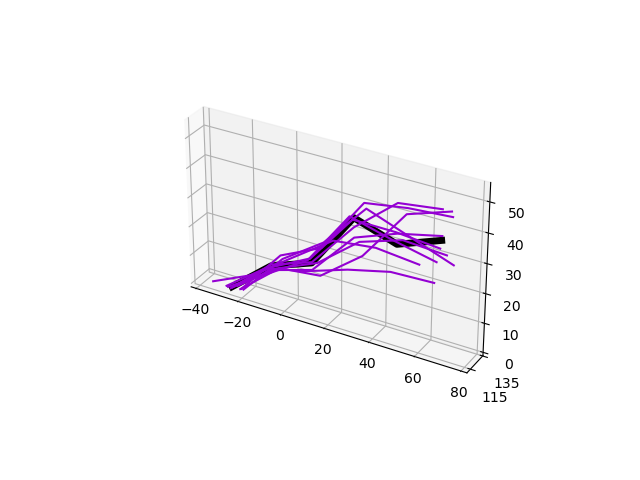} }}%
    \caption{Reconstructed plume centre-lines for a turbulent plume and a scattering particle concentration corresponding to a 7.2 nFOV:wFOV ratio in the measured data (line 1 in \cref{tab:errors})}%
    \label{fig:release_error_vv}%
\end{figure}

\begin{figure}[]%
    \centering
    \includegraphics[scale=.667,trim={0cm .75cm 0cm 1.45cm},clip]{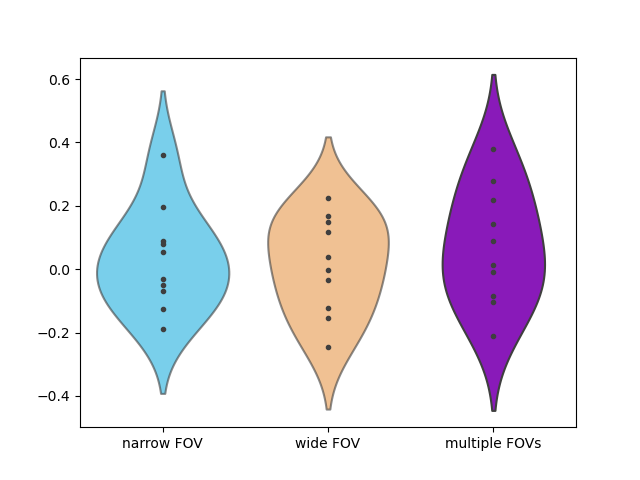}
    \caption{Relative deviation of reconstructed release rates corresponding to the same instance as shown in \cref{fig:release_error_vv} (line 1 in \cref{tab:errors})}%
    \label{fig:release_error_vv_violin}%
\end{figure}

An instance with less improvement is shown in \cref{fig:release_error_vv} and \cref{fig:release_error_vv_violin} (see line 1 in \cref{tab:errors}). Note that the absolute errors are rather large are due to the presence of turbulence which is particularly evident in the centre lines in \cref{fig:release_error_vv}. The relatively poor performance of multiple FOVs (relative to the narrow FOV data) is caused by a lack of scattering. In such instances the additional uncertainty due to a lack of signal, resulting in higher optical noise component, as well as unknown parameters in the wider FOV outweighs the potential benefits. This instance roughly corresponds to the bottom left corner of \cref{fig:known_test} and although the difference is small, narrow FOV measurements can be expected to outperform wide FOV in ``extreme`` cases.  

It is also worth noticing that even though the condition from \cref{eq:kernel_tails} has not been implemented specifically and injectivity of the forward operator is not guaranteed by \Cref{thm:relaxed_uniqueness}, whose proof heavily relies on the tail behaviour of the kernels, the perturbation of the plume which causes particularly large (relative) errors seems to affect the narrow FOV data more so than the wide. Such an observation can be intuitively expected as condition \cref{eq:kernel_tails} is a technical assumption relevant only in what can be thought of as a best case scenario whereas the narrow FOV in instances of high optical thickness (see line 3 and 4 in \cref{tab:errors}) collects more light scattered at the ``outside`` of the plume and is therefore to a greater degree affected by the model errors in that region. 

\section{Conclusions}
\label{sec:conclusions}

We have presented a method for imaging dispersion processes using differential absorption based measurements from a single Lidar source that is based on the RTE and can therefore utilise data from detectors that offer a wider FOV. We showed that any of these measurements has enough information to recover an image and, under the assumption of \(\mathsf{Poisson}\) noise in the data, argued that different aspects of the image are more efficiently captured by different detector types suggesting that data from wider FOVs can aid in determining gas concentration. Our reconstruction method is based on a semi-parametric model which is chosen such that it can capture variability due to unknown scattering and allows for computationally feasible reconstruction. Preliminary simulations confirm that the chosen parameterisation of the problem is able to capture scattering effects that cannot be explained by the dispersion parameters alone and that a plume shaped concentration profile can be recovered from knowledge of absorption behaviour from noise corrupted data that is unsuitable for the use in direct inversion formulas of classical DIAL.
\newpage
\appendix
\section{Supplementary proofs and remarks}
\label{sec:proofs}
\subsection{Continuity of RTE solutions}
The proof consists of a relatively straight-forward change of variable argument but relies on a few explicit computation which makes the expressions complicated and the overall argument rather hard to follow. Most of the ideas and techniques used are not necessarily new and are given mostly for completeness.
\begin{proof}[Proof of \cref{lem:path_continuity}]
    We give a proof for a constant \(h=1\), i.e. we drop it from the formulas, as the argument is virtually identical for other choices of \(h\). The ballistic term \(\mathfrak{m}_0\) vanishes on \(M\) and the claim can be verified easily for \(k = 1\). Let's start by fixing \(k \geq 2\). We can express
    \begin{align*}
          \int_{\partial V_{(+)}(x)} \mathcal{A}_k[g](x,v,t) dv &= \int_{(0,\infty)^k \times (\mathbb{S}^2)^{k+1}} g\Big(x_1 -\ell_{v_0}\left(x_1 \right)v_0,v_0,t - \sum_{l=1}^{k} s_{l}-\ell_{v_0}\left(x_1\right)\Big) \\ 
          &\phantom{=} \times F_k(x;v_{k}, \dots, v_0 ,s_k, \dots, s_1)  dv_{k} \dots dv_0 d s_k \dots d s_1
    \end{align*}
    where we have abbreviated
    \begin{align*}
        F_k(x; v_{k}, \dots, v_0 ,s_k, \dots, s_1) &= \exp\Big(- \int_{\Gamma(x_{0:k+1})} \sigma_{a+s}(r)dr\Big) \prod_{l=1}^{k} \sigma_{s}(x_l) f_p(x_l, v_{l-1}\cdot v_l)\\ 
        &\phantom{=} \times \vert \Vec{n}(x) \cdot v_k \vert 1_{\partial V_{(+)}(x)}(v_k) 1_{X^k}(x_k, \dots, x_1) 
    \end{align*}
    and introduced for any \(l = 1, \dots, k\) the implicit quantities \(y_{l} = s_{l}v_{l}\), \(x_{l} = x - \sum_{m=l}^{k} y_{m}\) as well as \(x_{k+1}=x\) and \(x_0 = x_1 -\ell_{v_0}\left(x_1 \right)v_0\). The set \(\Gamma(x_{0:k+1})\) is the graph corresponding to the shortest piece-wise linear curve that connects the surface points \(x_0\) and \(x_{k+1}\) while passing through \(x_k, \dots,x_1 \in X\) in the order of the indexing, i.e.
    \begin{align}\label{eq:path_definition}
        \Gamma(x_{0:k+1}) = \bigcup_{l=0}^k\big\{\lambda x_l + (1-\lambda)x_{l+1}: 0<\lambda <1\big\}.
    \end{align}
    Note that if we had surface scattering \(F_k\) would be more complicated as the integration over surface points in \(\partial X\) is quite different to that in the interior \(X\). In pursuit of our goal to obtain an expression for \(\mathfrak{m}_k\), which implies changing the variables in such a way that the integration is performed w.r.t. \(x_1 -\ell_{v_0}\left(x_1 \right)v_0\) and \(\sum_{l=1}^{k} s_{l}+\ell_{v_0}\left(x_1\right)\), we first make the change of variables \((s_{m},v_m) \mapsto y_m\) followed by \(y_m \mapsto x_1\), which is a simple shift, and utilise \cite{ChoulliStefanov1999} Lemma 2.1 to obtain
    \begin{align*}
          \int_{\partial V_{(+)}} \mathcal{A}_k[g](x,v,t) dv =& \int_{(0,\infty)^k \times (\mathbb{S}^2)^{k-1} \times \partial D_{(-)}} \frac{F_k(x; \dots)\vert \Vec{n}(x_0) \cdot v_0 \vert}{\big\|x-x_0 - \sum_{l=0,l \neq m}^k s_{l}v_l\big\|^2_2} \\ 
          & \times g\left(x_0,v_0,t -\textstyle \sum_{l=0,l \neq m}^k s_{l} -\|x-x_0 - \sum_{l=0,l \neq m}^k s_{l}v_l\|_2 \right) \\
          & dv_0 dx_0 dv_{k} \dots dv_{m+1} dv_{m-1} \dots dv_{1} d s_k \dots d s_{m+1}d s_{m-1}\dots ds_0.
    \end{align*}
    We have used the dots in \(F_k(x;\argdot)\) in order to abbreviate the dependence of \(F_k\) on the independent variables \(v_0, x_0, v_{k}, \dots, v_{m+1}, v_{m-1}, \dots, v_{1}, s_k \dots, s_{m+1}, s_{m-1},\dots s_0\). Note that in this change we have
    \begin{align}\label{eq:surface_variable_transformation}
        s_m = \Big\|x-x_0 - \sum_{l=0,l \neq m}^k s_{l}v_l\Big\|_2 \quad , \qquad v_m = \frac{x-x_0 - \sum_{l=0,l \neq m}^k s_{l}v_l}{\big\|x-x_0 - \sum_{l=0,l \neq m}^k s_{l}v_l\big\|_2}
    \end{align}
    as the inverse maps of the coordinate transform for the new variables. In order to obtain an expression for the Schwartz kernel \(\mathfrak{m}_k\) it remains to resolve the time component of the arguments in \(g\). The necessary transformation takes the form
    \begin{align*}
        s_n \mapsto s = \sum_{l=0,l \neq m}^k s_{l} + \Bigg\|x-x_0 - \sum_{l=0,l \neq m}^k s_{l}v_l\Bigg\|_2 
    \end{align*}
    and maps \((0,\infty)\) to \((\sum_{l=0,l \neq m,n}^k s_{l}+\|x-x_0 - \textstyle \sum_{l=0,l \neq m,n}^k s_{l}v_l\|_2,\infty)\). It is invertible everywhere but on a set of measure zero and the inverse is given by
    \begin{align}\label{eq:time_variable_transformation}
        s_n = \frac{\big(s - \sum_{l=0,l \neq m,n}^k s_{l}\big)^2 - \big\|x-x_0 - \sum_{l=0,l \neq m,n}^k s_{l}v_l\big\|_2^2}{2 \left(s - \sum_{l=0,l \neq m,n}^k s_{l} - v_n \cdot (x-x_0 - \sum_{l=0,l \neq m,n}^k s_{l}v_l)\right)}.
    \end{align}
    It is easily verified that
    \begin{align*}
        \frac{d}{d s_n} s = 1 - v_n \cdot v_m \implies  \frac{d}{d s} s_n = \frac{1}{1 - v_n \cdot v_m}
    \end{align*}
    where \(v_m\) is as in \Cref{eq:surface_variable_transformation} but should be considered a function of \(s\) instead of \(s_n\) in the latter expression. Making that change results in
    \begin{align*}
          \int_{\partial V_{(+)}} \mathcal{A}_k[g](x,v,t) dv = &\int_{\partial D_{(-)} \times (0,\infty)}\int_{(0,\infty)^{k-1} \times (\mathbb{S}^2)^{k-1}}\frac{g\left(x_0,v_0,t - s \right) \vert \Vec{n}(x_0) \cdot v_0 \vert}{s_m(\argdot)^2 (1-v_{n} \cdot v_m(\argdot))} \\ 
          &\times 1_{(0,s)}\left(\textstyle \sum_{l=0,l \neq m,n}^k s_{l} + \|(x-x_0 - \sum_{l=0,l \neq m,n}^k s_{l}v_l)\|_2\right)\\ 
          &\times F_k(x;\argdot)  dv_{k} \dots dv_{m+1} dv_{m-1} \dots dv_{1} \\ 
          & d s_k \dots d s_{m+1}d s_{m-1}\dots d s_{n+1}d s_{n-1}\dots ds_0 ds dv_0 dx_0 
    \end{align*}
    and \((\argdot)\) has been used as a placeholder to indicate that a quantity is implicit, i.e. a function of the variables that the integration is performed over. Although this process in principle yields an explicit expression for \(\mathfrak{m}_k\) this isn't good enough because the change of variables introduces singularities and the integrand is unbounded making it difficult to analyse. However, we don't have to always use the same variables \(s_m,s_n,v_m\) for the change but instead we can partition \((0,\infty)^k \times  (\mathbb{S}^2)^{k+1}\) and perform the change with respect to different coordinates on each subset. We would like to pick \(s_m,v_m\) such that \(s_m\) is not too small, e.g. we can take the maximum, and \(v_n\) such that it isn't too aligned with \(v_m\) which may be accomplished by taking the most misaligned direction. The latter will avoid singularities whenever \(s > \|x-x_0\|\). In the situation we will be interested in we have \(s \gg \|x-x_0\|\) since \(x = x_0\) and \(s \gg 0\) for Lidar measurements. 
    First, notice that \(s = \sum_{l=1}^{k} s_{l}+\ell_{v_0}\left(x_1\right) > 0\) so taking \(\Gamma(x_{0:k+1})\) and connecting the endpoints yields at a closed polygon of length \(s+\|x-x_0\|\) which must have at least two segments, other than the piece between \(x\) and \(x_0\), that are longer than \(\frac{s-\|x-x_0\|}{2k}\). To see that note that due to the triangle inequality we must have
    \begin{align*}
        \frac{s+\|x-x_0\|}{2k} \leq \max\{s_k, \dots, s_1, \ell_{v_0}(x_1)\} \leq \frac{s+\|x-x_0\|}{2}.
    \end{align*}
    Taking the total length of the polygon and subtracting \(\max(s_k, \dots, s_1, \ell_{v_0}(x_1))\) as well as \(\|x-x_0\|\) from \(s+\|x-x_0\|\) to account for the added segment we obtain the lower bound for the second largest value of \(s_k, \dots, s_1, \ell_{v_0}(x_1)\) which particularly implies that
    \begin{align*}
        \max\{s_k, \dots, s_1\} \geq \frac{1}{k} \big(s- \max\{s_k, \dots, s_1, \ell_{v_0}(x_1)\}\big) \geq \frac{s-\|x-x_0\|}{2k}.
    \end{align*}
    Since we also have that 
    \begin{align*}
        s^2 - \|x_0 - x\|^2 = \sum_{l=1}^{k}\sum_{m=1}^{k} s_m s_l (1-v_m \cdot v_{l}) + 2 \ell_{v_0}(x_1)\sum_{l=1}^{k}s_l(1-v_0 \cdot v_l) 
    \end{align*}
    it must hold that
    \begin{align*}
        \max_{l,m=0}^k 1-v_m \cdot v_{l} \geq \frac{s^2-\|x-x_0\|_2^2}{s^2} = 1-\frac{\|x-x_0\|_2^2}{s^2}.
    \end{align*}
    Further we have that \(\sqrt{2-2a \cdot b} = \|a-b\|_2\) for \(a,b \in \mathbb{S}^2\) which means we can apply the triangle inequality to obtain for any \(m=0,\dots,k\)
    \begin{align*}
        \max_{l=0}^k \sqrt{1-v_m \cdot v_{l}} \geq \frac{1}{2}\sqrt{1-\frac{\|x-x_0\|_2^2}{s^2}}.
    \end{align*}
    In order to see how we can exploit these findings consider for any pair \((m,n)\), where \(m=1,\dots,k\) and \(n = 0,\dots,k\), the indicator functions
    \begin{align*}
        J_{m,n}(v_k, \dots, v_0, s_k, \dots, s_1) = \prod_{p=1}^k 1_{[0,\infty)}(s_m - s_p)\prod_{q=0,q\neq n}^k 1_{(0,\infty)}(v_m\cdot v_q - v_m\cdot v_n)
    \end{align*}
    and observe that we have \(dv_{k} \dots dv_0 d s_k \dots d s_1\) almost everywhere
    \begin{align*}
        \sum_{m = 1}^k \sum_{n = 0, n\neq m}^k J_{m,n} = 1.
    \end{align*}
    In fact, it is not hard to see that \(J_{m,n}=1\) if \(s_m \geq s_p\) for any \(p =1,\dots,k\) and at the same time \(1-v_m\cdot v_n > 1-v_m\cdot v_q\) for any \(q \neq n\). Some \(s_m\) is inevitably the largest, which is also unique except on a set of measure zero, while some \(v_n\) must be the most misaligned to the corresponding direction \(v_m\). Consequently we may also write
     \begin{align*}
          \int_{\partial V_{(+)}(x)} \mathcal{A}_k[g](x,v,t) dv &=  \sum_{m = 1}^k \sum_{n = 0, n\neq m}^k \int_{(0,\infty)^k \times (\mathbb{S}^2)^{k+1}} J_{m,n}(v_k, \dots, v_0, s_k, \dots, s_1)\\
          &\phantom{=} \times g\Big(x_1 -\ell_{v_0}\left(x_1 \right)v_0,v_0,t - \sum_{l=1}^{k} s_{l}-\ell_{v_0}\left(x_1\right)\Big) \\ 
          &\phantom{=} \times F_k(x;v_{k}, \dots, v_0 ,s_k, \dots, s_1)  dv_{k} \dots dv_0 d s_k \dots d s_1
    \end{align*}
    which is nothing but the first identity that we used with \(J_{m,n}\) added. By performing the previous transform but choosing carefully the variables depending on the summation index we may arrive at
    \begin{align*}
          \int_{\partial V_{(+)}} \mathcal{A}_k[g](x,v,t) dv = &\int_{\partial D_{(-)}\times (0,\infty) } \sum_{m = 1}^k \sum_{n = 0, n\neq m}^k \int_{(0,\infty)^{k-1} \times (\mathbb{S}^2)^{k-1}} J_{m,n}(\argdot)\\ 
          &\phantom{\times} F_k(x;\argdot) \frac{g\left(x_0,v_0,t - s \right) \vert \Vec{n}(x_0) \cdot v_0 \vert}{s_m(\argdot)^2 (1-v_{n} \cdot v_m(\argdot))} \\
          & \phantom{\times} 1_{(0,s)}\left(\textstyle \sum_{l=0,l \neq m,n}^k s_{l} + \|(x-x_0 - \sum_{l=0,l \neq m,n}^k s_{l}v_l)\|_2\right)\\
          &\phantom{\times}  dv_{k} \dots dv_{m+1} dv_{m-1} \dots dv_{1} \\ 
          &\phantom{\times} d s_k \dots d s_{m+1}d s_{m-1}\dots d s_{n+1}d s_{n-1}\dots ds_0 ds dv_0 dx_0 
    \end{align*}
    which means that the Schwartz kernel \(\mathfrak{m}_k\) can be expressed as
    \begin{align}\label{eq:bounded_measurement_density}
    \begin{split}
        \mathfrak{m}_k(x,x_0,v_0,s) &= \sum_{m = 1}^k \sum_{n = 0, n\neq m}^k \int_{(0,\infty)^{k-1} \times (\mathbb{S}^2)^{k-1}} J_{m,n}(\argdot)\\ 
        & \phantom{\times} 1_{(0,s)}\left(\textstyle \sum_{l=0,l \neq m,n}^k s_{l} + \|(x-x_0 - \sum_{l=0,l \neq m,n}^k s_{l}v_l)\|_2\right)\\
        &\phantom{\times}  \frac{F_k(x;\argdot)\vert \Vec{n}(x_0) \cdot v_0 \vert}{s_m(\argdot)^2 (1-v_{n} \cdot v_m(\argdot))}dv_{k} \dots dv_{m+1} dv_{m-1} \dots dv_{1} \\ 
        &\phantom{\times} d s_k \dots d s_{m+1}d s_{m-1}\dots d s_{n+1}d s_{n-1}\dots ds_0
    \end{split}
    \end{align}
    Due to our previous developments we know that there is a constant \(C=C(\sigma_a, \sigma_s,f_p)\) that only depends on the optical parameters such that
    \begin{align}\label{eq:integrand_bound}
        C^k \left(\frac{k}{s - \|x-x_0\|_2}\right)^2 \frac{s^2}{s^2 - \|x-x_0\|_2^2}&\geq  \frac{F_k(x;\argdot)J_{m,n}(\argdot)}{s_m(\argdot)^2 (1-v_{n} \cdot v_m(\argdot))} 
    \end{align}
    for any pair \(m,n\). This means that for fixed \(k\) the integrand is bounded by a polynomial of degree \(4\) in \(k\) multiplied by a term that grows exponentially in \(k\). This shows that \(\mathfrak{m}_k\) is a function when restricted to the set \(M\). For continuity it is sufficient to show that each \(\mathfrak{m}_k\) is continuous on \(M\). In order to see the claim for the infinite sum \(\mathfrak{m}\) we may use formulas for the volume of unit balls, see e.g. \cite{Wang:2005}, to obtain the bound 
    \begin{align*}
        \int_{\{z \in (0,\infty)^{k} : \sum_{j=1}^k z_j \leq 1\}}dz_{k} \dots dz_1 \leq \int_{\{z \in \mathbb{R}^{k} : \sum_{j=1}^k \vert z_j \vert \leq 1\}}dz_{k} \dots dz_1 = \frac{2^k}{k!}.
    \end{align*}
    which means that the polynomial and exponential factors introduced in the coordinate change will not affect summability in any way and the claim follows by the dominated convergence theorem applied to the infinite sum. It remains to obtain continuity of \(\mathfrak{m}_k\) for each \(k \geq 2\). For that take any \((x^*,x^*_0,v^*_0,s^*_0) \in M\) and observe that due to \Cref{eq:integrand_bound} we can bound
    \begin{align*}
        \frac{F_k(x;\argdot)J_{m,n}(\argdot)}{s_m(\argdot)^2 (1-v_{n} \cdot v_m(\argdot))} = \mathcal{O}\left(\frac{C^k k^2 s^*_0}{(s^*_0 - \|x^*-x^*_0\|_2)^3}\right)
    \end{align*}
    uniformly, i.e. as a function of \(s^*_0\) and \(\|x^*-x^*_0\|_2\) only, in a sufficiently small \(\varepsilon\)-neighbourhood of \((x^*,x^*_0,v^*_0,s^*_0) \in M\). Each of the mappings
    \begin{align*}
        (x,x_0,v_0,s) &\mapsto J_{m,n}\\
        (x,x_0,v_0,s) &\mapsto 1_{(0,s)}\left(\textstyle \sum_{l=0,l \neq m,n}^k s_{l} + \|(x-x_0 - \sum_{l=0,l \neq m,n}^k s_{l}v_l)\|_2\right)\\
        (x,x_0,v_0,s) &\mapsto F_k(x;\argdot)
    \end{align*} 
    is \(dv_{k} \dots dv_{m+1} dv_{m-1} \dots dv_{1} d s_k \dots d s_{m+1}d s_{m-1}\dots d s_{n+1}d s_{n-1}\dots ds_0 \) almost everywhere continuous at \((x^*,x^*_0,v^*_0,s)\) because all indicator functions, which are the only sources of discontinuities, are continuous almost everywhere with respect to that measure. By bounded convergence we get that \(\mathfrak{m}_k\) must be continuous in \((x^*,x^*_0,v^*_0,s^*_0)\) and since the point was arbitrary the claim is true on \(M\).
    \end{proof}
\subsection{Proof of \texorpdfstring{\cref{lem:fixed_scatter_uniqueness}}{} \& \texorpdfstring{\cref{lem:fixed_phase_uniqueness}}{}}\label{proof:uniqueness}
    The following facts about the kernel spaces will be sueful. Due to the nature of our kernel functions, it must be true that for any \(p > 0\) we have
    \begin{align*}
        1-\sup \{z \in [0,1]: (1-z)^p \leq \phi(z)\} =: \varepsilon(p) > 0.
    \end{align*}
    Note that \(\varepsilon(p) \to 0\) as \(p \to \infty\) and by continuity also \(\phi(1-\varepsilon(p)) = \varepsilon(p)^p\). It follows for \(p\) sufficiently large that
    \begin{align}\label{eq:varepsilon_p_identity}
        \int_{0}^{\varepsilon(p)} \phi(1-z) dz \leq \int_{0}^{\varepsilon(p)} (1-z)^p dz = \frac{\varepsilon(p)^{p+1}}{p+1} = \frac{\varepsilon(p)}{p+1}  \phi(1-\varepsilon(p)).
    \end{align}
    Also note that
    \begin{align*}
       \frac{1}{h} \int_{0}^{\varepsilon} \phi\left(1-\frac{z}{h}\right) dz = \int_{0}^{\frac{\varepsilon}{h}} \phi\left(1-z\right) dz
    \end{align*}
    which, given that \(\phi\) is monotone on \((0,\infty)\), implies for any \(h_0 < h_1\) and any \(\varepsilon \in (0,h_0)\)
    \begin{align}\label{eq:kernel_integral_bandwidth_bound}
    \begin{split}
        \varepsilon \left(\frac{1}{h_0}-\frac{1}{h_1}\right) \phi\left(1-\frac{\varepsilon}{h_1}\right) &\leq \int_{\frac{\varepsilon}{h_1}}^{\frac{\varepsilon}{h_0}} \phi\left(1-z\right) dz\\
        &= \frac{1}{h_0} \int_{0}^{\varepsilon} \phi\left(1-\frac{z}{h_0}\right) dz - \frac{1}{h_1}\int_{0}^{\varepsilon} \phi\left(1-\frac{z}{h_1}\right) dz.
    \end{split}
    \end{align}
\begin{proof}[Proof of \cref{lem:fixed_scatter_uniqueness}]\label{proof:fixed_scatter_uniqueness}
    The direction \(\Rightarrow\) is trivial so we only have to show \(\Leftarrow\), i.e. that the measurement is different for different differential absorption functions. We may assume without loss of generality that \(x_D = 0\), \(\phi(0)=1\) and, because \(\phi\) is monotone, \(\phi(z) > 0\) for any \(z \in [0,1)\). Let \(\alpha_1 \neq \alpha_2\) be two different perturbations of \(\sigma_{a}\) as in the lemma. Further let \(\tau_{1/2}:= \mathrm{dist}(0, \supp(\alpha_{1/2}))\) and \(v_{1/2}\) directions such that \(\tau_{1/2} v_{1/2} \in \partial \mathrm{supp} (\alpha_{1/2})\). If \(\tau_{1} > \tau_2\) or \(\tau_{1} < \tau_2\) then obviously the measurements will differ at times between \(2 \tau_1\) and \(2 \tau_2\) in direction \(v_1\) or \(v_2\) as the one corresponding to the smaller first impact time will be influenced by an absorbing perturbation and differ from the offline measurement while the other won't. Therefore we assume \(\tau_{1} = \tau_2\) and define \(\alpha_\delta:=\alpha_1-\alpha_2\) as well as the first impact time \(\tau_{\delta} = \mathrm{dist}(0, \supp(\alpha_\delta))\). 

    Note that \(\tau_{\delta} \geq \tau_1 = \tau_2\) may have increased and \(\tau_\delta > 0\) due to \(\alpha_{1/2} \in K(\phi \mid x_D)\). There might be multiple possibilities for the first impact direction, i.e. the direction \(v\) such that \(v \tau_{\delta} \in \partial \supp (\alpha_\delta)\). The constant term \(w_0\) in \(\alpha_i\) can be identified from any measurement at time \(s < \tau_i\) for \(i=1,2\). Thus, if the constants are different there is nothing more to show. In case that the constant terms are equal we can write
    \begin{align*}
        \alpha_{\delta}(x) = \sum_{j=1}^N w_j \phi\left(\frac{\|x - m_j\|}{h_j}\right)
    \end{align*}
    with \(w_j \in \mathbb{R}\setminus \{0\}\) positive and negative such that \((m_k,h_k) \neq (m_j,h_j)\) for all \(j \neq k\) which makes the representation unique. Combining \Cref{eq:varepsilon_p_identity} and \Cref{eq:kernel_integral_bandwidth_bound} we can see that it is not possible to replicate the tail of a kernel will tails of larger width which implies that every possibility for \(v\) which satisfies \(v \tau_{\delta} \in \partial \supp (\alpha_\delta)\) must point towards a kernel centre \(m_j\). We choose a direction \(v_{\delta} \in \mathbb{S}^2\) with that property such that it points towards the kernel centre \(m_j\) with the smallest possible bandwidth parameter \(h_j\). If this isn't unique, any of them is good. Let \(h_{\min}, m_{\min}\) and \(w_{\min}\) be the kernel parameters that correspond to the kernel that was used in the choice of \(v_\delta\) and without loss of generality assume that the corresponding weight is positive, i.e. \(w_{\min} > 0\). We can define
    \begin{align*}
        A(\varepsilon)&=\int_{0}^{\tau_{\delta}+\varepsilon} \alpha_\delta(sv_\delta) ds
    \end{align*}
    and now claim that for \(\mathfrak{m}_{\mathrm{on}}^{(1/2)}\) corresponding to \(\alpha_{1/2}\) we have
    \begin{align}\label{eq:uniqueness_limit}
        \limsup_{\varepsilon \searrow 0} \left \lvert\frac{\mathfrak{m}_{\mathrm{on}}^{(1)}(v_\delta, 2(\tau_\delta + \varepsilon) \mid x_D)-\mathfrak{m}_{\mathrm{on}}^{(2)}(v_\delta, 2(\tau_\delta + \varepsilon) \mid x_D)}{A(\varepsilon)} \right \rvert > 0
    \end{align}
    which would prove our uniqueness statement. Let \(h_{\mathrm{2nd}}:= \inf\{h_j : h_j > h_{\min}\}\) be the smallest width parameter larger than \(h_{\min}\). From the above we get for \(\varepsilon\) sufficiently small,
    \begin{align*}
        A(\varepsilon)-w_{\min}\int_{0}^{\varepsilon} &\phi\left(1-\frac{z}{h_{\min}}\right) dz \\
        &\geq - \int_{0}^{\tau_{\delta}+\varepsilon} \left| \alpha_\delta(sv_\delta)  - w_{\min} \phi\left(1-\frac{\|sv_\delta - m_{\min}\|}{h_{\min}}\right) \right|ds.\\
        &\geq -\sum_{j=1}^N |w_{j}|\int_{0}^{\varepsilon} \phi\left(1-\frac{z}{h_{\mathrm{2nd}}}\right) dz 
    \end{align*}
    where the second inequality holds because for positive weights the integral is maximised when the integration is done in the direction that points from \(x_D\) towards \(m_j\) and decreases with the kernel width. Consequently,
    \begin{align*}
        \sum_{j=1}^N |w_{j}|\int_{0}^{\varepsilon} \phi\left(1-\frac{z}{h_{\mathrm{2nd}}}\right) dz \geq \left| A(\varepsilon)-w_{\min}\int_{0}^{\varepsilon} \phi\left(1-\frac{z}{h_{\min}}\right) dz \right|
    \end{align*}
    For \(\varepsilon\) sufficiently small we have \(\|m_j - x_D - (\varepsilon + \tau_{\delta})v_{\delta} \| < h_j\) implies that \(m_{\min} = m_j\) and \(h_{\min}=h_j\) or \(h_{\min} < h_j\) which means that only larger widths than \(h_{\min}\) must be taken into account. Taking \(\tilde{\varepsilon}(p) = h_{\mathrm{2nd}} \varepsilon(p)\) we see that as \(p \to \infty\) 
    \begin{align*}
        \frac{\int_{0}^{\tilde{\varepsilon}(p)} \phi\left(1-\frac{z}{h_{\mathrm{2nd}}}\right) dz }{\int_{0}^{\tilde{\varepsilon}(p)} \phi\left(1-\frac{z}{h_{\min}}\right) dz } \leq \frac{\int_{0}^{\varepsilon(p)} \phi\left(1-z\right) h_{\mathrm{2nd}} dz }{\varepsilon(p) (1-h_{\min}/h_{\mathrm{2nd}})\phi\left(1-\varepsilon(p)\right)} \leq \frac{h^2_{\mathrm{2nd}}}{(h_{\mathrm{2nd}}-h_{\min})(p+1)}    
    \end{align*}
    the first summand in \(A(\varepsilon)\) dominates, i.e.
    \begin{align}\label{eq:kernel_bound_direct}
        \lim_{p \to \infty} \frac{A(\tilde{\varepsilon}(p))}{w_{\min}\int_{0}^{\tilde{\varepsilon}(p)} \phi\left(1-\frac{z}{h_{\min}}\right) dz } = 1        
    \end{align}
    Let \(\Gamma(\varepsilon)\) be the set of all continuous piece-wise linear curves that start and end in \(x_D\) with length \(2(\tau_\delta + \varepsilon)\). Pick \(\gamma \in \Gamma(\varepsilon)\) and assume that \(\gamma:[0,2(\tau_\delta + \varepsilon)] \to \mathbb{R}^3\) satisfies \(\gamma(0)=\gamma(2(\tau_\delta + \varepsilon)) = x_D\) and \(
    \|\gamma'(s)\|=1\) almost everywhere. Then we have, again for \(\varepsilon\) sufficiently small,
    \begin{align}\label{eq:kernel_bound_curve}
        \begin{split}
        \left \lvert \int_{0}^{2(\tau_\delta + \varepsilon)} \alpha_\delta(\gamma(s)) ds \right \rvert &\leq 2 \sum_{j:\|m_j - x_D - (\varepsilon + \tau_{\delta})v_{\delta} \| < h_j} |w_{j}|\int_{0}^{\varepsilon} \phi\left(1-\frac{z}{h_{j}}\right) dz \\
        &\leq 2 \sum_{j=1}^{N} |w_{j}|\int_{0}^{\varepsilon} \phi\left(1-\frac{z}{h_{\min}}\right) dz.
        \end{split}
    \end{align}
    We next analyse the decay of first and higher order scattering when scaled by \(A(\varepsilon)\). From \cref{lem:path_continuity} we know that \(\mathfrak{m}_{k}\) is continuous in a neighbourhood of \((x_D,x_D,v_{\delta},\tau_{\delta})\). Using the same notation, i.e. \(x_l\) for \(l=0,\dots,k+1\) as in \cref{lem:path_continuity}, we know that
    \begin{align*}
        \mathfrak{m}_{\mathrm{off},k}(v,s \mid x_D) =  \sum_{m = 1}^k \sum_{n = 0, n\neq m}^k &\int_{(0,\infty)^{k-1} \times (\mathbb{S}^2)^{k-1}} R^{\mathrm{off}}_{m,n}(v,s; \argdot) d \nu_{m,n}
    \end{align*}
    for \(d \nu_{m,n} = dv_{k} \dots dv_{m+1} dv_{m-1} \dots dv_{1} d s_k \dots d s_{m+1}d s_{m-1}\dots d s_{n+1}d s_{n-1}\dots ds_0\) and functions \(R^{\mathrm{off}}_{m,n}\) given by \Cref{eq:bounded_measurement_density} that are bounded via \Cref{eq:integrand_bound} uniformly in a neighbourhood of \((x_D,x_D,v_\delta,2\tau_\delta)\). Note that the dependence of \(R\) w.r.t. \(k\) has been omitted to simplify notation. If we put
    \begin{align*}
        R^{\mathrm{on},i}_{m,n} = \exp\Big(- \int_{\Gamma(x_{0:k+1})} \alpha_i(r)dr\Big) R^{\mathrm{off}}_{m,n}
    \end{align*}
    for \(i = 1,2\) and set \(\Delta\mathfrak{m}_{\mathrm{on},k} = \mathfrak{m}_{\mathrm{on},k}^{(1)}-\mathfrak{m}_{\mathrm{on},k}^{(2)}\) to be the difference in the online component of the \(k\)-th scattering contribution for all \(k \geq 2\), then it is not hard to see that
    \begin{align*}
         \lim_{p \to \infty} \left \lvert\frac{\Delta\mathfrak{m}_{\mathrm{on},1}(v_\delta, 2(\tau_\delta + \tilde{\varepsilon}(p)) \mid x_D)}{2 A(\tilde{\varepsilon}(p))} \right \rvert = \mathfrak{m}_{\mathrm{off},1}(v_\delta, 2\tau_\delta \mid x_D).
    \end{align*}
    We claim that at the same distance \(\tau_\delta\) any higher order contributions \(\Delta \mathfrak{m}_{\mathrm{on},k}\) are much less affected by the absorption and we have for any \(k \geq 2\)
    \begin{align*}
         \lim_{p \to \infty} \left \lvert\frac{\Delta \mathfrak{m}_{\mathrm{on},k}(v_\delta, 2(\tau_\delta + \tilde{\varepsilon}(p)) \mid x_D)}{2 A(\tilde{\varepsilon}(p))} \right \rvert = 0.
    \end{align*}
    Since \(\mathfrak{m}_{\mathrm{off},1}(v_\delta, 2\tau_\delta \mid x_D) > 0\) we would obtain \Cref{eq:uniqueness_limit} by bounded convergence from the decay rate of \(\mathfrak{m}_{\mathrm{off},k}\) which was shown in the proof of \cref{lem:path_continuity} towards the end. To see this it is enough to show that 
    \begin{align}\label{eq:higher_order_scattering_limit}
        \frac{R^{\mathrm{on},1}_{m,n}(v_\delta, 2(\tau_\delta + \tilde{\varepsilon}(p)), \argdot)-R^{\mathrm{on},2}_{m,n}(v_\delta, 2(\tau_\delta + \tilde{\varepsilon}(p)); \argdot)}{A(\tilde{\varepsilon}(p))} \xrightarrow{p \to \infty} 0
    \end{align}
    in \(L^1((0,\infty)^{k-1} \times (\mathbb{S}^2)^{k-1}, d\nu_{m,n})\). Thanks to \Cref{eq:kernel_bound_curve} and \Cref{eq:kernel_bound_direct} we know that the left-hand side of \Cref{eq:higher_order_scattering_limit} is bounded by \(C R^{\mathrm{off}}_{m,n}(v_\delta, 2\tau_\delta)\) for some \(p\)-independent constant \(C>0\) that may depend on \(\tau_\delta\). Therefore it is enough to show the limit in \Cref{eq:higher_order_scattering_limit} \(\nu_{m,n}\) almost everywhere. We may also notice that the left-hand side of \Cref{eq:higher_order_scattering_limit} is non-zero only when 
    \(1_{(\tau_\delta, \infty)}(\|x_l(v_\delta, 2(\tau_\delta + \tilde{\varepsilon}(p); \argdot)\|_2)\) is non-zero for some \(l = 1, \dots, k\), i.e. when \(\Gamma(x_{0:k+1})\) leaves the ball with radius \(\tau_\delta\) around the light source. For \(k \geq 3\) there is an index \(q=1,\dots,k\) such that \(q\neq m\) and \(q \neq n\). If \(n\neq 0\) and \(s_0 s_q (1-\lvert v_\delta\cdot v_q\rvert) > z\) for some \(p\)-independent \(z > 0\) then \(1_{(\tau_\delta, \infty)}(\|x_l(v_\delta, 2(\tau_\delta + \tilde{\varepsilon}(p); \argdot)\|_2)=0\) for \(p\) sufficiently large whereas if \(n=0\) then there are two indices \(q_1 \neq q_2\) that differ from \(m\) and \(n\) and we have the same result for \(s_{q_1} s_{q_2} (1-\lvert v_{q_1}\cdot v_{q_2}\rvert) > z\). Note that both of these conditions are necessary for \(\Gamma(x_{0:k+1})\) to be sufficiently straight so that it can leave the ball with radius \(\tau_\delta\) around the light source. This shows \Cref{eq:higher_order_scattering_limit} for \(k \geq 3\) because, e.g. in the case of \(n=0\), we have
    \begin{align*}
        1_{(0,z)}(s_{q_1} s_{q_2} (1-\lvert v_{q_1}\cdot v_{q_2}\rvert)) \xrightarrow{z \to 0} 0
    \end{align*}
    \(\nu_{m,0}\) almost everywhere. Similar conditions can be found when \(k=2\). Given that \Cref{eq:higher_order_scattering_limit} for all \(k \geq 2\) implies \Cref{eq:uniqueness_limit} this concludes the proof.
    \end{proof}
    The proof of \cref{lem:fixed_phase_uniqueness} is similar. Note that it is enough to show the result for \(\sigma_a = 0\) since we have
    \begin{align*}
        \frac{\int_{0}^{\varepsilon}\phi(1-z) dz}{\phi(1-\varepsilon)} \to 0 \qquad \mathrm{for} \qquad \varepsilon \to 0
    \end{align*}
    and scattering locally at the tails dominates non-differential absorption with respect to its contribution in \(\mathfrak{m}_{\mathrm{off}}\). Note that the proportionality requirement, which implies that \(\sigma_a\) is uniquely determined by \(\sigma_s\), can be dropped under certain conditions on the kernel functions. Similar, although more restrictive, assumptions are made in \cref{thm:relaxed_uniqueness} whose proof is subject of the next section. 
    \begin{proof}[Proof of \cref{lem:fixed_phase_uniqueness} - sketch for \(\sigma_a = 0\)]\label{proof:fixed_phase_uniqueness}
    Using the same notation as we did for \cref{lem:fixed_scatter_uniqueness} we have in the same way as in \cref{eq:uniqueness_limit} that
    \begin{align}
        \limsup_{\varepsilon \searrow 0} \left \lvert\frac{\mathfrak{m}_{\mathrm{off}}^{(1)}(v_\delta, 2(\tau_\delta + \varepsilon) \mid x_D)-\mathfrak{m}_{\mathrm{off}}^{(2)}(v_\delta, 2(\tau_\delta + \varepsilon) \mid x_D)}{S(\varepsilon)} \right \rvert > 0
    \end{align}
    where \(\mathfrak{m}_{\mathrm{off}}^{(1/2)}\) correspond to different values \(\sigma^{(1/2)}_s\) and \[S(\varepsilon) = \sigma^{(1)}_s(x_D + v_\delta(\tau_\delta + \varepsilon)) - \sigma^{(2)}_s(x_D + v_\delta(\tau_\delta + \varepsilon)).\] As before the above follows from the fact that 
    \begin{align}
        \limsup_{\varepsilon \searrow 0} \frac{S(\varepsilon)}{w_{\min}\phi\left(1-\frac{\varepsilon}{h_{\min}}\right) } > 0        
    \end{align}
    which can be shown in virtually the same fashion as in the previous proof.
    \end{proof}
    
    \subsection{Semi-parametric relaxation and proof of \texorpdfstring{\cref{thm:relaxed_uniqueness}}{}}
    The proof of \cref{thm:relaxed_uniqueness} is somewhat more involved as it only assumes access to the quotient of our data which is equivalent to injectivity in a more general forward model of the form \((\tilde{\rho}_{\mathrm{on}}(v,t),\tilde{\rho}_{\mathrm{off}}(v,t))\) as per \cref{eq:off_constraint,eq:on_constraint} where \cref{eq:psi_constraint} has been dropped. Despite the relaxation, much of the previously developed ideas translate to the proof. 
    \begin{proof}\label{prf:relaxed_uniqueness}
        As before we may assume without loss of generality that \(x_D = 0\) as well as \(\phi_{\mathfrak{s}}(0)=\phi_{\mathfrak{a}}(0)=1\) and \(\phi_{\mathfrak{s}}(z), \phi_{\mathfrak{a}}(z)> 0\) for any \(z \in [0,1)\). Let \(\alpha_1 \neq \alpha_2\) be two different perturbations of \(\sigma_{a}\) as in the proof of \cref{lem:fixed_scatter_uniqueness} and define \(\tau_{1/2}, v_{1/2}\) as well as \(\alpha_\delta\) and \(\tau_{\delta} = \mathrm{dist}(0, \supp(\alpha_\delta))\) in the same way as before. If \(\alpha_\delta = 0\) we set \(\tau_\delta = \infty\). Note that If \(\tau_{1} > \tau_2\) or \(\tau_{1} < \tau_2\) then obviously the measurements will differ at times between \(2 \tau_1\) and \(2 \tau_2\) in direction \(v_1\) or \(v_2\) as the one corresponding to the smaller first impact time will be influenced by an absorbing perturbation and differ from the offline measurement while the other won't. Therefore we assume without loss of generality that \(\tau_{1} = \tau_2\). Note that \(\tau_{\delta} \geq \tau_1 = \tau_2\) may have increased and \(\tau_\delta > 0\) due to \(\alpha_{1/2} \in K(\phi_{\mathfrak{a}} \mid x_D)\). There might be multiple possibilities for the first impact direction, i.e. the direction \(v\) such that \(v \tau_{\delta} \in \partial \supp (\alpha_\delta)\) and we want to choose a direction \(v_{\delta} \in \mathbb{S}^2\) with that property such that it points towards the kernel centre \(m_{\min}\) with the smallest possible bandwidth parameter \(h_{\min}\) and without loss of generality assume that the corresponding weight is positive, i.e. \(w_{\min} > 0\), and corresponds to a summand in \(\alpha_1\). 
        
        Case 1: \(\tau_1 = \tau_2 = \tau_\delta < \infty\). Define 
        \begin{align}
            P^{(1/2)}_1(\varepsilon) &= \frac{\mathfrak{m}_{\mathrm{on},1}^{(1/2)}(v_\delta, 2(\tau_\delta + \varepsilon) \mid x_D)}{\mathfrak{m}_{\mathrm{off}}^{(1/2)}(v_\delta, 2(\tau_\delta + \varepsilon) \mid x_D)} \\
            P^{(1/2)}_2(\varepsilon) &= \frac{\sum_{j=2}^\infty \mathfrak{m}_{\mathrm{on},j}^{(1/2)}(v_\delta, 2(\tau_\delta + \varepsilon) \mid x_D)}{\mathfrak{m}_{\mathrm{off}}^{(1/2)}(v_\delta, 2(\tau_\delta + \varepsilon) \mid x_D)}
        \end{align}
        and note that the differential absorption, i.e. our data, can be expressed as
        \begin{align*}
            P^{(1/2)}_1 + P^{(1/2)}_2 = \frac{\mathfrak{m}_{\mathrm{on}}^{(1/2)}(v_\delta, 2(\tau_\delta + \varepsilon) \mid x_D)}{\mathfrak{m}_{\mathrm{off}}^{(1/2)}(v_\delta, 2(\tau_\delta + \varepsilon) \mid x_D)} .
        \end{align*}
        A similar argument that was used to show \cref{eq:uniqueness_limit} can be used to show that 
        \[
        P^{(1/2)}_1(\varepsilon)\to \frac{\mathfrak{m}_{\mathrm{off},1}^{(1)}(v_\delta, 2\tau_\delta\mid x_D)}{\mathfrak{m}_{\mathrm{off}}^{(1)}(v_\delta, 2\tau_\delta \mid x_D)}=\frac{\mathfrak{m}_{\mathrm{off},1}^{(2)}(v_\delta, 2\tau_\delta\mid x_D)}{\mathfrak{m}_{\mathrm{off}}^{(2)}(v_\delta, 2\tau_\delta \mid x_D)} \qquad \mathrm{as}\qquad \varepsilon \to 0
        \] 
        slower than 
        \[
        P^{(1/2)}_2(\varepsilon)\to 1-\frac{\mathfrak{m}_{\mathrm{off},1}^{(1/2)}(v_\delta, 2\tau_\delta\mid x_D)}{\mathfrak{m}_{\mathrm{off}}^{(1/2)}(v_\delta, 2\tau_\delta \mid x_D)} \qquad \mathrm{as}\qquad \varepsilon \to 0.
        \]
        Further, analogously to the proof of \cref{lem:fixed_scatter_uniqueness}, we have that \(P^{(1)}_1(\varepsilon)\) converges at a rate of \(\int_0^\varepsilon \phi_\mathfrak{a}(1-\frac{z}{h_{\min}})dz\) slower than \(P^{(2)}_1(\varepsilon)\) (due to \(w_{\min}\) corresponding to \(\alpha_1\)). As such
        \[
        \limsup_{\varepsilon \to 0} \left|\frac{P^{(1)}(\varepsilon)-P^{(2)}(\varepsilon)}{\int_0^\varepsilon \phi_\mathfrak{a}(1-\frac{z}{h_{\min}})dz} \right| > 0
        \]
        which implies \(P_1 \neq P_2\) as desired. Note that scattering from the kernels is locally negligible for the limit due to \Cref{eq:kernel_tails}.
        
        Case 2: \(\tau_1 = \tau_2 < \tau_\delta < \infty\) but identical scattering proportionality constants. Essentially identical to case 1 where  \(\tau_1 = \tau_2 = \tau_\delta < \infty\) except that 
        \[
        P^{(1/2)}_1(\varepsilon) \to \frac{\mathfrak{m}_{\mathrm{on},1}^{1}(v_\delta, 2\tau_\delta\mid x_D)}{\mathfrak{m}_{\mathrm{off}}^{1}(v_\delta, 2\tau_\delta \mid x_D)}=\frac{\mathfrak{m}_{\mathrm{on},1}^{2}(v_\delta, 2\tau_\delta\mid x_D)}{\mathfrak{m}_{\mathrm{off}}^{2}(v_\delta, 2\tau_\delta \mid x_D)}
        \]
        and again we have that \(P^{(1)}_1(\varepsilon)\) converges at a rate of \(\int_0^\varepsilon \phi_\mathfrak{a}(1-\frac{z}{h_{\min}})dz\) slower than \(P^{(2)}_1(\varepsilon)\). 

        Case 3: Differing proportionality constants. Note that \(\tau_\delta > \tau_1 = \tau_2\) or we can fall back to case 1. First assume that the proportionality between \(\sigma_s\) and \(\alpha\) is different. Let \(\bar{\tau}:=\tau_1=\tau_2\) and note that this is the smallest time where a difference in scattering can be observed but the differential absorption fields \(\alpha_{1/2}\) are identical until \(\tau_\delta > \bar{\tau}\). In a similar way to \(v_\delta\) let \(\bar{v}\) be the direction pointing to the smallest kernel at a distance of \(\bar{\tau}\) from \(x_D\). As before denote this width with \(h_{\min}\) which is necessarily (i.e. by assumption of the theorem) identical for both sets of non-differential parameters. Consider \(\tilde{\mathfrak{m}}_{\mathrm{on/off}}\) which matches both \(\mathfrak{m}_{\mathrm{on}}^{1}\) and \(\mathfrak{m}_{\mathrm{on}}^{2}\) in the ambient constants as well as the differential absorption for times \(t < 2\bar{\tau} + \varepsilon\) for \(\varepsilon\) sufficiently small but no scattering from the kernels (in other words non-differential proportionality constants identical to \(0\)). By construction it must be true that 
        \[
        \frac{\mathfrak{m}_{\mathrm{on}}^{(1/2)}(\bar{v}, 2(\bar{\tau} + \varepsilon) \mid x_D)}{\mathfrak{m}_{\mathrm{off}}^{(1/2)}(\bar{v}, 2(\bar{\tau} + \varepsilon) \mid x_D)} - \frac{\tilde{\mathfrak{m}}_{\mathrm{on}}(\bar{v}, 2(\bar{\tau} + \varepsilon) \mid x_D)}{\tilde{\mathfrak{m}}_{\mathrm{off}}(\bar{v}, 2(\bar{\tau} + \varepsilon) \mid x_D)} \to 0 \qquad \mathrm{as}\qquad \varepsilon \to 0.
        \] 
        Upon closer inspection it is evident that 
        \[
        \frac{\frac{\mathfrak{m}_{\mathrm{on}}^{(1/2)}(\bar{v}, 2(\bar{\tau} + \varepsilon) \mid x_D)}{\mathfrak{m}_{\mathrm{off}}^{(1/2)}(\bar{v}, 2(\bar{\tau} + \varepsilon) \mid x_D)} - \frac{\tilde{\mathfrak{m}}_{\mathrm{on}}(\bar{v}, 2(\bar{\tau} + \varepsilon) \mid x_D)}{\tilde{\mathfrak{m}}_{\mathrm{off}}(\bar{v}, 2(\bar{\tau} + \varepsilon) \mid x_D)}}{\phi_{\mathfrak{s}}(1-\frac{\varepsilon}{h_{\min}}) \int_0^\varepsilon \phi_\mathfrak{a}(1-\frac{z}{h_{\min}})dz }
        \]
        has a non-trivial limit which uniquely determines the proportionality constant. Similarly, if the proportionality only differs in the \(\sigma_a\) component we additionally match \(\tilde{\mathfrak{m}}_{\mathrm{on/off}}\) in the \(\sigma_s\) related constant and consider \[\frac{\frac{\mathfrak{m}_{\mathrm{on}}^{(1/2)}(v_\delta, 2(\tau_\delta + \varepsilon) \mid x_D)}{\mathfrak{m}_{\mathrm{off}}^{(1/2)}(v_\delta, 2(\tau_\delta + \varepsilon) \mid x_D)} - \frac{\tilde{\mathfrak{m}}_{\mathrm{on}}(v_\delta, 2(\tau_\delta + \varepsilon) \mid x_D)}{\tilde{\mathfrak{m}}_{\mathrm{off}}(v_\delta, 2(\tau_\delta + \varepsilon) \mid x_D)}}{\int_0^\varepsilon \phi_\mathfrak{s}(1-\frac{z}{h_{\min}})dz \int_0^\varepsilon \phi_\mathfrak{a}(1-\frac{z}{h_{\min}})dz }\]\
        which, as this covers all possible cases, concludes the proof.
    \end{proof}
    \section{Implementation \& simulation details}\label{sec:implementation_details}
    As mentioned in \cref{sec:experiments} we have not only considered smooth plumes that satisfy \cref{eq:general_plume_equation} for a set of highly regular parameters but also conducted simulations where the ground truth was perturbed in a way akin to turbulence. The low-dimensional dispersion model, denoted henceforth by \(\alpha_0\), has 14 parameters. A scalar release rate, a ground based source term (\(x\) and \(y\) coordinate) as well as 10 parameters that model the width and height of the centre-line by means of a linear spline (with 5 components each) and a proportionality constant for \(\sigma_s\) as described in \cref{thm:relaxed_uniqueness} (and we assume \(\sigma_a=0\)). Note that the centre-line of the plume is uniquely determined by the source, the (known) wind velocity and its \(z\) component (upwards drift). The kernel functions are Gaussian which are numerically very similar to good kernels as in \cref{def:kernel_space}.

    A random perturbation of the average model can be obtained in many ways although we should mention that this would be non-trivial if a standard model for turbulence was used as the non-linearities would mean that matching the parameters isn't straight forward. This is, together with the computational burden that a simulation of turbulence paired with an RTE solution would bear, a major factor for the following simplified modelling approach. If we describe the perturbed gas as
    \begin{align*}
        \hat{\alpha} = \frac{\sum_{k=1}^{K}\sum_{l=1}^{M_k} \hat{\alpha}_{k,l}}{\sum_{k=1}^{K} M_k} 
    \end{align*}
    then \(\mathbb{E}(\hat{\alpha}) = \alpha_0\) follows from \(\mathbb{E}(\hat{\alpha}_{k,l}) = \alpha_0\) for all \(k,l\). We assume that for fixed \(k\) the \(\hat\alpha_{k,l}\) are iid. The condition \(\mathbb{E}(\hat{\alpha}_{k,l}) = \alpha_0\) can be achieved by adding an appropriately scaled mean zero path with Gaussian marginal distributions to the plume's centre-line and adjusting the kernel widths accordingly. Indeed, if the marginal standard deviation of the centre-line perturbation is a fraction \(\sqrt{\phi_{V,k}} \in (0,1)\) of the plume width then we can set the adjusting kernel width to \(\sqrt{1-\phi_{V,k}}\) in order to achieve the desired effect. Similarly we may introduce (randomly occurring) jumps as samples from the adjusted kernel scaled by a fraction \(\sqrt{\phi_{J,k}} \in (0,1)\) and (recursively after each jump) reduce the kernel width by a fraction of \(\sqrt{1-\phi_{J,k}}\) and thereby introduce smaller eddys. By redistributing the kernel weights we may equally introduce a splitting process at every jump after which two independent trajectories emerge from a the previous location (each offset by an independently sampled jump). 

    For the purpose of our simulation we picked \(K=5\) and an AR(1) process with auto-covariance parameter \(0.9\) as the centre line perturbation in each case. Furthermore we chose \(M=(1,1,3,2,1)\) and \(\phi_{V}=(0, 8/9, 2/3, 2/3, 4/9)\) with no jumps in \(k=1,2,3\) while we have \(2\) jumps for \(k=4\) and \(4\) for \(k=5\). The jumps are uniformly distributed over the plume's centre line, i.e. w.r.t. time/down-wind distance, the jump widths are \(\phi_{J}=(0, 0,0,0.55,0.4)\) and the path splits in two at each jump with weights allocated uniformly at random.
    \begin{figure}[H]
        \centering
        \includegraphics[scale=0.45,trim={3cm 0.25cm 2.75cm 1.3cm},clip]{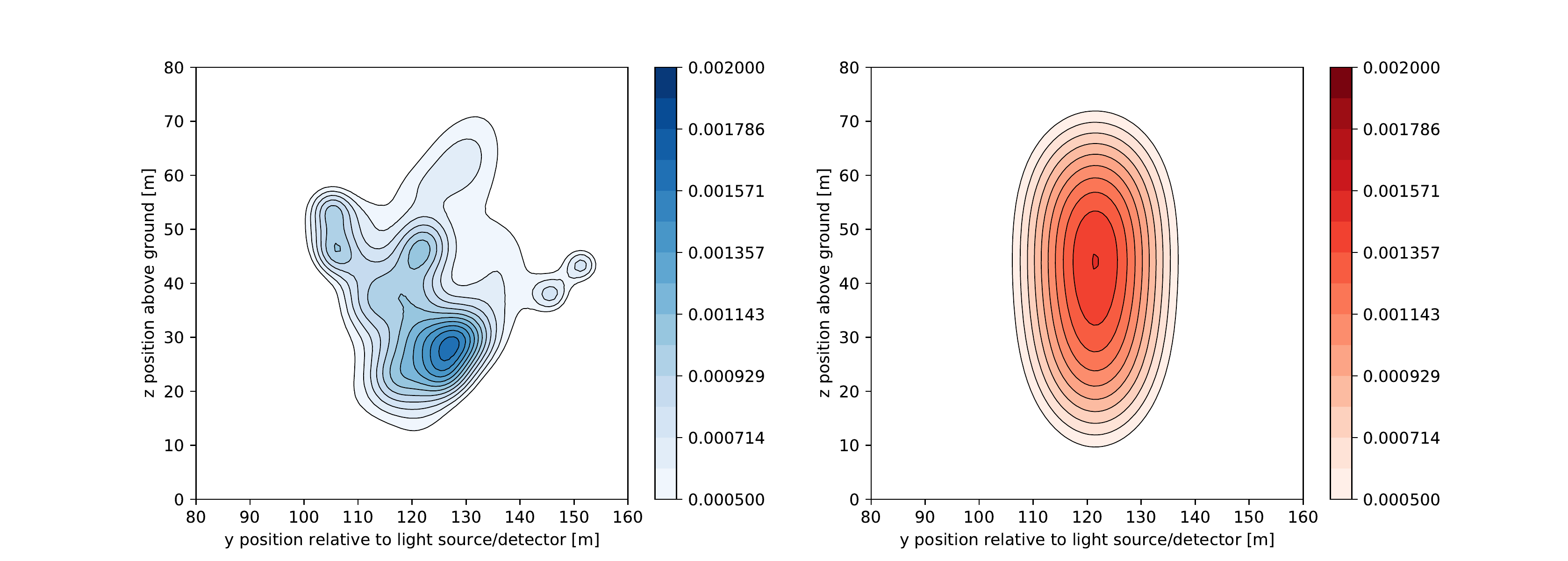}
        \caption{Cross sections perpendicular to the wind direction of the gas distribution approximately 60m away from the source. Note that the oval shape is due to the upwards drift of the gas. The \(L_1\) difference between the smooth and turbulent plume is approximately \(50 \%\).}
    \end{figure}
    It should be noted that the plume is assumed finite, roughly 120m in length, because at greater distances the increased width results in concentrations that don't allow any degrees of freedom. As such our reconstructions would be dominated by the regularisation term \(\mathsf{R}\) which in the presented case applies a finite difference operator to both the width and upwards drift parameters. Motivated by our findings in \cref{fig:free_test} and \cref{fig:known_test} we selected a quadratic penalty term that restricts the scattering parameter to values that correspond to levels similar as those discussed in \cref{sec:stability}. Despite the regularisation we observe that the reconstructed parameter varies more than its absorption counterpart which is not regularised. This is most likely due to the low sensitivity of the data to the scattering as suggested in \cref{sec:stability} as well as the proof of \cref{thm:relaxed_uniqueness}.
    \begin{figure}[H]
        \centering
        \includegraphics[scale=0.05,trim={5cm 35cm 5cm 40cm},clip]{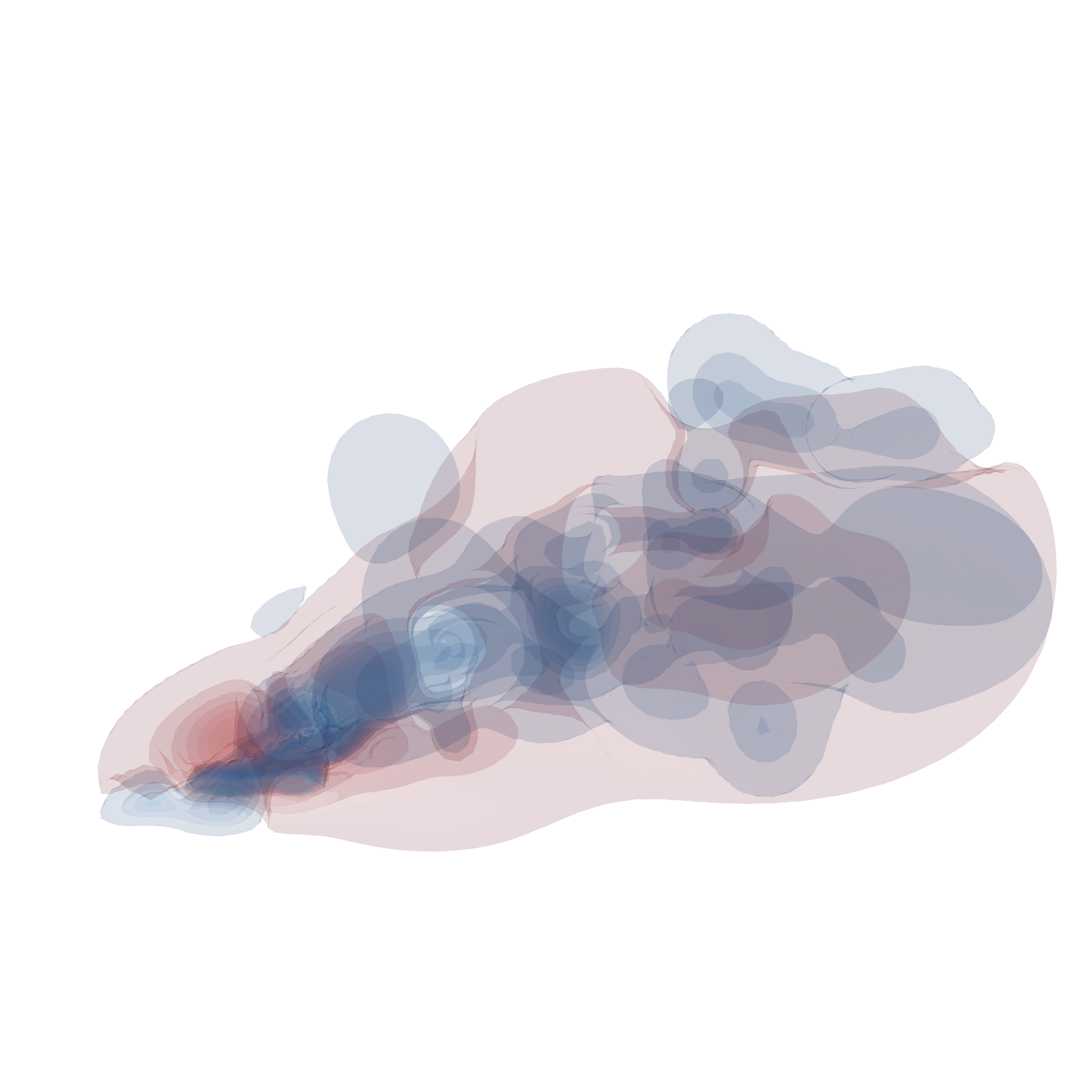}
        \caption{Qualitative structure of the turbulence: In the blue regions the turbulent gas is larger than its smooth counterpart, red indicates the opposite. The \(L_1\) difference between the smooth and turbulent plume is approximately \(50 \%\).}
    \end{figure}

\newpage
\bibliographystyle{siamplain}
\bibliography{references}
\end{document}